\RequirePackage{fix-cm}
\documentclass[smallextended]{svjour3} 
\smartqed  
\usepackage{graphicx,epsfig}
\usepackage{enumerate}
\usepackage{amssymb}
\usepackage{amsmath}
\usepackage{centernot}
\usepackage{xcolor}
\usepackage{todonotes}
\usepackage{times}
\usepackage{mathrsfs}
\usepackage[utf8]{inputenc}
\usepackage{multirow}

\usepackage{pifont}
\usepackage{subfig}

\usepackage{mathtools}
\usepackage{appendix}
\usepackage{cite}
\usepackage{wrapfig}
\usepackage{setspace}
\usepackage{booktabs}

\usepackage{lscape}
\usepackage{algorithm,algorithmicx}
\usepackage{algcompatible}
\usepackage{bm}
\usepackage{caption}
\usepackage{xfrac}
\usepackage{url}
\usepackage{longtable}
\usepackage{appendix}

\newtheorem{assumption}{Assumption}
   
 \newtheorem{thm}[theorem]{Theorem}{\bfseries}{\itshape}   
\spnewtheorem{lem}[theorem]{Lemma}{\bfseries}{\itshape} 
 \newtheorem{rem}[theorem]{Remark}{\bfseries}{\itshape} 
  \newtheorem{prop}[theorem]{Proposition}{\bfseries}{\itshape}  
\counterwithin{theorem}{section}

\numberwithin{Subcase}{Subcase}
\usepackage[linkcolor=blue, urlcolor=blue, citecolor=blue,
colorlinks, bookmarks]{hyperref}
\journalname{}
\begin{document}
\title{A Three-Operator Splitting Scheme Derived from Three-Block ADMM}

\titlerunning{}

\author{Anshika Anshika \and Jiaxing Li \and Debdas Ghosh \and Xiangxiong Zhang$^*$\thanks{* Corresponding author: Xiangxiong Zhang (zhan1966@purdue.edu)} }


\institute{Anshika Anshika (anshika.rs.mat19@itbhu.ac.in) \at
Department of Mathematical Sciences, Indian Institute of Technology (Banaras Hindu University) Varanasi, Uttar Pradesh  221005, India
\and
Jiaxing Li (li4944@purdue.edu) \at
Department of Mathematics, Purdue University, 150 North University Street, West Lafayette, Indiana 47907, USA
\and
Debdas Ghosh (debdas.mat@itbhu.ac.in) \at
Department of Mathematical Sciences, Indian Institute of Technology (Banaras Hindu University) Varanasi, Uttar Pradesh  221005, India
\and
Xiangxiong Zhang (zhan1966@purdue.edu) \at 
Department of Mathematics, Purdue University, 150 North University Street, West Lafayette, Indiana 47907, USA
}
 \date{}

\maketitle
\begin{abstract}
This work presents a new three-operator splitting method to handle monotone inclusion and convex optimization problems. \textcolor{black}{The proposed splitting serves as another natural extension of the Douglas-Rachford splitting technique to problems involving three operators.   For solving a composite convex minimization of a sum of three functions, 
its formula resembles but is different from Davis-Yin splitting and the dual formulation of  the classical three-block ADMM.
Numerical tests suggest that such a splitting  scheme is robust in the sense of allowing larger step sizes. When  two functions have orthogonal domains, the splitting operator can be proven $1/2$-averaged, which implies convergence of the iteration scheme  using any positive step size. }  
\end{abstract}


\section{Introduction}
\textcolor{black}{
Operator splitting methods break down complex problems into smaller, more manageable subproblems that can be tackled either one after another or at the same time. While these methods have been around for over six decades, their relevance has grown considerably in recent years. They have proven to be highly effective across multiple fields ranging—from partial differential equations, control theory, and high dimensional problems in machine learning, image analysis, signal processing, statistical estimation, image processing, compressive sensing, finance, matrix completion, and control\cite{shi2014linear,wei2012distributed,boyd2011distributed}.}\\

\textcolor{black}{In 1950s, Peaceman-Rachford  \cite{peaceman1955numerical} and Douglas-Rachford \cite{douglas1956numerical} operator splitting techniques were first introduced for solving the heat equation in two dimensions. In \cite{lions1979splitting}, Lions and Mercier extended these methods to address the formulation involving the sum of two operators, each being maximal and monotone, and the method in \cite{lions1979splitting} is named  Douglas-Rachford splitting (DRS). In \cite{raguet2013generalized,raguet2019note}, it was shown that Douglas-Rachford and Forward-Backward splitting can be integrated, followed by its extension in \cite{briceno2015forward}. } \textcolor{black}{Such a Forward-Douglas-Rachford (FDR) splitting was proven convergent by Davis and Yin in \cite{davis2017three}, and also referred to as Davis-Yin splitting. It is also well known that DRS is equivalent to the Alternating Direction Method of Multipliers (ADMM), which was studied in the dual setting by \cite{gabay1976dual,glowinski1975approximation}.  DRS and ADMM have been a popular tool in many applications. In particular, Goldstein and Osher \cite{goldstein2009split} later advanced the application of ADMM through the development of the split Bregman method, particularly in imaging and sparse recovery contexts. \\
}\\ 

\textcolor{black}{In this paper, we propose 
a splitting scheme which can be used for solving an inclusion problem 
$0\in (\mathbb{A}+\mathbb{B}+\mathbb{C})x,$
with operators $\mathbb{A},\mathbb{B},\mathbb{C}$ \textcolor{black}{being three maximal monotone mappings over a real Hilbert space} $\mathcal{X}$ and the operator $\mathbb{C}$ is cocoerceive with parameter $\beta$. 
Let $\mathbb{J}_{\mathbb A}=(I+\mathbb{A})^{-1}$ denote
the resolvent of a monotone operator $\textcolor{black}{\mathbb{A}}$. 
The proposed splitting operator is given as
\begin{equation}
    T=\mathbb{J}_{\gamma\mathbb{C}}\circ(\mathbb{J}_{\gamma \mathbb{A}}\circ(2\mathbb{J}_{\gamma \mathbb{B}}-I-\gamma \mathbb{C}\circ \mathbb{J}_{\gamma \mathbb{B}})+\gamma \mathbb{C}\circ \mathbb{J}_{\gamma \mathbb{B}})+ (I-\mathbb{J}_{\gamma \mathbb{B}}).
    \label{intro-new-operator}
\end{equation}
Such a splitting is inspired by  the Davis-Yin splitting \cite{davis2017three}, which can be written as
\[ T_{\text{DY}}= \mathbb{J}_{\gamma \mathbb{A}}\circ(2\mathbb{J}_{\gamma \mathbb{B}}-I-\gamma \mathbb{C}\circ \mathbb{J}_{\gamma \mathbb{B}}) + (I-\mathbb{J}_{\gamma \mathbb{B}}).\]
 \\
 The main application of such a splitting is to solve a  composite convex minimization in the form of} 
\[
\min_x d_1(x)+d_2(x)+d_3(x),
\]
\textcolor{black}{where  $d_i(x)$ are convex closed and proper functions, and their proximal mapping $$\operatorname{prox}_{\gamma d_i}(y):=\operatorname{argmin}_x d_i(x)+\frac{1}{2\gamma}\|x-y\|^2$$ can be efficiently computed.} When one of the functions has Lipschitz-continuous gradient, e.g., assume $\nabla d_2$ is Lipschitz continuous, characterized by a Lipschitz constant $L$, 
\textcolor{black}{our proposed splitting scheme with a step size $\gamma>0$ is written as}
\begin{eqnarray}
    x^{k+\frac12}&=&\operatorname{prox}_{\gamma d_3} (z^k)\notag \\
\mbox{Proposed Splitting}:\quad      p^{k+1}&=& \operatorname{prox}_{\gamma d_1} (2x^{k+\frac12}-z^k-\gamma \nabla d_2(x^{k+\frac12})) \label{intro:news}\\
    x^{k+1}&=&\operatorname{prox}_{\gamma d_2} (p^{k+1}+\gamma\nabla d_2(x^{k+\frac12})) \notag \\
    z^{k+1}&=&z^k+(x^{k+1}-x^{k+\frac12}).\notag
\end{eqnarray} 
The Davis-Yin splitting (or Forward Douglas-Rachford splitting) \cite{davis2017three,raguet2019note}  can be written as
\begin{eqnarray}   x^{k+\frac12}&=&\operatorname{prox}_{\gamma d_3}(z^k)\notag \\
\mbox{Davis-Yin Splitting}:\quad    x^{k+1}&=& \operatorname{prox}_{\gamma d_1}  (2x^{k+\frac12}-z^k-\gamma \nabla d_2(x^{k+\frac12})) \label{intro:dys}\\
    z^{k+1}&=&z^k+(x^{k+1}-x^{k+\frac12}).\notag
\end{eqnarray} 
\textcolor{black}{As will be shown in Section \ref{sec:admm-derivation}, the dual form of the classical three-block ADMM method with a penalty parameter or step size $\gamma>0$ can be written as 
\begin{eqnarray}
    x^{k+\frac12}&=&\operatorname{prox}_{\gamma d_3} (z^k)\notag \\
\mbox{ADMM Dual Form}:\quad      p^{k+1}&=& \operatorname{prox}_{\gamma d_1} (2x^{k+\frac12}-z^k-\gamma \nabla d_2(x^{k})) \\
    x^{k+1}&=&\operatorname{prox}_{\gamma d_2} (p^{k+1}+\gamma\nabla d_2(x^{k})) \notag \\
    z^{k+1}&=&z^k+(x^{k+1}-x^{k+\frac12}).\label{intro:admm}\notag
\end{eqnarray} 
}

\textcolor{black}{Our proposed scheme is motivated by the dual form of the empirically efficient, classical three-block ADMM, from which we derive a new method to solve \eqref{three_operator_prob}. In numerical experiments, the original dual form converges only when the step size is constrained within a very limited range. By introducing a simple modification to this dual form, our new scheme significantly widens the permissible range of the step size, thereby improving its robustness and convergence properties.} We can see that the proposed new scheme \eqref{intro:news} is similar to but different from the Davis-Yin scheme \eqref{intro:dys} and the ADMM dual form \eqref{intro:admm}.  
The proposed splitting method may appear less efficient than the Davis-Yin scheme, as it involves evaluating three proximal operators instead of two. On the other hand, the extra computation of $\operatorname{prox}_{\gamma d_2}$ might improve the robustness of the splitting.
In  \cite{davis2017three}, the Davis-Yin splitting was proven to converge for any constant step size $\gamma\in (0,\frac{2}{L})$. Numerically, when the step size $\gamma$ is much larger than $\frac{2}{L}$, the Davis-Yin splitting \eqref{intro:dys} and \textcolor{black}{and the ADMM dual form \eqref{intro:admm}} will not converge, and the proposed splitting method \eqref{intro:news} can still converge, as will be shown in numerical examples in Section \ref{section_numerical}.

\textcolor{black}{Notice that 
the proposed new scheme \eqref{intro:news} can be regarded as a modified version of the ADMM dual form \eqref{intro:admm}.
Thus the proposed splitting scheme \eqref{intro:news} can be also be formally extended to multiple operators by similarly modifying the dual form of multiple-block ADMM methods as will be shown in Section \ref{section_multiblock}.  
In the literature, there are several studies of ADMM methods that address the generic case where the number of blocks $m\geq3$. 
 Despite ADMM's widespread success in two-block convex optimization, the direct extension of ADMM to three or more blocks was shown to be not necessarily convergent \cite{chen2016directNO}. A key result in this direction is the provision of a concrete example demonstrating divergence of the direct three-block ADMM, marking a fundamental limitation of naive extensions. To address this, Cai et al. \cite{cai2017convergence} proposed a convergent three-block ADMM, assuming strong convexity of one function and establishing global convergence under additional structural conditions. Later, He and Yuan introduced a generalized symmetric ADMM \cite{he2018gsadmm} which permits larger dual step sizes while retaining convergence guarantees for separable multi-block convex programs. Beyond ADMM, notable progress has also been made in Forward-Backward-type primal-dual algorithms. The classical Chambolle-Pock (CP) algorithm \cite{chambolle2011first} was widely used for imaging problems and saddle-point formulations, but its behavior under inertial parameters in the range 
$(0,1)$ was not well understood. Recent works such as \cite{tian2022gpd} and \cite{he2025generalizedadjoint} addressed this by introducing corrected Forward-Backward-Adjoint frameworks. These not only clarified convergence in the inertial regime but also resolved long-standing questions about the robustness of primal-dual methods when extrapolation or momentum steps are introduced. These developments form a critical part of the ongoing effort to restore convergence in multi-block settings—an area of direct relevance to our work.}
In \cite{han2012note}, the strong convexity assumption has been assumed on all the given objective functions. Lin et al. in \cite{lin2015global} considered $(m-1)$ functions to be strongly convex and established global convergence without imposing any restrictions on the penalty parameter. \textcolor{black}{Furthermore, Lin et al. \cite{lin2015sublinear} established linear convergence assuming that the objective functions satisfy Lipschitz continuity.} In a sequel, the work in \cite{hong2017linear} showed that linear convergence is ensured if the step size in each updating step is sufficiently reduced and \textcolor{black}{an appropriate error-bound assumption} holds. Moreover, it has been discussed in \cite{cai2017convergence,chen2013convergence,han2012note,li2015convergent,lin2015sublinear} \textcolor{black}{that the penalty parameter must be suitably controlled for favourable convergence behaviour}. The importance of restricting the penalty parameter to ensure faster convergence has been discussed extensively in \cite{cai2017convergence,chen2013convergence,han2012note,li2015convergent,lin2015sublinear}. Although the restriction can be conservative to ensure convergence, it may be relaxed to achieve faster rates. \textcolor{black}{In \cite{davis2016convergence}, Davis and Yin  proposed a convergent three-block ADMM variant by assuming strong convexity in one objective component and the step size parameter is bounded by a threshold value.} In \cite{lin2016iteration}, Lin et al. proposed several alternative methods for three-block ADMM without any restriction on penalty parameters to solve regularized least square decomposition problems.

 \textcolor{black}{
The Davis-Yin splitting operator $ T_{\text{DY}}$ is an averaged operator for any $\gamma\in(0,2\beta)$ thus converges under suitable assumptions, as proven in \cite{davis2017three}. 
 It is however nontrivial to prove the averagedness of the operator $T$ defined in \eqref{intro-new-operator} under the same assumptions, due to its more complicated structure.
 In Section \ref{sec:specialcase}, we consider a special assumption that $d_1(x)$ and $d_2(x)$ have orthogonal domains, under which the new splitting operator can be proven averaged thus the proposed scheme \eqref{intro:news} converges with any step size $\gamma>0$.
  Such a special assumption is similar to but weaker than the  sufficient condition $A_1^\top A_2=0$ to ensure convergence of three-block ADMM in \cite{chen2016directNO}, as will be explained Remark \ref{rem-convergence}.
  On the other hand, the Davis-Yin splitting still needs the step size constraint $\gamma<\frac{2}{L}$ for convergence even under this special assumption, as will be explained in Remark \ref{rem:advantage}.}

 \textcolor{black}{The main contributions of this paper include the introduction to a new three-operator splitting scheme.  As will be reviewed in Section \ref{subsection_comparison}, there are other three-operator splitting schemes. To the best of our knowledge,  \eqref{intro-new-operator} is different from the existing three-operator splitting schemes.  As will be shown in Section \ref{section_numerical}, the scheme \eqref{intro:news} uses only one parameter $\gamma$ and it is numerically more robust when $\gamma>\frac{2}{L}$ than existing one parameter three-operator splitting schemes such as Davis-Yin splitting. Though there are other three-operator splitting schemes which can also allow much  larger step sizes, these schemes usually involve  tuning another parameter for convergence. In other words, the main practical advantage of \eqref{intro:news}  is its robustness without further tuning parameters. Allowing large step sizes could be useful when the Lipschitz constant  $L$ is unknown and hard to estimate.  We also prove its convergence using any step size $\gamma>0$ for the composite minimization under one special  assumption that $d_1$ and $d_2$ have orthogonal domains.} 
\textcolor{black}{In the Appendix, we derive the primal form of the proposed splitting scheme \eqref{intro:news}, which gives a 3-block ADMM type scheme. Compared to many existing 3-block ADMM variants, this particular variant does not seem to have any advantage, and the Appendix only serves the purpose of comparing it with existing methods for interested readers. On the other hand, the proposed operator splitting \eqref{intro-new-operator} can also be used to solve an operator inclusion problem like $0\in (\mathbb{A}+\mathbb{B}+\mathbb{C})x,$ which however ADMM type methods do not  directly solve. 
We emphasize that there are some important modern data science applications which can only be formulated as such a monotone operator inclusion problem but not a convex minimization problem, e.g., \cite{nurbekyan2024monotone}.}

\textcolor{black}{
The rest of this paper is organized as follows. Section \ref{section_basic_notation} describes the notation and symbols. In Section \ref{section_three_operator_splitting}, we first review some existing three-operator splitting schemes then introduce the proposed new splitting operator.   Section \ref{section_weak_convergence_rates} includes standard results under the assumption that the splitting operator is averaged. In Section \ref{section_multiblock}, we prove the splitting operator is averaged for the composite convex minimization under the assumption that two functions have orthogonal domains. Section \ref{section_numerical} presents numerical examples, followed by concluding remarks in Section \ref{section_conclusion}.}

\section{Basic Notation and Fundamental Results}\label{section_basic_notation}
\textcolor{black}{We adopt the following notation for use throughout the paper.}
\begin{enumerate}[$\bullet$]
    \item $\mathcal{X}$ denotes an infinite dimensional Hilbert space
    \item \textcolor{black}{$\langle \cdot ,\cdot \rangle$} denotes inner product associated to $\mathcal{X}$
    \item $(\lambda_j)_{j\geq 0}\subseteq\mathbb{R}_{++}$ denotes a stepsize sequence.
\end{enumerate}
\textcolor{black}{We now present some standard definitions and basic facts that will be used throughout.\\
\\
A map $F:S\to\mathcal{X}$, defined on a nonempty subset $S\subseteq \mathcal{X}$ is said to be $L$-Lipschitz for $L\geq0$ if for every $x,y\in S$
\[
\lVert F(x)-F(y) \rVert \leq L \lVert x-y\rVert. 
\]
$F$ is said to be nonexpansive if the above inequality holds for $L=1$.}\\
\\
Let $I_{\mathcal{X}}$ be the identity map. A map $F_{\alpha}:S\to\mathcal{X}$ is called $\alpha$-averaged if it can be written as 
\[
F_{\alpha}=(1-\alpha)I_{\mathcal{X}}+\alpha F,
\]
where $F$ is some nonexpansive map.
Moreover, \textcolor{black}{if $F$ is $(1/2)$-averaged, then it is said to be firmly nonexpansive. For convenience, we may use $I$ in place of $I_{\mathcal{X}}$.\\
\\
Let $2^{\mathcal{X}}$ denotes the power set of $\mathcal{X}$. An operator $\mathbb{\textcolor{black}{A}}:\mathcal{X}\to 2^{\mathcal{X}}$ is said to be monotone if for every $x,y\in\mathcal{X},~u\in \textcolor{black}{\mathbb{A}}x,~v\in \textcolor{black}{\mathbb{A}}y$, the inequality
\[
\langle x-y,u-v \rangle\geq0 \text{ holds}.
\]
The solution set (also called the zeroes set) of a monotone operator $\mathbb{\textcolor{black}{A}}$ is defined as}
\[
\text{zer}(\mathbb{\textcolor{black}{A}})=\{x\in\mathcal{X}~|~0\in \mathbb{\textcolor{black}{A}}x\}.
\]
\textcolor{black}{An operator $\mathbb{\textcolor{black}{A}}$ is strongly monotone with parameter $\beta>0$ if
\[
\langle x-y,u-v \rangle\geq \beta \lVert x-y \rVert^2,
\]
for all $x,y\in\mathcal{X}$ and $u\in \mathbb{\textcolor{black}{A}}x,~v\in \mathbb{\textcolor{black}{A}}y$.}\\
An operator $\mathbb{\textcolor{black}{A}}$ is called cocoerceive with constant $\beta>0$ if
\[
\langle u-v,x-y\rangle \geq \beta \lVert u-v\rVert^2, \text{ for all }x,y\in\mathcal{X},~u\in \mathbb{\textcolor{black}{A}}x,~v\in \mathbb{\textcolor{black}{A}}y.
\]
Additionally, from the definition of $\mathbb{\textcolor{black}{A}}$, Cauchy-Schwarz inequality and $\beta$-cocoerceive of $\mathbb{\textcolor{black}{A}}$, it follows that 
\begin{equation}\label{cocoerceive_lipschitz}
\lVert x-y\rVert\geq \beta \lVert u-v\rVert, \text{ for all }x,y\in\mathcal{X},~u\in \mathbb{\textcolor{black}{A}}x,~v\in \mathbb{\textcolor{black}{A}}y,~\beta>0.
\end{equation}
Moreover, if a convex function $f$ has $L$-Lipschitz gradient, then $\nabla f$ is $1/L$-cocoerceive.\\
\\
The operator $\mathbb{\textcolor{black}{A}}$ has a graph defined as
\[
\text{gra}(\mathbb{\textcolor{black}{A}})=\{(x,y)|~x\in\mathcal{X}, y\in \mathbb{\textcolor{black}{A}}x\}.
\]
This $\text{gra}(\mathbb{\textcolor{black}{A}})$ is called maximal monotone if it is not a proper subset of the graph of any other monotone operator.\\
\\
The graph of the inverse operator is defined by
\[
\text{gra}(\mathbb{\textcolor{black}{A}}^{-1})=\{(y,x)|x\in\mathcal{X}, y\in \mathbb{\textcolor{black}{A}}x\}.
\]
The resolvent and reflection of a monotone operator $\mathbb{\textcolor{black}{A}}$ is denoted by $\mathbb{J}_{\mathbb{\textcolor{black}{A}}}$ and $\textcolor{black}{\mathbf{R}}_\mathbb{\textcolor{black}{A}}$, respectively, and defined by
\[
\mathbb{J}_{\mathbb{\textcolor{black}{A}}}=(I+\mathbb{\textcolor{black}{A}})^{-1} \text{ and }\textcolor{black}{\textcolor{black}{\mathbf{R}}_\mathbb{\textcolor{black}{A}}}=2\mathbb{J}_{\mathbb{\textcolor{black}{A}}}-I.
\]
For maximal monotone $\mathbb{\textcolor{black}{A}}$, $\textcolor{black}{\mathbf{R}}_\mathbb{\textcolor{black}{A}}$ is nonexpansive.\\
\\
\textcolor{black}{For an extended function $f:\mathcal{X}\to(-\infty, \infty]$ that is proper, closed (or equivalently {\it lower semi-continuous}), and convex, its subdifferential set at $x$ is denoted by a map $\partial f:\mathcal{X}\to2^{\mathcal{X}}$,   
\[
\partial f(x)=\{g\in\mathcal{X}|f(y)\geq f(x)+\textcolor{black}{\langle y-x,g\rangle }, \text{ for all }y\in\mathcal{X} \}.
\]
Moreover, if $f$ is differentiable at $x$, then
$\nabla f(x)\in \partial f(x)$.\\
The convex (or \textcolor{black}{Fenchel}) conjugate of $f$ is given by
\[
f^*(y)=\textcolor{black}{\underset{x\in\mathcal{X}}{\sup} \{ \langle y,x \rangle-f(x) \}.}
\]
The indicator function $i_C(x)$ of a closed convex set $C\subseteq \mathcal{X}$ is
\[
i_C(x)=\begin{cases}
    0 &x\in C\\
    +\infty & x\in\mathcal{X}/C,
\end{cases}
\]
The proximal and reflection operators for $f$ with a step size $\lambda>0$ are 
\[
\operatorname{prox}_{\lambda f}(x)=\underset{y\in\mathcal{X}}{\text{arg min}}(f(y)+\tfrac{1}{2\lambda}\lVert y-x \rVert^2)\text{ and }\textcolor{black}{\mathbf{R}}_{\textcolor{black}{\lambda f}}=2\operatorname{prox}_{\lambda f}-I. 
\]}
In this article, we use the following \emph{cosine rule} and \emph{Young's Inequality} given by:
\begin{align}
&2\langle y-x,z-x  \rangle =\lVert y-x\rVert^2+\lVert z-x\rVert^2-\lVert y-z\rVert^2\text{ for all }x,y,z\in\mathcal{X}\label{cosine_rule}\\
&\text{and }ab\leq \frac{a^2}{2\varepsilon}+\frac{b^2\varepsilon}{2}\text{ for all }a,b\geq0,~\varepsilon>0,\text{ respectively}\label{inequality_rule}.
\end{align}

\section{A Three-Operator Splitting Scheme}\label{section_three_operator_splitting}
In this section, we introduce a \textcolor{black}{new} three-operator splitting scheme, which can be used to solve nonsmooth and monotone inclusion optimization problems of many different forms. We consider the problem
\begin{eqnarray}\label{three_operator_prob}
 \text{ find }x\in\mathcal{X} \text{ such that }0\in (\mathbb{A}+\mathbb{B}+\mathbb{C})x,   
\end{eqnarray}
where $\mathbb{A},\mathbb{B},\mathbb{C}$ are three maximal monotone operators defined on a \textcolor{black}{real} Hilbert space $\mathcal{X}$ and the operator $\mathbb{C}$ is cocoerceive with parameter $\beta$.

\subsection{Existing three-operator splitting schemes}\label{subsection_comparison}
We first review a few existing three-operator splitting schemes.
For solving the problem \eqref{three_operator_prob}, there are at least the following three-operator splitting methods:
\begin{enumerate}
    \item The Forward Douglas-Rachford splitting is also known as  Davis-Yin splitting $ T_{\text{DY}}$, and its convergence was proven by Davis and Yin in  \cite{davis2017three}. See also \cite{raguet2013generalized,raguet2019note} for the extension of Forward Douglas-Rachford splitting to multiple operators. 
\item In \cite{malitsky2020forward}, Malitsky and Tam 
introduced a Forward-Reflected-Backward splitting method with a generalization for three operators. This method requires $\mathbb{A}$ to be maximal monotone, $\mathbb{B}$ to be monotone and Lipschitz, and $\mathbb{C}$ to be cocoercive. The iteration scheme can be written as
\begin{equation}
x_{k+1} = \textcolor{black}{\mathbb{J}}_{\gamma \mathbb{A}}(x_k - 2\gamma \mathbb{B}(x_k) + \gamma \mathbb{B}(x_{k-1}) - \gamma \mathbb{C}(x_k)).    
\end{equation} 
\item In \cite{ryu2020finding},
Ryu and Vu combined Douglas-Rachford splitting with Forward-Backward-Forward splitting (FDRF) under the assumption that two operators are maximal monotone and the other is monotone and Lipschitz, 
\begin{equation}
 T_{\text{FDRF}} := (I - \gamma \mathbb{C})\circ \textcolor{black}{\mathbb{J}}_{\gamma \mathbb{A}} \circ(2 \textcolor{black}{\mathbb{J}}_{\gamma \mathbb{B}} - I - \gamma \mathbb{C} \circ \textcolor{black}{\mathbb{J}}_{\gamma \mathbb{B}}) + I - (I - \gamma \mathbb{C}) \circ \textcolor{black}{\mathbb{J}}_{\gamma \mathbb{B}}.
\end{equation}
The  FDRF iteration can be written as
\begin{equation}
\begin{aligned}
x_{n+1} &= \textcolor{black}{\mathbb{J}}_{\gamma \mathbb{B}} z_n\\
y_{n+1} &= \textcolor{black}{\mathbb{J}}_{\gamma \mathbb{A}}(2x_{n+1} - z_n - \gamma \mathbb{C} x_{n+1}) \\
z_{n+1} &= z_n + y_{n+1} - x_{n+1} - \gamma (\mathbb{C} y_{n+1} - \mathbb{C} x_{n+1}) \quad \text{(FDRF).} 
\end{aligned}
\label{FDRF}
\end{equation}
For relaxing assumptions needed for convergence, in \cite{ryu2020finding} Ryu and Vu also
 proposed  a method that combines DR with Forward-Reflected-Backward Splitting (FRDR) with an extra parameter $\beta$, which can be written as
\begin{equation}
\begin{aligned}
x_{n+1} &= \textcolor{black}{\mathbb{J}}_{\gamma \mathbb{B}}\big(x_n - \gamma u_n - \gamma (2\mathbb{C} x_n - \mathbb{C} x_{n-1})\big)\\
y_{n+1} &= \textcolor{black}{\mathbb{J}}_{\beta \mathbb{A}}(2x_{n+1} - x_n + \beta u_n) \\
u_{n+1} &= u_n + \tfrac{1}{\beta}(2x_{n+1} - x_n - y_{n+1})  \quad \text{(FRDR).} 
\end{aligned}
\label{FRDR}
\end{equation}
\end{enumerate}
One application for these three-operator splitting schemes is to   solve a composite convex minimization problem in the following form
\begin{equation}
    x^* \in \operatorname*{arg\,min}_{x \in \mathcal{X}} \{ d_1(x) + d_2(x) + D_3(x) \},
\end{equation}
where $d_1,~d_2$, and $D_3$ are convex, closed   and proper. Operator splitting schemes apply to such a problem with $\mathbb{A} = \partial d_1$, $\mathbb{B} = \partial d_2$, $\mathbb{C} = \partial D_3$, and $\textcolor{black}{\mathbb{J}}_{\gamma \partial D_3}$ is the proximal operator $\operatorname{prox}_{\gamma D_3}$ for the function $D_3$.
For the case $D_3(x)=d_3(Bx)$ where $B$ is a matrix, Yan introduced the following primal-dual algorithm (PD3O) in \cite{yan2018new}:
\begin{align}
x_k     &= \operatorname{prox}_{\gamma d_1}(z_k) \notag\\
s_{k+1} &= \operatorname{prox}_{\delta d_3^*}\left( (I - \gamma \delta BB^\top)s_k + \delta B^\top(2x_k - z_k - \gamma \nabla d_2(x_k)) \right) \label{scheme-PD3O}\\
z_{k+1} &= x_k - \gamma \nabla d_2(x_k) - \gamma B^\top s_{k+1}  \notag \quad \text{(PD3O)}.
\end{align} 
The PD3O method reduces to Davis-Yin splitting 
 when $\textcolor{black}{B} = I$, and it is preferred when it is much easier to compute $\textcolor{black}{\mathbb{J}}_{\partial d_3}$ than $\textcolor{black}{\mathbb{J}}_{\partial D_3}$, e.g., $D_3$ is the total variation (TV) norm function and $d_3$ is the $\ell^1$-norm function.

\subsection{A new three-operator splitting}

\textcolor{black}{Next, we derive a new  three-operator splitting. To the best of our knowledge, it is different from all existing three-operator splitting methods in the literature}.
For any $\gamma>0$, we have
\begin{eqnarray}\label{operator_T}
&&0\in \gamma(\mathbb{A}+\mathbb{B}+\mathbb{C})x\nonumber\\ 
&&\iff 0\in (I+\gamma \mathbb{A})x-(I-\gamma \mathbb{B})x+\gamma \mathbb{C}x\nonumber\\
&&\iff0\in(I+\gamma \mathbb{A})x-\textcolor{black}{\mathbf{R}}_{\gamma \mathbb{B}}(I+\gamma \mathbb{B})x+\gamma \mathbb{C}x~~[\text{since }\textcolor{black}{\mathbf{R}}_{\gamma\mathbb{B}}(I+\gamma\mathbb{B})=(I-\gamma\mathbb{B})]\nonumber\\
&&\iff0\in(I+\gamma \mathbb{A})x-\textcolor{black}{\mathbf{R}}_{\gamma \mathbb{B}}z+\gamma \mathbb{C}x~~[\text{assume }z\in(I+\gamma\mathbb{B})x]\nonumber\\
&&\iff\textcolor{black}{\mathbf{R}}_{\gamma \mathbb{B}}z-\gamma \mathbb{C}x\in(I+\gamma \mathbb{A})x\nonumber\\
&&\iff2\mathbb{J}_{\gamma \mathbb{B}}z-z-\gamma \mathbb{C}x\in(I+\gamma \mathbb{A})x~~[\text{since }\textcolor{black}{\mathbf{R}}_{\gamma \mathbb{B}}z=(2\mathbb{J}_{\gamma \mathbb{B}}-I)z]\nonumber\\
&&\iff\mathbb{J}_{\gamma \mathbb{A}}(2\mathbb{J}_{\gamma \mathbb{B}}z-z-\gamma \mathbb{C}x)= x~~[\text{since }\mathbb{J}_{\gamma \mathbb{A}}=(I+\gamma\mathbb{A})^{-1}]\nonumber\\
&&\iff\mathbb{J}_{\gamma \mathbb{A}}(2\mathbb{J}_{\gamma \mathbb{B}}z-z-\gamma \mathbb{C}\circ \mathbb{J}_{\gamma \mathbb{B}}z)= x\nonumber\\
&&\iff\mathbb{J}_{\gamma \mathbb{A}}(2\mathbb{J}_{\gamma \mathbb{B}}z-z-\gamma \mathbb{C}\circ \mathbb{J}_{\gamma \mathbb{B}}z)+\gamma \mathbb{C}\circ \mathbb{J}_{\gamma \mathbb{B}}z= x+\gamma \mathbb{C}\circ \mathbb{J}_{\gamma \mathbb{B}}z\nonumber\\
&&\iff\mathbb{J}_{\gamma \mathbb{A}}(2\mathbb{J}_{\gamma \mathbb{B}}z-z-\gamma \mathbb{C}\circ \mathbb{J}_{\gamma \mathbb{B}}z)+\gamma \mathbb{C}\circ \mathbb{J}_{\gamma \mathbb{B}}z= (I+\gamma \mathbb{C})x\nonumber\\
&&\iff\mathbb{J}_{\gamma\mathbb{C}}\left(\mathbb{J}_{\gamma \mathbb{A}}(2\mathbb{J}_{\gamma \mathbb{B}}-I-\gamma \mathbb{C}\circ \mathbb{J}_{\gamma \mathbb{B}})+\gamma \mathbb{C}\circ \mathbb{J}_{\gamma \mathbb{B}}\right)z= x~~[\text{since }\mathbb{J}_{\gamma \mathbb{C}}=(I+\gamma\mathbb{C})^{-1}]\nonumber\\
&&\iff\mathbb{J}_{\gamma\mathbb{C}}\left(\mathbb{J}_{\gamma \mathbb{A}}(2\mathbb{J}_{\gamma \mathbb{B}}-I-\gamma \mathbb{C}\circ \mathbb{J}_{\gamma \mathbb{B}})+\gamma \mathbb{C}\circ \mathbb{J}_{\gamma \mathbb{B}}\right)z= \mathbb{J}_{\gamma \mathbb{B}}z~~[\text{assumption }z\in(I+\gamma\mathbb{B})x]\nonumber\\
&&\iff\left(\mathbb{J}_{\gamma\mathbb{C}}\circ(\mathbb{J}_{\gamma \mathbb{A}}\circ(2\mathbb{J}_{\gamma \mathbb{B}}-I-\gamma \mathbb{C}\circ \mathbb{J}_{\gamma \mathbb{B}})+\gamma \mathbb{C}\circ \mathbb{J}_{\gamma \mathbb{B}})+(I-\mathbb{J}_{\gamma \mathbb{B}})\right)z=z,
\end{eqnarray}
where $T=\mathbb{J}_{\gamma\mathbb{C}}\circ(\mathbb{J}_{\gamma \mathbb{A}}\circ(2\mathbb{J}_{\gamma \mathbb{B}}-I-\gamma \mathbb{C}\circ \mathbb{J}_{\gamma \mathbb{B}})+\gamma \mathbb{C}\circ \mathbb{J}_{\gamma \mathbb{B}})+ (I-\mathbb{J}_{\gamma \mathbb{B}})$.
 
The proposed operator $T$ in \eqref{operator_T} splits the three-operator sum problem given in \eqref{three_operator_prob} into simpler sub-problems. \textcolor{black}{For special cases, Algorithm \ref{operator_T_algorithm1} reduces to Douglas-Rachford splitting (DRS) and a splitting closely related to Forward-Backward splitting (FBS)}. The proposed operator $T$  reduces to two-operator splittings for the two-operator sum problem given by 
\[
\text{Find }x\in\mathcal{X} \text{ such that }0\in \mathbb{A}x+\mathbb{B}x.
\]
The two special cases for the proposed operator $T$ are as follows:
\begin{enumerate}[(i)]
    \item \label{operator_part1} If $\mathbb{A}=0$,  \eqref{operator_T} reduces to
    \[
    T=\mathbb{J}_{\gamma\mathbb{C}}\circ(2\mathbb{J}_{\gamma \mathbb{B}}-I)+ (I-\mathbb{J}_{\gamma \mathbb{B}}), \text{ which is DRS  \cite{lions1979splitting}.}
    \]
    
    \item \label{operator_part2} If $\mathbb{C}=0$, \eqref{operator_T} reduces to
    \[
    T=\mathbb{J}_{\gamma \mathbb{A}}\circ(2\mathbb{J}_{\gamma \mathbb{B}}-I)+ (I-\mathbb{J}_{\gamma \mathbb{B}}), \text{ which is DRS  \cite{lions1979splitting}.}
    \]  
    \item \label{operator_part3} If $\mathbb{B}=0$, \eqref{operator_T} reduces to
\[\mathbb{J}_{\gamma\mathbb{C}}\circ(\mathbb{J}_{\gamma \mathbb{A}}\circ(I-\gamma \mathbb{C})+\gamma \mathbb{C})z=z.\]
The splitting above is not FBS \cite{lions1979splitting,passty1979ergodic} but it is closely related to FBS since     
    \begin{eqnarray*}
     & & \mathbb{J}_{\gamma\mathbb{C}}\circ(\mathbb{J}_{\gamma \mathbb{A}}\circ(I-\gamma \mathbb{C})+\gamma \mathbb{C})z=z\\
     &\iff& \mathbb{J}_{\gamma \mathbb{A}}\circ(I-\gamma \mathbb{C})z+\gamma \mathbb{C}z=(I+\gamma \mathbb{C})z\\
     &\iff& \mathbb{J}_{\gamma \mathbb{A}}\circ(I-\gamma \mathbb{C})z=z.
    \end{eqnarray*}
\end{enumerate}

For solving \eqref{three_operator_prob}, consider the  Krasnosel’ski{\u\i}--Mann (KM) iteration with the operator $T$ above:
\begin{eqnarray}
&&T_{\lambda}:=(1-\lambda)I+\lambda T\nonumber\\
&\text{and }& z^{k+1}=(1-\lambda_k)z^k+\lambda_k T z^k.\label{operator_z^k+1}
\end{eqnarray}
If $T$ is $\alpha$-averaged, then the classical fixed point iteration theorem states that the iteration $z_{k+1}=T_\lambda(z_k)$ converges for any constant $\lambda\in(0,\frac{1}{\alpha}].$ 
See \cite{boct2023fast} and references therein for some of the latest developments of strategies for designing $\lambda_k$.
\textcolor{black}{The steps of the iterative procedure are given in Algorithm \ref{operator_T_algorithm1}:}
\begin{algorithm}[!htb]
\caption{Initialize $z^0\in\mathcal{X},~\gamma\in(0,2\beta]$, and sequence $(\lambda_k)_{k\geq0}\in(0,(4\beta-\gamma)/2\beta)$. For $k=0,1,2,\ldots$}\label{operator_T_algorithm1}
\begin{enumerate}
\item Compute $ x^{k}_\mathbb{B}=\mathbb{J}_{\gamma \mathbb{B}}(z^k)$;
\item Compute $x^{k}_\mathbb{A}=\mathbb{J}_{\gamma \mathbb{A}}(2x^{k}_\mathbb{B}-z^k-\gamma\mathbb{C}x^{k}_\mathbb{B})$ \hspace{0.5cm}\text{ or }~~ $\mathbb{J}_{\gamma \mathbb{A}}(2\mathbb{J}_{\gamma \mathbb{B}}(z^k)-z^k-\gamma\mathbb{C}\circ \mathbb{J}_{\gamma \mathbb{B}}(z^k))$;
\item Compute $x^{k}_\mathbb{C}=\mathbb{J}_{\gamma \mathbb{C}}( x^{k}_\mathbb{A}+\gamma\mathbb{C}x^{k}_\mathbb{B})$\hspace{1.52cm}\text{ or }~~$\mathbb{J}_{\gamma \mathbb{C}}( x^{k}_\mathbb{A}+\gamma\mathbb{C}\circ \mathbb{J}_{\gamma \mathbb{B}}(z^k))$;
\item Update $z^{k+1}=z^k+\lambda_k(x^{k}_\mathbb{C}-x^{k}_\mathbb{B}).$
\end{enumerate}
\end{algorithm}

\textcolor{black}{As a remark, for solving an inclusion problem like 
$$ x^* \in \operatorname*{arg\,min}_{x \in \mathcal{X}} \{ d_1(x) + d_2(x) + d_3(Bx) \},$$
when $B = I$, PD3O and other primal–dual algorithms such as Condat-Vũ \cite{condat2013aprimal,vu2011asplitting}, the Primal-Dual Fixed–Point algorithm \cite{chen2016aprimal}, and the Asymmetric Forward–Backward–Adjoint splitting \cite{latafat2017asymmetric} 
naturally give algorithms for a general inclusion problem \eqref{three_operator_prob}. On the other hand, these primal–dual methods were designed with motivations to exploit the explicit separation between $B$ and $d_3$, allowing the proximal operator of $d_3$ and the linear mapping $B$ to be handled efficiently.  Nonetheless, for solving a general inclusion problem \eqref{three_operator_prob}, the proposed new splitting is indeed different from these existing methods. Numerical tests in Section \ref{section_numerical} suggest that the proposed new splitting for solving composite convex minimization problems is quite robust with any positive step size. Furthermore, it can be regarded as a variant of the Davis-Yin splitting, and this variant has the numerical advantage of allowing much larger step sizes, which can be proven for a special problem in Section \ref{sec:specialcase}. }


\section{Weak Convergence and Rates of Proposed Three-Block 
Operator Splitting Scheme}\label{section_weak_convergence_rates}
 
\begin{figure}[H]
    \begin{center}
         \includegraphics[scale=0.36]{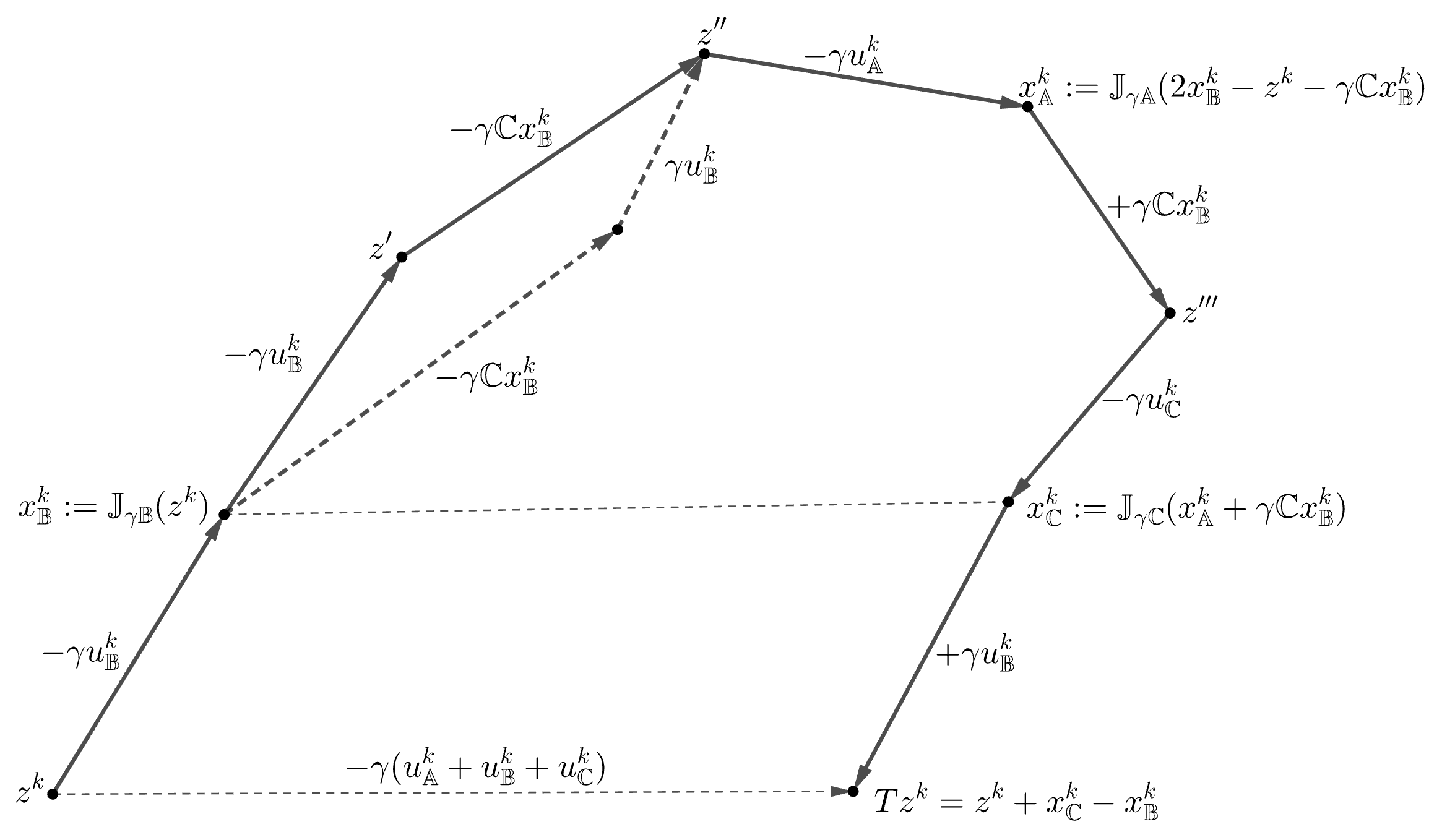}
    \caption{The proposed mapping $T:z^k\to Tz^k$ and vectors $u^k_{\mathbb{B}}\in \mathbb{B}x^k_{\mathbb{B}}$, $u^k_{\mathbb{A}}\in \mathbb{A}x^k_{\mathbb{A}}$, and $u^k_{\mathbb{C}}\in \mathbb{C}x^k_{\mathbb{C}}$ as given in Lemma \ref{weak_convergence_supporting_lem}.}
    \label{figure1}
    \end{center}  
\end{figure}

In this section, we discuss some properties of the operator $T$ defined in \eqref{operator_T}. Figure \ref{figure1} provides a geometric illustration of how $T$ acts on a point $z\in\mathcal{X}$ corresponding to the points defined in Lemma \ref{weak_convergence_supporting_lem}.

\begin{lem}\label{weak_convergence_supporting_lem}
For $z\in\mathcal{X}$, define the points:
\begin{eqnarray*}
x^k_{\mathbb{B}}=\mathbb{J}_{\gamma \mathbb{B}}(z^k),~ & z'=2x^k_{\mathbb{B}}-z^k,~ &u^k_{\mathbb{B}}=\gamma^{-1}(z^k-x^k_{\mathbb{B}})\in \mathbb{B}x^k_{\mathbb{B}}   \\
x^k_{\mathbb{A}}=\mathbb{J}_{\gamma \mathbb{A}}(z''),~ &z'''=x^k_{\mathbb{A}}+\gamma\mathbb{C}x^k_{\mathbb{B}},~ & u^k_{\mathbb{A}}=\gamma^{-1}(z''-x^k_{\mathbb{A}})\in \mathbb{A}x^k_{\mathbb{A}}\\    
z''=z'-\gamma \mathbb{C} x^k_{\mathbb{B}},~&x^k_{\mathbb{C}}=\mathbb{J}_{\gamma \mathbb{C}}(z''') ,~& u^k_{\mathbb{C}}=\gamma^{-1}(x^k_{\mathbb{A}}+\gamma \mathbb{C}x^k_{\mathbb{B}}-x^k_{\mathbb{C}})\in \mathbb{C}x^k_{\mathbb{C}}.
\end{eqnarray*}
\textcolor{black}{Based on the above relations, the following equations hold:}
\[
Tz^k- z^k= x^k_{\mathbb{C}}-x^k_{\mathbb{B}}=-\gamma  (u^k_{\mathbb{A}}+u^k_{\mathbb{B}}+u^k_{\mathbb{C}})\text{ and }Tz^k= x^k_{\mathbb{C}}+\gamma u^k_{\mathbb{B}}.
\]
\end{lem}

\begin{proof}
In view of the definition of $T$, we observe that 
\[
Tz^k= z^k + x^k_{\mathbb{C}}-x^k_{\mathbb{B}}= x^k_{\mathbb{C}}+\gamma u^k_{\mathbb{B}}.
\]
Now, we conclude that
\begin{align*}
Tz^k-z^k= x^k_{\mathbb{C}}-x^k_{\mathbb{B}}= x^k_{\mathbb{A}}+\gamma 
\mathbb{C} x^k_{\mathbb{B}}-\gamma u^k_{\mathbb{C}}-x^k_{\mathbb{B}}
&= 2x^k_{\mathbb{B}}-z^k-\gamma u^k_{\mathbb{A}}-\gamma u^k_{\mathbb{C}}-x^k_{\mathbb{B}}\\
&= x^k_{\mathbb{B}}-z^k-\gamma u^k_{\mathbb{A}}-\gamma u^k_{\mathbb{C}}\\
&=-\gamma  (u^k_{\mathbb{A}}+u^k_{\mathbb{B}}+u^k_{\mathbb{C}}).
\end{align*}
\end{proof}

Next, we show that the fixed point relation associated with the operator $T$ holds true. Moreover, with the help of any fixed point $z^*$ of $T$ and $\mathbb{J}_{\gamma\mathbb{B}}z^*$, a zero of $\mathbb{A}+\mathbb{B}+\mathbb{C}$ can be obtained.

\begin{lem}
Let $\mathbb{A},~\mathbb{B}$, and $\mathbb{C}$ be three operators. Then, we have the following set equality:
\begin{eqnarray*}
\text{zer}(\mathbb{A}+\mathbb{B}+\mathbb{C})=\mathbb{J}_{\gamma\mathbb{B}}(\text{Fix }T),
\end{eqnarray*}
where Fix $T=\{x+\gamma u|0\in (\mathbb{A}+\mathbb{B}+\mathbb{C})x,~u\in(\mathbb{B}x)\cap(-\mathbb{A}x-\mathbb{C}x)\}$.
\end{lem}

\begin{proof}
Let $x\in \text{zer}(\mathbb{A}+\mathbb{B}+\mathbb{C})$, that is, $0\in(\mathbb{A}+\mathbb{B}+\mathbb{C})x$. Let $u_{\mathbb{A}}\in \mathbb{A}x$ and $u_{\mathbb{B}}\in \mathbb{B}x$ such that $\gamma(u_{\mathbb{A}}+u_{\mathbb{B}}+\mathbb{C}x)=0$ and $z=x+\gamma u_{\mathbb{B}}$. We first show that $z$ is a fixed point of $T$. Notice that $x=\mathbb{J}_{\gamma\mathbb{B}}(z)$. We have
\begin{eqnarray*}
&&2\mathbb{J}_{\gamma \mathbb{B}}(z)-z-\gamma \mathbb{C}\circ \mathbb{J}_{\gamma \mathbb{B}}(z)=2x-z-\gamma\mathbb{C}x=x-\gamma u_{\mathbb{B}}-\gamma\mathbb{C}x=x+\gamma u_{\mathbb{A}} \\
& \Rightarrow &\mathbb{J}_{\gamma \mathbb{A}}\circ(2\mathbb{J}_{\gamma \mathbb{B}}(z)-z-\gamma \mathbb{C}\circ \mathbb{J}_{\gamma \mathbb{B}}(z))=x\\
& \Rightarrow  &\mathbb{J}_{\gamma\mathbb{C}}\circ(\mathbb{J}_{\gamma \mathbb{A}}\circ(2\mathbb{J}_{\gamma \mathbb{B}}(z)-z-\gamma \mathbb{C}\circ \mathbb{J}_{\gamma \mathbb{B}}(z))+\gamma \mathbb{C}\circ \mathbb{J}_{\gamma \mathbb{B}}(z))=\mathbb{J}_{\gamma\mathbb{C}}(x+\gamma \mathbb{C}x)=x.
\end{eqnarray*}
With all the  identities above, we conclude that
\begin{eqnarray*}
Tz &=&  \mathbb{J}_{\gamma\mathbb{C}}\circ(\mathbb{J}_{\gamma \mathbb{A}}\circ(2\mathbb{J}_{\gamma \mathbb{B}}(z)-z-\gamma \mathbb{C}\circ \mathbb{J}_{\gamma \mathbb{B}}(z))+\gamma \mathbb{C}\circ \mathbb{J}_{\gamma \mathbb{B}}(z))+\gamma u_{\mathbb{B}}=x+\gamma u_{\mathbb{B}}=z.
\end{eqnarray*}   

Next, assume that $z=x+\gamma u_{\mathbb{B}}\in \text{Fix }T$, then  $x=\mathbb{J}_{\gamma\mathbb{B}}(z)\in \text{zer}(\mathbb{A}+\mathbb{B}+\mathbb{C})$ since
\begin{eqnarray*}
x+\gamma u_{\mathbb{B}}=Tz &=&  \mathbb{J}_{\gamma\mathbb{C}}\circ(\mathbb{J}_{\gamma \mathbb{A}}\circ(2\mathbb{J}_{\gamma \mathbb{B}}(z)-z-\gamma \mathbb{C}\circ \mathbb{J}_{\gamma \mathbb{B}}(z))+\gamma \mathbb{C}\circ \mathbb{J}_{\gamma \mathbb{B}}(z))+\gamma u_{\mathbb{B}}\\
\Rightarrow x&=&  \mathbb{J}_{\gamma\mathbb{C}}\circ(\mathbb{J}_{\gamma \mathbb{A}}\circ(2x-z-\gamma \mathbb{C}x)+\gamma \mathbb{C}x)\\
\Rightarrow x+\gamma \mathbb{C} x&=&   \mathbb{J}_{\gamma \mathbb{A}}\circ(2x-z-\gamma \mathbb{C}x)+\gamma \mathbb{C}x\\
\Rightarrow x+\gamma u_{\mathbb{A}}  &=&    x-\gamma u_{\mathbb{B}}-\gamma \mathbb{C}x \\
\Rightarrow   0 &=& u_{\mathbb{A}}+ u_{\mathbb{B}}+ \mathbb{C}x.
\end{eqnarray*}   
\end{proof}

 Next we discuss the convergence of the proposed Algorithm \ref{operator_T_algorithm1} under the assumption that $T$ is an averaged operator. The following are some standard convergence results under the assumption that $T$ is an averaged operator, following  Corollary 2.1 and Theorem 2.1 in \cite{davis2017three}.
 
\begin{thm}\label{convergence_theorem1}
Assume that $T:\mathcal{X}\to\mathcal{X}$ is $\alpha$-averaged with $\alpha=\frac{2\beta}{4\beta-\gamma}<1$ and suppose that $z^*\in \text{ Fix} \,T$. Consider a sequence of relaxation parameters  $(\lambda_j)_{j\geq0}\subset(0,\tfrac{1}{\alpha})$, where $\alpha=1/(2-\varepsilon)$ and $\alpha<2\beta/(4\beta-\gamma)$. Let $\sum_{j=0}^{\infty}\tau_j=\sum_{j=0}^{\infty}(1-\lambda_j/\alpha)\lambda_j/\alpha=\infty$. Let $z^0\in\mathcal{X}$ and $(z^j)_{j\geq0}\subseteq\mathcal{X}$ be the sequence generated by Algorithm \ref{operator_T_algorithm1}. Then, the following hold:
\begin{enumerate}[(i)]
    \item \label{convergence_theorem1a} The norm sequence $\lVert z^j-z^*\rVert_{j\geq 0}$ is monotonically nonincreasing for any $z^*\in\text{Fix }T$.
    \item \label{convergence_theorem1b} The residuals $\lVert Tz^j-z^j\rVert_{j\geq 0}$ form a monotonically nonincreasing sequence which converges to $0$.
    \item \label{convergence_theorem1c} The sequence $(z^j)_{j\geq0}$ converges weakly to a fixed point of $T$.
    \item Assume that $\underline{\tau}=\underset{j\geq 0}{\inf} \tau_j >0$ for $\tau>0$. Then, for any $z^*\in\text{Fix}\, T$ and for all $k\geq 0$, we have the following convergence rate:
    \begin{eqnarray*}
        \lVert Tz^k-z^k \rVert^2 \leq \frac{\lVert z^0-z^*\rVert^2}{\underline{\tau}(k+1)}\text{ and } \lVert Tz^k-z^k \rVert^2=o \left( \tfrac{1}{k+1}\right).
    \end{eqnarray*}
\end{enumerate}
\end{thm}
\begin{proof} 
 The proof for parts (i)-(iii) follows from Proposition 5.15 of \cite{bauschke2017correction}, and the proof of part (iv) follows from \cite{davis2016convergence}.  
\end{proof}

\begin{thm}\emph{(Convergence theorem)}. Let $T:\mathcal{X}\to\mathcal{X}$ be an $\alpha$-averaged operator with $\alpha=\frac{2\beta}{4\beta-\gamma}<1$. Let $(\lambda_j)_{j\geq0}\subseteq \left(0,\frac{1}{\alpha}\right)$ be a sequence of relaxation parameters, where $\alpha=1/(2-\varepsilon)<2\beta/(4\beta-\gamma)$. Suppose   $\sum_{j=0}^{\infty}\tau_j=\sum_{j=0}^{\infty}(1-\lambda_j/\alpha)\lambda_j/\alpha=\infty$ is satisfied. Let $z^0\in\mathcal{X}$ and $(z^j)_{j\geq0}\subseteq\mathcal{X}$ be the sequence generated by Algorithm \ref{operator_T_algorithm1}. Then, the following results hold.
\begin{enumerate}
    \item Assume $\inf_{j\geq0}\lambda_j>0$ and $z^*$ denote a weak limit point of $z^k$. Then:
    \begin{enumerate}[(a)]
        \item The sequence \label{convergence_thm_1a} $(\mathbb{C}x^j_{\mathbb{B}})_{j\geq0}$ converges strongly to $\mathbb{C}x^*$, where $x^*\in\text{zer}(\mathbb{A}+\mathbb{B}+\mathbb{C})$, 
    \item \label{convergence_thm_1b} The iterates $\mathbb{J}_{\gamma\mathbb{B}}(z^j)_{j\geq 0}$ converges weakly to $\mathbb{J}_{\gamma\mathbb{B}}(z^*)\in\text{zer}(\mathbb{A}+\mathbb{B}+\mathbb{C})$.
     \item \label{convergence_thm_1c} The sequence $(\mathbb{J}_{\gamma\mathbb{A}}\circ(2\mathbb{J}_{\gamma\mathbb{B}}-I-\gamma\mathbb{C}\circ\mathbb{J}_{\gamma\mathbb{B}})(z^j)_{j\geq 0}$ weakly converges to $\mathbb{J}_{\gamma\mathbb{B}}(z^*)\in\text{zer}(\mathbb{A}+\mathbb{B}+\mathbb{C})$.
    \item \label{convergence_thm_1d} the sequence $(\mathbb{J}_{\gamma\mathbb{C}}\circ(\mathbb{J}_{\gamma\mathbb{A}}\circ(2\mathbb{J}_{\gamma\mathbb{B}}-I-\gamma\mathbb{C}\circ\mathbb{J}_{\gamma\mathbb{B}})+\gamma\mathbb{C}\circ \mathbb{J}_{\gamma \mathbb{B}})(z^j)_{j\geq 0}$ weakly converges to $\mathbb{J}_{\gamma\mathbb{B}}(z^*)\in\text{zer}(\mathbb{A}+\mathbb{B}+\mathbb{C})$.
    \end{enumerate}
    \item The sequences $\mathbb{J}_{\gamma\mathbb{B}}(z^j)_{j\geq 0}$ and $(\mathbb{J}_{\gamma\mathbb{C}}\circ(\mathbb{J}_{\gamma\mathbb{A}}\circ(2\mathbb{J}_{\gamma\mathbb{B}}-I-\gamma\mathbb{C}\circ\mathbb{J}_{\gamma\mathbb{B}})+\gamma\mathbb{C}\circ \mathbb{J}_{\gamma \mathbb{B}})(z^j)_{j\geq 0}$ converges strongly to a solution in $\text{zer}(\mathbb{A}+\mathbb{B}+\mathbb{C})$ if at least one of the following holds:
    \begin{enumerate}[(a)]
        \item Operator $\mathbb{A}$ is uniformly monotone on every nonempty bounded subset of its domain.
        \item Operator $\mathbb{B}$ is uniformly monotone on every nonempty bounded subset of its domain.
        \item Operatoe $\mathbb{C}$ is demiregular at each solution point in $\text{zer}(\mathbb{A}+\mathbb{B}+\mathbb{C})$.
    \end{enumerate}
\end{enumerate}
    
\end{thm}
\begin{proof}
\begin{enumerate}
    \item \begin{enumerate}[(a)]
    \item 
Let $k\geq 0$. Then, using Corollary 2.14 in \cite{bauschke2017correction}, we observe that
\begin{eqnarray}\label{convergence_thm}
&&\lVert z^{k+1}-z^* \rVert^2= \lVert (1-\lambda_k)(z^{k}-z^*)+ \lambda_k(Tz^k-z^*)\rVert^2   \nonumber \\
&&=(1-\lambda_k)\lVert z^{k}-z^*\rVert^2+\lambda_k \lVert Tz^k-z^*\rVert^2-\lambda_k(1-\lambda_k)\lVert Tz^k-z^k\rVert^2.
\end{eqnarray}
In view of Theorem 2.1 of \cite{davis2017three}, we get
\begin{eqnarray*}
\sum_{i=k}^{\infty}   \lVert \mathbb{C} x^k_{\mathbb{B}}-\mathbb{C} \mathbb{J}_{\gamma\mathbb{B}}(z^*) \rVert^2 \leq \frac{\lVert z^k-z^* \rVert^2}{\gamma\lambda_k\left(2\beta-\tfrac{\gamma}{\varepsilon}\right)} , 
\end{eqnarray*}
which gives $\lVert \mathbb{C} x^k_{\mathbb{B}}-\mathbb{C} \mathbb{J}_{\gamma\mathbb{B}}(z^*) \rVert^2 \to 0$ as $k\to\infty$.
\item We recall the notations from Lemma \ref{weak_convergence_supporting_lem} given by
\begin{eqnarray*}
&x^k_{\mathbb{B}}=\mathbb{J}_{\gamma \mathbb{B}}(z^k) ,~ &u^k_{\mathbb{B}}=\gamma^{-1}(z^k-x^k_{\mathbb{B}})\in \mathbb{B}x^k_{\mathbb{B}}   \\
&x^k_{\mathbb{A}}=\mathbb{J}_{\gamma \mathbb{A}}(2x^k_{\mathbb{B}}-z^k-\gamma \mathbb{C} x^k_{\mathbb{B}}),~ & u^k_{\mathbb{A}}=\gamma^{-1}(2x^k_{\mathbb{B}}-z^k-\gamma \mathbb{C}x^k_{\mathbb{B}}-x^k_{\mathbb{A}})\in \mathbb{A}x^k_{\mathbb{A}}\\    
&x^k_{\mathbb{C}}=\mathbb{J}_{\gamma \mathbb{C}}(x^k_{\mathbb{A}}+\gamma \mathbb{C} x^k_{\mathbb{B}}),~& u^k_{\mathbb{C}}=\gamma^{-1}(x^k_{\mathbb{A}}+\gamma \mathbb{C}x^k_{\mathbb{B}}-x^k_{\mathbb{C}})\in \mathbb{C}x^k_{\mathbb{C}}.    
\end{eqnarray*}
Note that for all $k\geq0$, we have
\[
\lVert x^k_{\mathbb{B}}-\mathbb{J}_{\gamma\mathbb{B}}(z^*)\rVert=\lVert\mathbb{J}_{\gamma\mathbb{B}}(z^k)-\mathbb{J}_{\gamma\mathbb{B}}(z^*)\rVert\leq \lVert z^k-z^*\rVert\leq \lVert z^0-z^*\rVert,
\]
therefore $(x^j_{\mathbb{B}})_{j\geq 0}$ is bounded and admits a weak sequential limit labelled as $\bar{x}$.\\

Now, assume that there exists a subsequence $(k_j)_{j\geq 0}$ such that $x^{k_j}_{\mathbb{B}}\rightharpoonup\bar{x}$ as $j\to\infty$. Let $x^*\in \text{zer}(\mathbb{A}+\mathbb{B}+\mathbb{C})$. Next, observe that $\mathbb{C}$ is maximal monotone and $\mathbb{C}x^k_{\mathbb{B}}\to\mathbb{C}x^*$, and $x^{k_j}_{\mathbb{B}}\to\bar{x}$, thus, in view of Proposition 20.33(ii) of \cite{bauschke2017correction} and weak-to-strong sequential closeness of $\mathbb{C}$, we have 
\[
\mathbb{C}\bar{x}=\mathbb{C}x^*\text{ and } \mathbb{C} x^{k_j}_{\mathbb{B}}=\mathbb{C}\bar{x}.
\]
Further, in view of \eqref{convergence_theorem1b} of Theorem \ref{convergence_theorem1} and Lemma \ref{weak_convergence_supporting_lem}, we have $x^{k}_{\mathbb{C}}-x^{k}_{\mathbb{B}}=Tz^k-z^k\to0$ as $k\to\infty$. Thus, with $j\to\infty$, we obtain
\begin{eqnarray*} &&x^{k_j}_{\mathbb{B}}\rightharpoonup\bar{x},~x^{k_j}_{\mathbb{A}}\rightharpoonup\bar{x},~x^{k_j}_{\mathbb{C}}\rightharpoonup\bar{x},~\mathbb{C}x^{k_j}_{\mathbb{B}}\rightharpoonup\mathbb{C}\bar{x}\\
&\text{and}&u^{k_j}_{\mathbb{B}}\rightharpoonup\tfrac{1}{\gamma}(z^*-\bar{x}),~u^{k_j}_{\mathbb{A}}\rightharpoonup\tfrac{1}{\gamma}(\bar{x}-z^*-\gamma\mathbb{C}\bar{x}),~u^{k_j}_{\mathbb{C}}\rightharpoonup\tfrac{1}{\gamma}(\bar{x}+\gamma\mathbb{C}\bar{x}). 
\end{eqnarray*}
On applying Proposition, 25.5 of \cite{bauschke2017correction} to $(x^{k_j}_{\mathbb{A}},u^{k_j}_{\mathbb{A}})\in\text{gra}\textcolor{black}{(}\mathbb{A}\textcolor{black}{)},~(x^{k_j}_{\mathbb{B}},u^{k_j}_{\mathbb{B}})\in\mathbb{B}$, and $(x^{k_j}_{\mathbb{B}},\mathbb{C}x^{k_j}_{\mathbb{B}})\in\mathbb{C}$, we observe that $\bar{x}\in\text{zer}(\mathbb{A}+\mathbb{B}+\mathbb{C})$, $z^*-\bar{x}\in\gamma\mathbb{B}\bar{x}$, $\bar{x}-z^*-\gamma\mathbb{C}\bar{x}\in\gamma\mathbb{A}\bar{x}$, and $\bar{x}+\gamma\mathbb{C}\bar{x}\in\gamma\mathbb{C}\bar{x}$. Thus, we obtain $\bar{x}=\mathbb{J}_{\gamma\mathbb{B}}(z^*)$, and hence $\bar{x}$ is a unique weak sequential cluster point of $(x^j_{\mathbb{B}})_{j\geq0}$. Therefore, in view of Lemma 2.38 of \cite{bauschke2017correction}, we conclude that $(x^j_{\mathbb{B}})_{j\geq0}$ converges weakly to $\mathbb{J}_{\gamma\mathbb{B}}(z^*)$. 
\item On proceeding in similar manner to part \ref{convergence_thm_1c} and part \ref{convergence_thm_1d}, we implies that
\[
x^k_{\mathbb{C}}-x^k_{\mathbb{B}}=Tz^k-z^k \to 0 \text{ as }k\to\infty \text{ and }x^k_{\mathbb{B}}\rightharpoonup \mathbb{J}_{\gamma\mathbb{B}}(z^*) \implies x^k_{\mathbb{C}}\rightharpoonup \mathbb{J}_{\gamma\mathbb{B}}(z^*). 
\]
\end{enumerate}
\item The proofs are in a similar manner to those given in Theorem 2.1 of \cite{davis2017three}.
\end{enumerate}
\end{proof}

Before we proceed to the analysis of convex optimization problems of using Algorithm \ref{operator_T_algorithm1} under several assumptions on the regularity of the problem, we discuss the following lemma:

\begin{lem}\label{closed_ball}
Assume that $T:\mathcal{X}\to\mathcal{X}$ is $\alpha$-averaged with $\alpha=\frac{2\beta}{4\beta-\gamma}<1$.
Let $(z^j)_{j\geq 0}$ is generated by Algorithm \ref{operator_T_algorithm1} and $\gamma>0$. Let $z^*$ be a fixed point of $T$ and $x^*=\mathbb{J}_{\gamma\beta}(z^*)$. Then, $(x^j_{\mathbb{A}})_{j\geq 0}$ and $(x^j_{\mathbb{B}})_{j\geq 0}$ are contained within the closed ball $\overline{B(x^*,(1+\gamma/\beta)\lVert z^0-z^*\rVert})$.    
\end{lem}
\begin{proof}
By (\ref{convergence_theorem1a}) of Theorem \ref{convergence_theorem1}, we have
\begin{eqnarray*}
  \lVert x^k_{\mathbb{B}}-x^*\rVert=\lVert \mathbb{J}_{\gamma\mathbb{B}}(z^k)-\mathbb{J}_{\gamma\mathbb{B}}(z^*)\rVert \leq \lVert z^k-z^*\rVert \leq \lVert z^0-z^*\rVert.  
\end{eqnarray*}
Similarly, we have
\begin{eqnarray*}
\lVert x^k_{\mathbb{A}}-x^*\rVert&=&\lVert \mathbb{J}_{\gamma\mathbb{A}}({\mathbb R}_{\gamma\mathbb{B}}(z^k)-\gamma\mathbb{C}x^k_{\mathbb{B}})-\mathbb{J}_{\gamma\mathbb{A}}({\mathbb R}_{\gamma\mathbb{B}}(z^*)-\gamma\mathbb{C}x^*)\rVert\\
&\leq& \lVert {\mathbb R}_{\gamma\mathbb{B}}(z^k)-{\mathbb R}_{\gamma\mathbb{B}}(z^*)+\gamma\mathbb{C}x^*-\gamma\mathbb{C}x^k_{\mathbb{B}}\rVert\\
&\leq& \lVert z^k-z^*\rVert+\frac{\gamma}{\beta}\lVert z^k-z^*\rVert \leq\left(1+\frac{\gamma}{\beta}\right)\lVert z^0-z^*\rVert.
\end{eqnarray*}
With the inequality above, we also have
\begin{eqnarray*}
\lVert x^k_{\mathbb{C}}-x^*\rVert&=&\lVert \mathbb{J}_{\gamma\mathbb{C}}(x^k_{\mathbb{A}}+\gamma\mathbb{C}x^k_{\mathbb{B}})-\mathbb{J}_{\gamma\mathbb{C}}(x^*+\gamma\mathbb{C}x^*)\rVert
=\lVert x^k_{\mathbb{A}}-x^*+\gamma\mathbb{C}x^k_{\mathbb{B}}-\gamma\mathbb{C}x^*\rVert\\
&\leq& \lVert z^k-z^*\rVert+\frac{2\gamma}{\beta}\lVert z^k-z^*\rVert \leq\left(1+\frac{2\gamma}{\beta}\right)\lVert z^0-z^*\rVert.    
\end{eqnarray*}
\end{proof}



\section{Application to composite convex minimization problems}\label{section_multiblock}
\label{sec:admm}

\textcolor{black}{In this section, we consider the application of the proposed three-operator splitting for composite convex minimization problems. We first review the classical three-block ADMM method, then  in Section \ref{sec:admm-derivation} we derive its dual form. We will also discuss the convergence of the proposed new splitting scheme.}

\subsection{The classical three-block ADMM}

We  analyze the convex optimization problems under several assumptions on the regularity of the problem in this section:
\begin{enumerate}[$\bullet$]
    \item Every considered function is proper, closed, and convex.
    \item Every differentiable function is Fr\'echet differentiable.
    \item The functions $f_i:\mathcal{X}\to(-\infty,+\infty]$, $i=1,2,3$, satisfy the existence of solution condition, i.e.,
    \[
    \text{zer}(\partial f_1+\partial f_2+\partial f_3)\neq \emptyset.
    \]
\end{enumerate}
\textcolor{black}{ The following are some well known facts. }
\begin{prop}\emph{(Optimality conditions of prox).}
Let $w\in\mathcal{X}$ and $f$ be a proper, closed, and convex function. Then, the following identity holds.
\[
\tilde{w}=\operatorname{prox}_{\gamma f}(w) \text{ if and only if }\frac{1}{\gamma}(w-\tilde{w})\in \partial f(\tilde{w}). 
\]
\end{prop}

\begin{prop}\emph{(Firm nonexpansiveness of prox).} Let $w,r\in\mathcal{X}$, and let $\tilde{w}=\operatorname{prox}_{\gamma f}(w)$ and $\tilde{r}=\operatorname{prox}_{\gamma f}(r)$. Then,
\[
\lVert\tilde{w}-\tilde{r}\rVert^2\leq \langle \tilde{w}-\tilde{r},w-r \rangle.
\]
In particular, $\operatorname{prox}_{\gamma f}$ is nonexpansive.    
\end{prop}

\begin{thm}\emph{(Descent Lemma).} Let $f$ be a differentiable function and $\nabla f$ is $\tfrac{1}{\gamma}$-Lipschitz. Then, for every $x,y\in\mathcal{X}$, we have
\[
f(x)\leq f(y)+\langle x-y,\nabla f(y)\rangle+\tfrac{1}{2\gamma}\lVert x-y \rVert^2.
\]  
\end{thm}

Consider a convex minimization problem with linear constraints and a separable objective function given by :
\begin{eqnarray}
\begin{rcases}\label{optimal_eq}
\min& f_1(x_1)+f_2(x_2)+f_3(x_3),\\
\text{subject to }&A_1x_1+A_2x_2+A_3x_3=b,\\
&x_1\in\textcolor{black}{\mathcal X_1=\mathbb{R}^{n_1}},~x_2\in\textcolor{black}{\mathcal X_2=\mathbb{R}^{n_2}},~x_3\in\textcolor{black}{\mathcal X_3=\mathbb{R}^{n_3}}
\end{rcases},
\end{eqnarray}
where $f_i:\mathbb{R}^{n_i}\to\mathbb{R}\cup\{+\infty\}$ are proper, closed, and convex functions (not necessarily smooth), $A_i\in\mathbb{R}^{m\times n_i}$ and $b\in\mathbb{R}^m$. Let $f^*_i$ denote the convex conjugate of function $f_i$, and let
\[
d_1(w)=f_1^*(A_1^\top w),~d_2(w)=f_2^*(A_2^\top w),~d_3(w)=f_3^*(A_3^\top w)-\langle w,b \rangle.
\]
The dual problem of \eqref{optimal_eq} is given by
\begin{eqnarray}
\label{dualproblemtoADMM}
    \underset{w\in \mathbb R^m}{\min} ~d_1(w)+d_2(w)+d_3(w)
\end{eqnarray}


\textcolor{black}{The original two-block ADMM method in \cite{glowinski1975approximation} has two parameters  $\gamma$ and $\sigma$, and it converges with certain parameter constraints for convex problems. If using only one parameter, i.e., taking  $\sigma=\gamma$, the classical two-block ADMM method in \cite{glowinski1975approximation}
is equivalent to the Douglas-Rachford splitting, thus such a single parameter two-block ADMM converges with any step size $\gamma>0$ for   convex problems, implied by the averagedness of the DRS operator. One could consider multiple-block ADMM methods with multiple parameters, but we only consider ADMM methods using a single step size here to compare its dual form to the proposed splitting scheme in this paper. } 
 A direct extension of the original ADMM in \cite{glowinski1975approximation} \textcolor{black}{with only one step size} for (\ref{optimal_eq}) is given by
\begin{subequations}
\label{ADMM-3block}
    \begin{eqnarray}\label{optimality_seq}
&&x_1^{k+1}=\underset{x_1\in\mathcal{X}_1}{\operatorname{argmin}}\{f_1(x_1)+\tfrac{\gamma}{2}\lVert (A_1x_1+A_2x_2^k+A_3x_3^k-b)-\tfrac{1}{\gamma}w^k\rVert^2\label{seq_1}  \\
&&x_2^{k+1}=\underset{x_2\in\mathcal{X}_2}{\operatorname{argmin}}\{f_2(x_2)+\tfrac{\gamma}{2}\lVert (A_1x_1^{k+1}+A_2x_2+A_3x_3^k-b)-\tfrac{1}{\gamma}w^k\rVert^2\}\label{seq_2}\\
&&x_3^{k+1}=\underset{x_3\in\mathcal{X}_3}{\operatorname{argmin}}\{f_3(x_3)+\tfrac{\gamma}{2}\lVert (A_1x_1^{k+1}+A_2x_2^{k+1}+A_3x_3-b)-\tfrac{1}{\gamma}w^k\rVert^2\} \label{seq_3}\\
&&w^{k+1}=w^k-\gamma(A_1x_1^{k+1}+A_2x_2^{k+1}+A_3x_3^{k+1}-b)\label{seq_4},
\end{eqnarray}
\end{subequations}
where $w^k\in\mathbb{R}^m$ is the Lagrange multiplier, $\gamma>0$ is the penalty parameter or step size. \\
\\

\subsection{The similarity compared with the dual formulation of the three-block ADMM method}
\label{sec:admm-derivation}

\textcolor{black}{In this subsection, we derive the dual formulation of \eqref{ADMM-3block}, which will be a scheme solving the dual problem \eqref{dualproblemtoADMM}. 
For convenience, define $\bar{x}^k:=x^{k+1}$.  We first list the change of variables needed in this subsection:
\begin{eqnarray}
    z^k&=&w^k-\gamma (A_1x_1^{k+1}+A_2 x_2^{k+1}) \label{definition-z} \\
    p^{k+1}&=&w^{k+1}-\gamma (A_1\bar{x}_1^{k+1}+A_2 x_2^{k+1}+A_3x_3^{k+1}\textcolor{black}{-b})\label{definition-p}\\
    v^{k+1}&=&w^{k+1}-\gamma (A_1\bar{x}_1^{k+1}+A_2 \bar{x}_2^{k+1}+A_3x_3^{k+1}-b) \label{definition-v}
\end{eqnarray}  
}

We start the derivation with \eqref{seq_3} as follows 
\begin{eqnarray}\label{operator_splitting_convex_1}
&&x_3^{k+1}=\underset{x_3\in\mathcal{X}_3}{\operatorname{argmin}}\{f_3(x_3)+\tfrac{\gamma}{2}\lVert (A_1x_1^{k+1}+A_2x_2^{k+1}+A_3x_3-b)-\tfrac{1}{\gamma}w^k\rVert^2\} \nonumber\\
&\iff& 0\in \partial f_3(x_3^{k+1})- A_3^\top (w^k-\gamma (A_1x_1^{k+1}+A_2 x_2^{k+1}+A_3x_3^{k+1}-b))\nonumber\\
&\iff& 0\in \partial f_3(x_3^{k+1})- A_3^\top w^{k+1} \quad \text{ [from \eqref{seq_4}]}\nonumber\\
&\iff&A_3^\top w^{k+1}\in \partial f_3(x_3^{k+1})\nonumber\\
&\iff& x_3^{k+1}\in \partial f_3^*(A_3^\top w^{k+1}).
\end{eqnarray}
In view of \eqref{seq_4}, observe that
\begin{eqnarray}
&& w^{k+1}=w^k-\gamma (A_1x_1^{k+1}+A_2 x_2^{k+1}+A_3x_3^{k+1}-b)\nonumber\\
&& w^{k+1}=z^k-\gamma (A_3x_3^{k+1}-b) \quad \textcolor{black}{\text{[by definition of $z^k$ in \eqref{definition-z}}]}\label{admm_extra_eq1}\\
&\iff& w^{k+1}= z^k-\gamma (A_3\partial f_3^*(A_3^\top w^{k+1})-b)\quad \text{ [from \eqref{operator_splitting_convex_1}]}\nonumber\\
&\iff& w^{k+1}=z^k-\gamma \partial d_3(w^{k+1})\label{admm_eq2}\\
&\iff&(I+\gamma \partial d_3)(w^{k+1})=z^k\nonumber\\
&\iff& w^{k+1}=\operatorname{prox}_{\gamma d_3}(z^k) \label{admm_eq1}.
\end{eqnarray}
Consider the equation given by \eqref{seq_1}, we have
\begin{eqnarray}\label{operator_splitting_convex_2}
&&x_1^{k+1}=\underset{x_1\in\mathcal{X}_1}{\operatorname{argmin}}\{f_1(x_1)+\tfrac{\gamma}{2}\lVert (A_1x_1+A_2x_2^{k}+A_3x_3^{k}-b)-\tfrac{1}{\gamma}w^k\rVert^2\} \nonumber\\
&\iff& 0\in \partial f_1(x_1^{k+1})- A_1^\top (w^k-\gamma (A_1x_1^{k+1}+A_2 x_2^{k}+A_3x_3^{k}-b))\nonumber\\
&\iff& 0\in \partial f_1(\bar{x}_1^{k})- A_1^\top (w^k-\gamma (A_1\bar{x}_1^{k}+A_2 x_2^{k}+A_3x_3^{k}-b))\nonumber\\
&\iff& 0\in \partial f_1(\bar{x}_1^{k+1})- A_1^\top (w^{k+1}-\gamma (A_1\bar{x}_1^{k+1}+A_2 x_2^{k+1}+A_3x_3^{k+1}-b))\text{ holds for any }k\in\mathbb{N}\cup\{0\}\nonumber\\
&\iff& 0\in \partial f_1(\bar{x}_1^{k+1})- A_1^\top p^{k+1} \quad \textcolor{black}{\text{[by definition of $p^{k+1}$ in \eqref{definition-p}}]}\nonumber\\
&\iff&A_1^\top p^{k+1}\in \partial f_1(\bar{x}_1^{k+1})\nonumber\\
&\iff& \bar{x}_1^{k+1}\in \partial f_1^*(A_1^\top p^{k+1}).  
\end{eqnarray}

Next, from \eqref{seq_2}, we have
\begin{eqnarray}\label{operator_splitting_convex_3}
&&x_2^{k+1}=\underset{x_2\in\mathcal{X}_2}{\operatorname{argmin}}\{f_2(x_2)+\tfrac{\gamma}{2}\lVert (A_1x_1^{k+1}+A_2x_2+A_3x_3^{k}-b)-\tfrac{1}{\gamma}w^k\rVert^2\} \nonumber\\
&\iff& 0\in \partial f_2(x_2^{k+1})- A_2^\top (w^k-\gamma (A_1x_1^{k+1}+A_2 x_2^{k+1}+A_3x_3^{k}-b))\nonumber\\
&\iff& 0\in \partial f_2(\bar{x}_2^{k})- A_2^\top (w^k-\gamma (A_1\bar{x}_1^{k}+A_2 \bar{x}_2^{k}+A_3x_3^{k}-b))\nonumber\\
&\iff& 0\in \partial f_2(\bar{x}_2^{k+1})- A_2^\top (w^{k+1}-\gamma (A_1\bar{x}_1^{k+1}+A_2 \bar{x}_2^{k+1}+A_3x_3^{k+1}-b))\text{ holds for any }k\in\mathbb{N}\cup\{0\}\nonumber\\
&\iff& 0\in \partial f_2(\bar{x}_2^{k+1})- A_2^\top v^{k+1}  \quad \textcolor{black}{\text{[by definition of $v^{k+1}$ in \eqref{definition-v}}]} \nonumber\\
&\iff&A_2^\top v^{k+1}\in \partial f_2(\bar{x}_2^{k+1})\nonumber\\
&\iff& \bar{x}_2^{k+1}\in \partial f_2^*(A_2^\top v^{k+1}).  
\end{eqnarray}
Now, from \eqref{operator_splitting_convex_1}, \eqref{operator_splitting_convex_2}, and \eqref{operator_splitting_convex_3}, observe that 
\begin{eqnarray}
&&v^{k+1}=w^{k+1}-\gamma (A_1 \bar{x}_1^{k+1}+A_2\bar{x}_2^{k+1}+A_3x_3^{k+1}-b)\notag \\
&\iff&v^{k+1}=w^{k+1}-\gamma (A_1 \partial f_1^*(A_1^\top p^{k+1})+A_2\partial f_2^*(A_2^\top v^{k+1})+A_3\partial f_3^*(A_3^\top w^{k+1})-b)\notag \\
&\iff&v^{k+1}=w^{k+1}-\gamma \partial d_1(p^{k+1})-\gamma\partial d_2(v^{k+1})-\gamma\partial d_3(w^{k+1}) \label{v_k-formula}
\end{eqnarray}
 
\textcolor{black}{
Note that,
\begin{eqnarray}
&&p^{k+1}=w^{k+1}-\gamma (A_1 \bar{x}_1^{k+1}+A_2x_2^{k+1}+A_3x_3^{k+1}-b)\nonumber\\
&\iff&p^{k+1}=w^{k+1}-\gamma (A_1 \partial f_1^*(A_1^\top p^{k+1})+A_2\partial f_2^*(A_2^\top v^{k})+A_3\partial f_3^*(A_3^\top w^{k+1})-b)\nonumber\\
&\iff&p^{k+1}=w^{k+1}-\gamma \partial d_1(p^{k+1})-\gamma\partial d_2(v^{k})-\gamma\partial d_3(w^{k+1})\label{admm_extra_eq3}\\
&\iff& (I+\gamma\partial d_1)(p^{k+1})=w^{k+1}-\gamma\partial d_2(v^{k})-\gamma\partial d_3(w^{k+1})\label{admm_eq3}\\
&\iff& (I+\gamma\partial d_1)(p^{k+1})=2w^{k+1}-z^k-\gamma\partial d_2(v^{k})\quad \text{ [from \eqref{admm_eq2}]}\nonumber\\
&\iff&p^{k+1}=\operatorname{prox}_{\gamma d_1}(2w^{k+1}-z^k-\gamma\partial d_2(v^{k}))\nonumber.
\end{eqnarray}
Now by \eqref{v_k-formula}, we have
\begin{eqnarray*}
&&v^{k+1}=w^{k+1}-\gamma \partial d_1(p^{k+1})-\gamma\partial d_2(v^{k+1})-\gamma\partial d_3(w^{k+1})\\
&\iff& (I+\gamma\partial d_2)(v^{k+1})=w^{k+1}-\gamma\partial d_1(p^{k+1})-\gamma\partial d_3(w^{k+1}) \\
    &\iff& (I+\gamma\partial d_2)(v^{k+1})=p^{k+1}+\gamma\partial d_2(v^{k})\quad \text{ [from \eqref{admm_extra_eq3}]}\\
&\iff&v^{k+1}=\operatorname{prox}_{\gamma d_2}(p^{k+1}+\gamma\partial d_2(v^{k})).
\end{eqnarray*}
Last, we have 
\begin{eqnarray*} 
  && z^k=w^k-\gamma(A_1x_1^{k+1}+A_2x_2^{k+1}) \quad\text{[by definition of $z^{k}$ in \eqref{definition-z}]} \notag \\
&\implies& z^{k+1}=w^{k+1}-\gamma(A_1\bar{x}_1^{k+1}+A_2\bar{x}_2^{k+1}) \\
&\implies& z^{k+1}=z^{k}-\gamma(A_1\bar{x}_1^{k+1}+A_2\bar{x}_2^{k+1}+A_3 x_3^{k+1}-b)
   \quad \text{[by \eqref{admm_extra_eq1}]} \\
   &\implies& z^{k+1}=z^k+v^{k+1}-w^{k+1} \quad\text{[by definition of $v^{k+1}$ in \eqref{definition-v}]}.
\end{eqnarray*}
}
Hence, we derive \textcolor{black}{the dual formulation of the three-block ADMM} as
\begin{eqnarray}\label{dual-ADMM}
\begin{rcases}
& w^{k+1}=\operatorname{prox}_{\gamma d_3}(z^k)\\
& p^{k+1}= \operatorname{prox}_{\gamma d_1} (2w^{k+1}-z^k-\textcolor{black}{\gamma\partial d_2(v^{k})})\\
&v^{k+1}=\operatorname{prox}_{\gamma d_2} (p^{k+1}+\textcolor{black}{\gamma\partial d_2(v^{k})})\\
&z^{k+1}=z^k+v^{k+1}-w^{k+1}
\end{rcases},
\end{eqnarray}
\textcolor{black}{
where we have abused notation to use $\partial d_2(v^{k})$ to denote any element in the subdifferential set $\partial d_2(v^{k})$.
As a comparison,
the proposed new splitting in the previous section for minimizing $d_1(w)+d_2(w)+d_3(w)$ is written as
\begin{eqnarray}\label{FPI_algorithm_convex_case}
\begin{rcases}
& w^{k+1}=\operatorname{prox}_{\gamma d_3}(z^k)\\
& p^{k+1}= \operatorname{prox}_{\gamma d_1} (2w^{k+1}-z^k- {\gamma\partial d_2(w^{k+1})})\\
&v^{k+1}=\operatorname{prox}_{\gamma d_2} (p^{k+1}+ {\gamma\partial d_2(w^{k+1})})\\
&z^{k+1}=z^k+v^{k+1}-w^{k+1}
\end{rcases}.
\end{eqnarray}
We can observe that the dual formulation to three-block ADMM \eqref{dual-ADMM} is similar to but different from the proposed scheme \eqref{FPI_algorithm_convex_case}.  
}
 
The splitting scheme \eqref{FPI_algorithm_convex_case} can be implemented as in Algorithm \ref{FPI_algorithm_modified_convex_case} for solving \eqref{dualproblemtoADMM} where $d_i$ are proper closed convex functions and $\nabla d_2$ is $L$-Lipshcitz continuous.
\begin{algorithm}[!htb]
\caption{Initialize $z^0\in\mathcal{X},~\gamma\in(0,2/L)$. For $k=0,1,2,\ldots$}\label{FPI_algorithm_modified_convex_case}
\begin{enumerate}
\item Compute $w^{k+1}=\operatorname{prox}_{\gamma d_3}(z^k)$;
\item Compute $p^{k+1}= \operatorname{prox}_{\gamma d_1} (2w^{k+1}-z^k-\gamma \nabla d_2(w^{k+1}))$;
\item Compute $v^{k+1}=\operatorname{prox}_{\gamma d_2} (p^{k+1}+\gamma\nabla d_2(w^{k+1}))$;
\item Update $z^{k+1}=z^k+(v^{k+1}-w^{k+1}).$
\end{enumerate}
\end{algorithm}

\subsection{The averagedness of the splitting operator for a special case}
\label{sec:specialcase}
\textcolor{black}{
In this subsection, in addition to the basic convexity assumption of functions $d_i$ for solving \eqref{dualproblemtoADMM},
we show that the proposed splitting operator $T$ corresponding to iteration \eqref{FPI_algorithm_convex_case} is $\frac{1}{2}$-averaged under the following special  assumption: 
\begin{assumption}
\label{assump}
Assume  $\textcolor{black}{\mathcal{X}}=X \bigoplus Y$ where $X,Y$ are orthogonal subspaces,
the domain of $d_1(w)$ is inside $X$ and the domain of $d_2(w)$ is inside  $Y$. 
\end{assumption}
\begin{rem}
\label{rem-convergence}
 In the literature,  in addition to basic convexity assumptions of $f_i$, the convergence of  the 
  three-block  ADMM iteration \eqref{ADMM-3block}  (or equivalently \eqref{dual-ADMM}) has been well studied. In \cite{cai2017convergence}, the convergence of   \eqref{ADMM-3block} to a saddle point can be proven if assuming  strong convexity of $f_3$. In \cite{chen2016directNO},   a sufficient condition for  the convergence of the 
  3-block ADMM iteration \eqref{ADMM-3block} is $A_1^\top A_2=0$, or $A_2^\top A_3=0$, or $A_1^\top A_3=0$. 
  For $d_1(w)=f_1^*(A^\top_1 w)$ and $d_2(w)=f_2^*(A^\top_2 w)$, $A^\top_1A_2=0$ implies Assumption \ref{assump}, thus the assumption we consider here is weaker than
    $A^\top_1A_2=0$.
\end{rem}
}

\textcolor{black}{
\begin{lem}
\label{lemma:specialcase}
Under Assumption \ref{assump}, if using the same subderivative $\textcolor{black}{\tilde\nabla d_2(w^{k+1})\in\partial d_2(w^{k+1})}$ in the second line and third line in \eqref{FPI_algorithm_convex_case}, the iteration 
   \eqref{FPI_algorithm_convex_case}  is equivalent to
   \begin{eqnarray} 
z^{k+1} &=&  z^k + \operatorname{prox}_{\gamma d_2}(\operatorname{prox}_{\gamma d_1}(2\operatorname{prox}_{\gamma d_3}(z^k) -z^{k}))-\operatorname{prox}_{\gamma d_3}(z^k). \label{eq:simplified}
\end{eqnarray}
\end{lem}
\begin{proof}
Without loss of generality, we consider the following special case of 
the domains of $d_1$ and $d_2$ being orthogonal complement to each other.
Let any $w\in \textcolor{black}{\mathcal{X}}$
 be written as $w=(x,y)\in X \bigoplus Y$ 
with $x\in X$ and $y\in Y$, 
and assume $d_1(x, y) = d_1(x, 0)$, $d_2(x, y) = d_2(0, y)$. Then we have $\textcolor{black}{\tilde\nabla_w} d_2(x, y) = \textcolor{black}{\tilde\nabla_y} d_2(0, y)$ for all $x\in X$, $y\in Y$.  
Let $x_w$, $y_w$ be the $X$ and $Y$ components of $w$. Let $x_z$, $y_z$ be the $X$ and $Y$ components of $z$. Let $\partial_x$, $\partial_y$ be the $x$ subdifferential and $y$ subdifferential, respectively. The second line
of \eqref{FPI_algorithm_convex_case} can be written as\\
\begin{eqnarray*}
&&p^{k+1} = \underset{h=(x,y)}{\operatorname{argmin}} \{ d_1(h) + \frac{1}{2\gamma} \lVert h - (2w^{k+1} - z^k - \gamma \textcolor{black}{\tilde\nabla} d_2(w^{k+1})) \rVert^2 \}\\
=&& \underset{(x, y) }{\operatorname{argmin}} \{ d_1(x, 0) + \frac{1}{2\gamma} \lVert (x, y) - [ 2(x_w^{k+1}, y_w^{k+1}) - (x_z^k, y_z^k) - \gamma ( 0, \textcolor{black}{\tilde\nabla}_y d_2(0, y_w^{k+1}) ) ] \rVert^2 \}\\ 
=&& \underset{x}{\operatorname{argmin}} \{ d_1(x, 0) + \frac{1}{2\gamma} \lVert (x, 0) - ( 2(x_w^{k+1}, 0) - (x_z^k, 0) ) \rVert^2 \}\\
&&\bigoplus \underset{y}{\operatorname{argmin}} \{ \frac{1}{2\gamma} \lVert (0, y) - [ 2(x_w^{k+1}, y_w^{k+1}) - (0, y_z^k) - \gamma ( 0, \textcolor{black}{\tilde\nabla}_y d_2(0, y_w^{k+1}) ) ] \rVert^2 \}\\
=&& \underset{x \in X}{\operatorname{argmin}} \{ d_1(x, 0) + \frac{1}{2\gamma} \lVert (x, 0) - ( 2(x_w^{k+1}, 0) - (x_z^k, 0) ) \rVert^2 \}\\
&&\bigoplus \underset{y \in Y}{\operatorname{argmin}} \{\frac{1}{2\gamma} \lVert (0, y) - ( 2(x_w^{k+1}, y_w^{k+1}) - (0, y_z^k)) \rVert^2 \} - \gamma ( 0, \textcolor{black}{\tilde\nabla}_y d_2(0, y_w^{k+1}) ) \\
=&& \underset{(x, y) }{\operatorname{argmin}} \{ d_1(x, y) + \frac{1}{2\gamma} \lVert (x, y) - ( 2(x_w^{k+1}, y_w^{k+1}) - (x_z^k, y_z^k)) \rVert^2 \} - \gamma ( 0, \textcolor{black}{\tilde\nabla}_y d_2(0, y_w^{k+1}) ) \\
=&& \operatorname{prox}_{\gamma d_1}\left(2w^{k+1} - z^k \right) - \gamma \textcolor{black}{\tilde\nabla} d_2(w^{k+1}).
\end{eqnarray*}
Thus under Assumption \ref{assump}, if $\textcolor{black}{\tilde\nabla} d_2(w^{k+1})$ in the second line and $\textcolor{black}{\tilde\nabla} d_2(w^{k+1})$ in the third line of \eqref{FPI_algorithm_convex_case} are taken to be the same subderivative, then \eqref{FPI_algorithm_convex_case}  becomes:
\begin{eqnarray*}
w^{k+1} &=& \operatorname{prox}_{\gamma d_3}(z^k) \notag \\
p^{k+1} &=& \operatorname{prox}_{\gamma d_1}\left(2w^{k+1} - z^k \right) - \gamma \partial d_2(w^{k+1})   \\
v^{k+1} &=& \operatorname{prox}_{\gamma d_2}\left(p^{k+1} + \gamma \partial d_2(w^{k+1})\right) = \operatorname{prox}_{\gamma d_2}\left( \operatorname{prox}_{\gamma d_1}\left(2w^{k+1} - z^k \right)\right) \notag \\
z^{k+1} &=& z^k + v^{k+1} - w^{k+1} = z^k + \operatorname{prox}_{\gamma d_2}(\operatorname{prox}_{\gamma d_1}(2w^{k+1}-z^{k}))-w^{k+1},\notag
\end{eqnarray*}
which simplifies to \eqref{eq:simplified} and the proof is concluded.
\end{proof}
}

\textcolor{black}{
In general, the composition of two proximal operators of convex functions is $2/3$-averaged and is not  $1/2$-averaged (i.e., firmly nonexpansive), see \cite{bauschke2017correction}.
\begin{lem}
\label{lemma:FNE}
Under Assumption \ref{assump}, the composition $\operatorname{prox}_{\gamma d_2} \circ \operatorname{prox}_{\gamma d_1}$ is firmly nonexpansive.
\end{lem}
\begin{proof}
Without loss of generality, we consider the following special case of 
the domains of $d_1$ and $d_2$ being orthogonal complement to each other.
Let any $w\in \textcolor{black}{\mathcal{X}}$
 be written as $w=(x,y)\in X \bigoplus Y$ 
with $x\in X$ and $y\in Y$, 
and assume $d_1(x, y) = d_1(x, 0)$, $d_2(x, y) = d_2(0, y)$. 
Let $w\in\textcolor{black}{\mathcal{X}}$ be given and $x_w$, $y_w$ be the $X$ and $Y$ components of $w$, then
\begin{eqnarray*}
&&\operatorname{prox}_{\gamma d_2} \circ \operatorname{prox}_{\gamma d_1}(w) = \operatorname{prox}_{\gamma d_2}(\underset{h=(x,y)}{\operatorname{argmin}} \{ d_1(h) + \frac{1}{2\gamma} \lVert h - w \rVert^2 \})\\
=&& \operatorname{prox}_{\gamma d_2}( \underset{x }{\operatorname{argmin}} \{ d_1(x, 0) + \frac{1}{2\gamma} \lVert (x, 0) - (x_w, 0) \rVert^2 \}, \underset{y}{\operatorname{argmin}} \{\frac{1}{2\gamma} \lVert (0, y) - (0, y_w) \rVert^2 \})\\
=&& \operatorname{prox}_{\gamma d_2}( \underset{x }{\operatorname{argmin}} \{ d_1(x, 0) + \frac{1}{2\gamma} \lVert (x, 0) - (x_w, 0) \rVert^2 \}, y_w)\\
=&& \underset{\hat{h} = (\hat{x},\hat{y})}{\operatorname{argmin}} \{d_2(\hat{h}) + \frac{1}{2\gamma}\lVert\hat{h} - ( \underset{x }{\operatorname{argmin}} \{ d_1(x, 0) + \frac{1}{2\gamma} \lVert (x, 0) - (x_w, 0) \rVert^2 \}, y_w) \rVert^2\}\\
=&& \underset{\hat{x}}{\operatorname{argmin}} \{d_2(0,\hat{y}) + \frac{1}{2\gamma}\lVert(\hat{x},0) - ( \underset{x }{\operatorname{argmin}} \{ d_1(x, 0) + \frac{1}{2\gamma} \lVert (x, 0) - (x_w, 0) \rVert^2 \}, 0) \rVert^2\}\\
&&\bigoplus \underset{\hat{y}}{\operatorname{argmin}} \{d_2(0,\hat{y}) + \frac{1}{2\gamma}\lVert(0,\hat{y}) - (0, y_w) \rVert^2\}\\
=&& \underset{\hat{x}}{\operatorname{argmin}} \{\frac{1}{2\gamma}\lVert(\hat{x},0) - ( \underset{x }{\operatorname{argmin}} \{ d_1(x, 0) + \frac{1}{2\gamma} \lVert (x, 0) - (x_w, 0) \rVert^2 \}, 0) \rVert^2\}\\
&&\bigoplus \underset{\hat{y}}{\operatorname{argmin}} \{d_2(0,\hat{y}) + \frac{1}{2\gamma}\lVert(0,\hat{y}) - (0, y_w) \rVert^2\}\end{eqnarray*}
\begin{eqnarray*}
=&& \underset{x}{\operatorname{argmin}} \{ d_1(x, y) + \frac{1}{2\gamma} \lVert (x, y) - (x_w, y_w) \rVert^2 \}\bigoplus \underset{y}{\operatorname{argmin}} \{d_2(x,y) + \frac{1}{2\gamma}\lVert(x,y) - (x_w, y_w) \rVert^2\}\\
=&& \underset{h=(x,y)}{\operatorname{argmin}}\{(d_1 + d_2)(x,y) + \frac{1}{2\gamma} \lVert (x, y) -(x_w,y_w)\rVert^2\}= \operatorname{prox}_{\gamma (d_1+d_2)}(w),
\end{eqnarray*}
which is firmly nonexpansive since it is the proximal operator of a convex closed proper function $d_1+d_2$. 
\end{proof}
}

\textcolor{black}{
\begin{lem}[Lemma 2.3 in \cite{davis2017three}]
\label{lemma:DY}
Let $\mathcal{X}$ be a Hilbert space.
Let $S := U + T_1 \circ V$, where $U$, $T_1 \colon \mathcal{X} \to \mathcal{X}$ are both firmly nonexpansive and $V \colon \mathcal{X} \to \mathcal{X}$. Let $W = I - (2U + V)$. Then we have for all $z,w \in \mathcal{X}$:
\begin{equation}
\lVert Sz - Sw \rVert^2 \leq \lVert z - w \rVert^2 - \lVert (I - S)z - (I - S)w \rVert^2 - 2\langle{T_1 \circ V z - T_1 \circ V w,Wz - Ww}\rangle. \label{eq:lemmaDY}
\end{equation}
\end{lem}
}

\textcolor{black}{
\begin{thm}
\label{thm-specialcase}
Under Assumption \ref{assump} and the assumption that $d_1, d_2, d_3$ are closed convex proper functions,
if using the same subderivative $\textcolor{black}{\tilde\nabla} d_2(w^{k+1})$ in the second line and third line in \eqref{FPI_algorithm_convex_case},
the splitting operator for the iteration \eqref{FPI_algorithm_convex_case} is $\frac12$-averaged thus \eqref{FPI_algorithm_convex_case} converges to the minimizer with any step size $\gamma>0$. 
\end{thm}
\begin{proof}
    Notice that \eqref{eq:simplified} derived in Lemma \ref{lemma:specialcase}
corresponds to the following splitting operator:
\[
T = \textcolor{black}{\mathbb{J}}_{\gamma \mathbb C} \circ \left( \textcolor{black}{\mathbb{J}}_{\gamma \mathbb A} \circ \left( 2\textcolor{black}{\mathbb{J}}_{\gamma \mathbb B} -  I \right) \right) + \left(   I - \textcolor{black}{\mathbb{J}}_{\gamma \mathbb B} \right),
\]
where $ \mathbb A=\partial d_1$, $ \mathbb B=\partial d_2$ and $ \mathbb C=\partial d_3$.
For any closed convex proper function $f$, its proximal operator $\operatorname{prox}_{\gamma f}$ is firmly nonexpansive, and so is $I-\operatorname{prox}_{\gamma f}$. 
Lemma \ref{lemma:FNE} implies $\textcolor{black}{\mathbb{J}}_{\gamma \mathbb C} \circ \textcolor{black}{\mathbb{J}}_{\gamma \mathbb A}$ is firmly nonexpansive.
Let $U = (I - \textcolor{black}{\mathbb{J}}_{\gamma \mathbb B})$, $T_1 = \textcolor{black}{\mathbb{J}}_{\gamma \mathbb C} \circ \textcolor{black}{\mathbb{J}}_{\gamma \mathbb A}$, $V = 2\textcolor{black}{\mathbb{J}}_{\gamma \mathbb B} - I$, then we have $S = T$ and $W=0$ in Lemma \ref{lemma:DY}, thus \eqref{eq:lemmaDY} gives
\begin{equation}
\lVert Tz-Tw \rVert^2 \leq \lVert z-w \rVert^2 - \frac{1-\frac{1}{2}}{\frac{1}{2}}\lVert (I-T)z - (I-T)w \rVert^2, \label{eq:averagedness}
\end{equation}
which implies $T$ is firmly nonexpansive thus $\frac{1}{2}$-averaged. Therefore,   results in
Section \ref{section_weak_convergence_rates}, apply to establish the convergence of the splitting scheme \eqref{FPI_algorithm_convex_case}.  
\end{proof}
\begin{rem}
\label{rem:advantage}
For a convex minimization $d_1(x)+d_2(x)+d_3(x)$ with differentiable $d_2$, Theorem \ref{thm-specialcase} shows that the proposed splitting scheme \eqref{intro:news} converges for any step size $\gamma>0$ under Assumption \ref{assump}. Even under Assumption \ref{assump}, the Davis-Yin splitting \eqref{intro:dys} still needs the step size to be small enough $\gamma<\frac{2}{L}$, since $\nabla d_2$ cannot be cancelled out like in Lemma \ref{lemma:specialcase} for the proposed splitting scheme \eqref{intro:news}. Thus for special problems satisfying Assumption \ref{assump},   the proposed new scheme is provably much more robust than Davis-Yin splitting in the sense of provable convergence using any positive step size.
\end{rem}
}

\subsection{Some basic properties}
\textcolor{black}{In this subsection, we discuss some basic properties of the proposed splitting scheme \eqref{FPI_algorithm_convex_case} under certain assumptions.}
We introduce the variables $w_{d_i}$ by writing \eqref{FPI_algorithm_convex_case} as follows: \begin{eqnarray}\label{FPI_iterations}
\begin{rcases}
&w^k_{d_3}:= \lambda^{k+1}=\operatorname{prox}_{\gamma d_3}(z^k)\\
& w^k_{d_1}:=p^{k+1}= \operatorname{prox}_{\gamma d_1} (2\lambda^{k+1}-z^k-\gamma\partial d_2(\lambda^{k+1}))\\
&w^k_{d_2}:=v^{k+1}=\operatorname{prox}_{\gamma d_2} (p^{k+1}+\gamma\partial d_2(\lambda^{k+1}))\\
&z^{k+1}=z^k+w^k_{d_2}-w^k_{d_3}.
\end{rcases}
\end{eqnarray}
For convenience we use
$\textcolor{black}{\tilde\nabla} f(x)\in \partial f(x)$
to denote a subgradient of $f$ at $x$.
\begin{prop}
\label{prop3}
Let $z^0\in\mathcal{X}$ and $(z_j)_{j\geq 0}$ be the sequence generated by \eqref{FPI_iterations}. Then the following identities hold: 
\begin{enumerate}[(i)]
\item $w^k_{d_3}=z^k-\gamma \textcolor{black}{\tilde\nabla} d_3(w^k_{d_3})$
\item $ w^k_{d_1}-w^k_{d_3}= -\gamma(\textcolor{black}{\tilde\nabla} d_1(w^k_{d_1})+\nabla d_2(w^k_{d_3})+\textcolor{black}{\tilde\nabla} d_3(w^k_{d_3}))$
 \item $ w^k_{d_2}-w^k_{d_3}= -\gamma(\textcolor{black}{\tilde\nabla} d_1(w^k_{d_1})+\nabla d_2(w^k_{d_2})+\textcolor{black}{\tilde\nabla} d_3(w^k_{d_3}))$
 \item $w^k_{d_1}-w^k_{d_3}= w^k_{d_2}-w^k_{d_3}+\gamma(\nabla d_2(w^{k}_{d_2})-\nabla d_2(w^k_{d_3}))$
\item $z^{k+1}-z^k=w^k_{d_2}-w^{k}_{d_3}= -\gamma(\textcolor{black}{\tilde\nabla} d_1(w^k_{d_1})+\nabla d_2(w^k_{d_2})+\textcolor{black}{\tilde\nabla} d_3(w^k_{d_3})).$
\end{enumerate}
\end{prop}

\begin{proof}
We start with the relation $w^k_{d_3}=\operatorname{prox}_{d_3}^{\gamma} (z^k)$, we get 
\begin{eqnarray}\label{modified_fpi_1}
w^k_{d_3}=z^k-\gamma \textcolor{black}{\tilde\nabla} d_3(w^k_{d_3}).
\end{eqnarray}
Next, we have the following relation 
\begin{eqnarray}\label{modified_fpi_2}
 &&w^k_{d_1}=\operatorname{prox}_{\gamma d_1} (2w^k_{d_3}-z^k-\gamma\nabla d_2(w^k_{d_3}))\nonumber\\
 &\iff& (I+\gamma\textcolor{black}{\tilde\nabla} d_1)(w^k_{d_1})= 2w^k_{d_3}-z^k-\gamma\nabla d_2(w^k_{d_3})\nonumber\\
  &\iff& w^k_{d_1}-w^k_{d_3}= w^k_{d_3}-z^k-\gamma\nabla d_2(w^k_{d_3})-\gamma\textcolor{black}{\tilde\nabla} d_1(w^k_{d_1})\nonumber\\
   &\iff& w^k_{d_1}-w^k_{d_3}= -\gamma(\textcolor{black}{\tilde\nabla} d_1(w^k_{d_1})+\nabla d_2(w^k_{d_3})+\textcolor{black}{\tilde\nabla} d_3(w^k_{d_3})).
\end{eqnarray} 
Again, we have the following relation
\begin{eqnarray}\label{modified_fpi_3}
 &&w^k_{d_2}=\operatorname{prox}_{\gamma d_2} (w^k_{d_1}+\gamma\nabla d_2(w^k_{d_3}))\nonumber\\
 &\iff&(I+\gamma \nabla d_2) (w^k_{d_2})= w^k_{d_1}+\gamma\nabla d_2(w^k_{d_3}) \nonumber\\
  &\iff& w^k_{d_2}= (w^k_{d_3}-\gamma\textcolor{black}{\tilde\nabla} d_1(w^k_{d_1})-\gamma\nabla d_2(w^k_{d_3})-\gamma\textcolor{black}{\tilde\nabla} d_3(w^k_{d_3}))+\gamma \nabla d_2(w^k_{d_3})-\gamma\nabla d_2(w^k_{d_2})\text{ from \eqref{modified_fpi_2}}\nonumber\\
   &\iff& w^k_{d_2}-w^k_{d_3}= -\gamma(\textcolor{black}{\tilde\nabla} d_1(w^k_{d_1})+\nabla d_2(w^k_{d_2})+\textcolor{black}{\tilde\nabla} d_3(w^k_{d_3})).
\end{eqnarray}
Finally, we get
\begin{eqnarray}\label{modified_fpi_4}
 w^k_{d_1}-w^k_{d_3}= w^k_{d_2}-w^k_{d_3}+\gamma(\nabla d_2(w^{k}_{d_2})-\nabla d_2(w^k_{d_3}))\text{ and } z^{k+1}-z^k=w^k_{d_2}-w^{k}_{d_3}.   
\end{eqnarray}
\end{proof}

\begin{prop}\emph{(Upper Inequality).}\label{upper_ineq} Let $w\in\mathcal{X}$ and $w^*$ be the fixed point of the FPI algorithm given in \eqref{FPI_algorithm_convex_case}. Then, the following inequality holds:
\begin{eqnarray*}
&& 2\gamma (d_1(w^k_{d_1})+d_2(w^k_{d_3})+d_3(w^k_{d_3})-(d_1+d_2+d_3)(w^*)) \\
&&\leq \lVert z^k-w^*\rVert^2-\lVert z^k-z^{k+1}\rVert^2-\lVert z^{k+1}-w^*\rVert^2+ 2\gamma\langle\nabla d_2(w^k_{d_3})- \nabla d_2 (w^k_{d_2}),z^k-w^*\rangle\\
&&-2\gamma \langle\nabla d_2(w^k_{d_3})- \nabla d_2 (w^k_{d_2}),z^k-z^{k+1}\rangle+2\gamma \langle z^k-z^{k+1},\nabla d_2(w^k_{d_2})\rangle\\
&&+2 \gamma^2\langle\nabla d_2(w^k_{d_3})-\nabla d_2 (w^k_{d_2}),\nabla d_2(w^k_{d_2})\rangle.   
\end{eqnarray*}
    
\end{prop}
\begin{proof}
We will show that the required inequality holds for every $k\geq 0$. In observance of subgradient inequality, we get
\begin{eqnarray}
&&2\gamma (d_1(w^k_{d_1})+d_2(w^k_{d_3})+d_3(w^k_{d_3})-(d_1+d_2+d_3)(w^*))\nonumber\\
&\leq& 2\gamma \left(\langle w^k_{d_1}-w^*,\textcolor{black}{\tilde\nabla} d_1(w^k_{d_1}) \rangle + \langle w^k_{d_3}-w^*,\nabla d_2(w^k_{d_3})+\textcolor{black}{\tilde\nabla} d_3(w^k_{d_3}) \rangle \right)\nonumber\\
&=& 2\gamma \langle w^k_{d_1}-w^k_{d_3},\textcolor{black}{\tilde\nabla} d_1(w^k_{d_1}) \rangle +2\gamma\langle w^k_{d_3}-w^*,\textcolor{black}{\tilde\nabla} d_1(w^k_{d_1})+\nabla d_2(w^k_{d_3})+\textcolor{black}{\tilde\nabla} d_3(w^k_{d_3}) \rangle\nonumber\\
&=&2 \langle w^k_{d_1}-w^k_{d_3},\gamma\textcolor{black}{\tilde\nabla} d_1(w^k_{d_1}) \rangle +2\langle w^k_{d_3}-w^*, w^k_{d_3}-w^k_{d_1}\rangle  \text{ from (ii) of Proposition \ref{prop3}}\nonumber \\
&=&2 \langle w^k_{d_3}-w^k_{d_1},w^k_{d_3}-w^*-\gamma\textcolor{black}{\tilde\nabla} d_1(w^k_{d_1}) \rangle \nonumber\\
&=& 2 \langle w^k_{d_3}-w^k_{d_2}+\gamma(\nabla d_2(w^{k}_{d_3})-\nabla d_2(w^k_{d_2})),w^k_{d_3}-w^*-\gamma\textcolor{black}{\tilde\nabla} d_1(w^k_{d_1})\rangle \text{ from (iv) of Proposition \ref{prop3}}\nonumber \\
&=& 2 \langle z^k-z^{k+1}+\gamma\nabla d_2(w^k_{d_3})-\gamma \nabla d_2 (w^k_{d_2}),w^k_{d_3}-w^*-\gamma\textcolor{black}{\tilde\nabla} d_1(w^k_{d_1})\rangle \text{ from (v) of Proposition \ref{prop3}} \nonumber\\
&=& 2 \langle z^k-z^{k+1}+\gamma\nabla d_2(w^k_{d_3})-\gamma \nabla d_2 (w^k_{d_2}),z^k+w^k_{d_3}-z^k-w^*-\gamma\textcolor{black}{\tilde\nabla} d_1(w^k_{d_1})\rangle \nonumber\\
&=& 2 \langle z^k-z^{k+1}+\gamma\nabla d_2(w^k_{d_3})-\gamma \nabla d_2 (w^k_{d_2}),z^k-\gamma\textcolor{black}{\tilde\nabla} d_3(w^k_{d_3})-w^*-\gamma\textcolor{black}{\tilde\nabla} d_1(w^k_{d_1})\rangle \nonumber\\
&=& 2 \langle z^k-z^{k+1}+\gamma\nabla d_2(w^k_{d_3})-\gamma \nabla d_2 (w^k_{d_2}),z^k-\gamma\textcolor{black}{\tilde\nabla} d_3(w^k_{d_3})-w^*-\gamma\textcolor{black}{\tilde\nabla} d_1(w^k_{d_1})-\gamma\nabla d_2(w^k_{d_2})+\gamma\nabla d_2(w^k_{d_2})\rangle \nonumber\\
&=& 2 \langle z^k-z^{k+1}+\gamma\nabla d_2(w^k_{d_3})-\gamma \nabla d_2 (w^k_{d_2}),z^k-(z^k-z^{k+1})-w^*\rangle \nonumber\\
&&+2 \langle z^k-z^{k+1}+\gamma\nabla d_2(w^k_{d_3})-\gamma \nabla d_2 (w^k_{d_2}),\gamma\nabla d_2(w^k_{d_2})\rangle \text{ from (v) of Proposition \ref{prop3}} \nonumber\\
&=& 2 \langle z^k-z^{k+1},z^k-w^*\rangle-2\langle z^k-z^{k+1},z^k-z^{k+1}\rangle+ 2\langle\gamma\nabla d_2(w^k_{d_3})-\gamma \nabla d_2 (w^k_{d_2}),z^k-w^*\rangle \nonumber\\
&&-2 \langle\gamma\nabla d_2(w^k_{d_3})-\gamma \nabla d_2 (w^k_{d_2}),z^k-z^{k+1}\rangle+2 \langle z^k-z^{k+1},\gamma\nabla d_2(w^k_{d_2})\rangle \nonumber\\
&&+2 \langle \gamma\nabla d_2(w^k_{d_3})-\gamma \nabla d_2 (w^k_{d_2}),\gamma\nabla d_2(w^k_{d_2})\rangle \label{upper_ineq_eqn_a}\\
\end{eqnarray}
\begin{eqnarray}    
&=& \lVert z^k-z^{k+1}\rVert^2+ \lVert z^k-w^*\rVert^2-\lVert z^{k+1}-w^*\rVert^2-2\lVert z^k-z^{k+1}\rVert^2+ 2\langle\gamma\nabla d_2(w^k_{d_3})-\gamma \nabla d_2 (w^k_{d_2}),z^k-w^*\rangle \nonumber\\
&&-2 \langle\gamma\nabla d_2(w^k_{d_3})-\gamma \nabla d_2 (w^k_{d_2}),z^k-z^{k+1}\rangle+2 \langle z^k-z^{k+1},\gamma\nabla d_2(w^k_{d_2})\rangle \nonumber\\
&&+2 \langle \gamma\nabla d_2(w^k_{d_3})-\gamma \nabla d_2 (w^k_{d_2}),\gamma\nabla d_2(w^k_{d_2})\rangle \text{ from \eqref{cosine_rule}}\nonumber \\
&=& \lVert z^k-w^*\rVert^2-\lVert z^k-z^{k+1}\rVert^2-\lVert z^{k+1}-w^*\rVert^2+ 2\gamma\langle\nabla d_2(w^k_{d_3})- \nabla d_2 (w^k_{d_2}),z^k-w^*\rangle \nonumber\\
&&-2\gamma \langle\nabla d_2(w^k_{d_3})- \nabla d_2 (w^k_{d_2}),z^k-z^{k+1}\rangle+2\gamma \langle z^k-z^{k+1},\nabla d_2(w^k_{d_2})\rangle \nonumber\\
&&+2 \gamma^2\langle\nabla d_2(w^k_{d_3})-\nabla d_2 (w^k_{d_2}),\nabla d_2(w^k_{d_2})\rangle. \nonumber 
\end{eqnarray}
\end{proof}

\begin{prop}\emph{(Lower Inequality)}.\label{lower_ineq} Let $w\in\mathcal{X}$ and $w^*$ be the fixed point of the FPI algorithm given in \eqref{FPI_algorithm_convex_case}. Then, the following inequality holds:  
\begin{eqnarray*}
 &&2\gamma(d_1(w^k_{d_1})+ d_2(w^k_{d_3})+ d_3(w^k_{d_3})- (d_1+d_2+d_3)(w^*))\\
 &&\geq \langle w^k_{d_2}-w^k_{d_3},\textcolor{black}{\tilde\nabla} d_1(w^*)\rangle +\langle\gamma \nabla d_2(w^k_{d_2})-\gamma\nabla d_2(w^k_{d_3}),\textcolor{black}{\tilde\nabla} d_1 (w^*)\rangle.    
\end{eqnarray*}
\end{prop}

\begin{proof}
  By subgradient inequality and Proposition \ref{prop3} (iii), we have
\begin{align*}
& d_1(w^k_{d_1})-d_1(w^*)\\
\geq& \langle w_{d_1}^k-w^*,\nabla d_1 (w^*)\rangle\\
=&\langle w_{d_1}^k-w_{d_3}^k,\nabla d_1 (w^*)\rangle+\langle w_{d_3}^k-w^*,\textcolor{black}{\tilde\nabla} d_1 (w^*)\rangle\\
=&\langle w_{d_2}^k-w_{d_3}^k+\gamma \nabla d_2(w^k_{d_2})-\gamma\nabla d_2(w^k_{d_3}),\textcolor{black}{\tilde\nabla} d_1 (w^*)\rangle+\langle w_{d_3}^k-w^*,\textcolor{black}{\tilde\nabla} d_1 (w^*)\rangle.
\end{align*}
In a similar manner, we obtain
\begin{eqnarray*}
&& d_2(w^k_{d_3})-d_2(w^*)\geq \langle w_{d_3}^k-w^*,\nabla d_2 (w^*)\rangle\\
\text{and }&& d_3(w^k_{d_3})-d_3(w^*)\geq \langle w_{d_3}^k-w^*,\textcolor{black}{\tilde\nabla} d_3 (w^*)\rangle.
\end{eqnarray*}
On adding the above three relations, we get
\begin{eqnarray}\label{lower_ineq_eqn1}
    &&d_1(w^k_{d_1})+ d_2(w^k_{d_3})+ d_3(w^k_{d_3})- (d_1+d_2+d_3)(w^*) \nonumber\\
    &&\geq \langle w^k_{d_2}-w^k_{d_3},\textcolor{black}{\tilde\nabla} d_1(w^*)\rangle +\langle w_{d_3}^k-w^*,\textcolor{black}{\tilde\nabla} d_1 (w^*)+\nabla d_2 (w^*)+\textcolor{black}{\tilde\nabla} d_3 (w^*)\rangle \nonumber\\
    &&+\langle\gamma \nabla d_2(w^k_{d_2})-\gamma\nabla d_2(w^k_{d_3}),\textcolor{black}{\tilde\nabla} d_1 (w^*)\rangle.
    \end{eqnarray}
Now, from the optimality condition of a subdifferential set, we can assume that
\[
\textcolor{black}{\tilde\nabla} d_1 (w^*)+\nabla d_2 (w^*)=-\textcolor{black}{\tilde\nabla} d_3 (w^*)\in \partial d_3(w^*).
\]
Therefore, on plugging the above relation in \eqref{lower_ineq_eqn1}, we get
\begin{eqnarray*}    
 &&d_1(w^k_{d_1})+ d_2(w^k_{d_3})+ d_3(w^k_{d_3})- (d_1+d_2+d_3)(w^*)\\
 &&\geq \langle w^k_{d_2}-w^k_{d_3},\textcolor{black}{\tilde\nabla} d_1(w^*)\rangle +\langle\gamma \nabla d_2(w^k_{d_2})-\gamma\nabla d_2(w^k_{d_3}),\textcolor{black}{\tilde\nabla} d_1 (w^*)\rangle.
 \end{eqnarray*}
\end{proof}

In the next theorem, we prove the convergence rate of the proposed Algorithm \eqref{FPI_algorithm_convex_case}.
\begin{thm}\label{convergence_rates}
\label{theorem-convergence-underassumption}
Consider the iteration \eqref{FPI_algorithm_convex_case}. Assume that the function $d_1$ is $L$-Lipschitz continuous on the closed ball $\overline{B(0,(1+\gamma/\beta)\lVert w^0-w^*\rVert})$. Under the same assumptions in Theorem \ref{convergence_theorem1}, we have  the following convergence 
\[
(d_1+d_2+d_3)(w^k_{d_3})-(d_1+d_2+d_3)(w^*)=o\left(\tfrac{1}{\sqrt{k+1}}\right).
\]
\end{thm}

\begin{proof}
Note that in view of Lemma \ref{closed_ball}, the sequences $w^k_{d_1}$ and $w^k_{d_3}$ are within the region where $d_1$ is $L$-Lipschitz continuous. Therefore, we have
\begin{eqnarray*}
&&2\gamma\left(d_1(w^k_{d_1})+d_2(w^k_{d_3})+d_3(w^k_{d_3})-(d_1+d_2+d_3)(w^*)\right)\\
&&\leq2\gamma\left(d_1(w^k_{d_3})+d_2(w^k_{d_3})+d_3(w^k_{d_3})-(d_1+d_2+d_3)(w^*)\right)+2\gamma L\lVert w^k_{d_1}-w^k_{d_3}\rVert.
\end{eqnarray*}
Now, employing relation \eqref{upper_ineq_eqn_a} of Proposition \ref{upper_ineq} in above relation, we get
\begin{eqnarray}\label{rate_1}
&&2\gamma\left(d_1(w^k_{d_1})+d_2(w^k_{d_3})+d_3(w^k_{d_3})-(d_1+d_2+d_3)(w^*)\right)\nonumber\\
&&\leq 2 \langle z^k-z^{k+1},z^k-w^*\rangle-2\langle z^k-z^{k+1},z^k-z^{k+1}\rangle+ 2\langle\gamma\nabla d_2(w^k_{d_3})-\gamma \nabla d_2 (w^k_{d_2}),z^k-w^*\rangle\nonumber \\
&&-2 \langle\gamma\nabla d_2(w^k_{d_3})-\gamma \nabla d_2 (w^k_{d_2}),z^k-z^{k+1}\rangle+2 \langle z^k-z^{k+1},\gamma\nabla d_2(w^k_{d_2})\rangle \nonumber\\
&&+2 \langle \gamma\nabla d_2(w^k_{d_3})-\gamma \nabla d_2 (w^k_{d_2}),\gamma\nabla d_2(w^k_{d_2})\rangle+2\gamma L\lVert w^k_{d_1}-w^k_{d_3}\rVert \text{ from \eqref{cosine_rule}} \nonumber\\
&&\leq 2 \|z^k-z^{k+1}\| \|z^k-w^*\|+2\|z^k-z^{k+1}\| \|z^k-z^{k+1}\|+ 2\gamma\|\nabla d_2(w^k_{d_3})-  \nabla d_2 (w^k_{d_2})\| \|z^k-w^*\|\nonumber \\
&&+2  \gamma \|\nabla d_2(w^k_{d_3})-  \nabla d_2 (w^k_{d_2})\| \|z^k-z^{k+1}\|+2 \gamma\|z^k-z^{k+1}\| \|\nabla d_2(w^k_{d_2})\| \nonumber\\
&&+2  \gamma^2\| \nabla d_2(w^k_{d_3})- \nabla d_2 (w^k_{d_2})\|  \|\nabla d_2(w^k_{d_2})\|+2\gamma L \lVert w^k_{d_1}-w^k_{d_3}\rVert  \nonumber
\end{eqnarray}
Now, note that in view of Theorem \ref{convergence_theorem1}, $z^k$ is bounded. Therefore, in view of Lipschitz continuity of $\operatorname{prox}_{\gamma d_3}$ and $\nabla d_2$, we have that $w^k_{d_3}=\operatorname{prox}_{\gamma d_3}(z^k)$ and $\nabla d_2(w^k_{d_3})$ are bounded. Similarly, $\nabla d_2(w^k_{d_2})$ is also bounded.\\
By Theorem \ref{convergence_theorem1},  $\lVert z^k-z^{k+1}\rVert=o\left(\tfrac{1}{\sqrt{k+1}}\right)$. By Proposition \ref{prop3} (v), $\lVert w^k_{d_2}-w^k_{d_3}\rVert=o\left(\tfrac{1}{\sqrt{k+1}}\right)$. The Lipschitz continuity of $\nabla d_2$, the gradients implies $\|\nabla d_2(w^k_{d_2})-\nabla d_2(w^k_{d_3})\|=o\left(\tfrac{1}{\sqrt{k+1}}\right)$. 
By Proposition \ref{prop3} (iv), we also have $\lVert w^k_{d_1}-w^k_{d_3}\rVert=o\left(\tfrac{1}{\sqrt{k+1}}\right)$. Thus the order $o\left(\tfrac{1}{\sqrt{k+1}}\right)$ is proven.
\end{proof}

\subsection{Multi-operator splitting schemes}
Algorithm \ref{FPI_algorithm_modified_convex_case} can be formally extended to a four-block problem. Consider the following problem given by 
\begin{eqnarray*}
 \begin{rcases}
\min& f_1(x_1)+f_2(x_2)+f_3(x_3)+f_4(x_4),\\
\text{subject to }&A_1x_1+A_2x_2+A_3x_3+A_4x_4=b,\\
&x_1\in\mathcal{X}_1,~x_2\in\mathcal{X}_2,x_3\in\mathcal{X}_3,x_4\in\mathcal{X}_4
\end{rcases},   
\end{eqnarray*}
where for each $i=1,2,3,4$, the functions $f_i:\mathbb{R}^{n_i}\to\mathbb{R}\cup\{+\infty\}$ are proper, closed, and convex (not necessarily smooth). 
For each $i=1,2,3,4$, $A_i\in\mathbb{R}^{m\times n_i}$  and $b\in\mathbb{R}^4$.
\textcolor{black}{
Notice that the proposed three-operator splitting scheme \eqref{FPI_algorithm_convex_case} can be regarded as a modification to the dual form of three-block ADMM \eqref{dual-ADMM}.
Thus similar to derivations in Section \ref{sec:admm-derivation}, we can first derive the dual form of the classical four-block ADMM for the problem above, then perform similar modifications to obtain a formal extension of the proposed scheme \eqref{FPI_algorithm_convex_case}  for 
 minimizing the dual problem $\sum_{i=1}^4 d_i(x)$ where each $d_i$ is a proper closed convex function.}\\
 The  extended $4$-operator splitting algorithm is given as follows:
\begin{algorithm}[H]
\caption{Initialize $z^0\in\mathcal{X},~\gamma\in(0,\min\{L_2, L_3\})$. $L_i$ is the Lipschitz constant of $\nabla d_i$.}
\begin{enumerate}
\item Compute $ w^{k+1}=\operatorname{prox}_{\gamma d_4}(z^k)$;
\item Compute $p^{k+1}= \operatorname{prox}_{\gamma d_1} (2w^{k+1}-z^k-\gamma \nabla d_2(w^{k+1})-\gamma \nabla d_3(w^{k+1}))$;
\item Compute $v_1^{k+1}=\operatorname{prox}_{\gamma d_2}  (p^{k+1}+\gamma\nabla d_2(w^{k+1}))$;
\item Compute $v_2^{k+1}=\operatorname{prox}_{\gamma d_3}  (v_1^{k+1}+\gamma\nabla d_3(w^{k+1}))$;
\item Update $z^{k+1}=z^k+v_{2}^{k+1}-w^{k+1}.$
\end{enumerate}
\end{algorithm}
\begin{rem}
    In general, a splitting scheme for $m$-operator can be similarly derived from \textcolor{black}{modifying the dual form of} $m$-block ADMM for minimizing $\sum_{i=1}^m d_i(x)$ where each $d_i$ is a proper closed convex function, and at least $m-2$ of them have Lipschitz continuous gradients.
\end{rem}

\section{Applications and Numerical Results}\label{section_numerical}
In this section, we demonstrate the working of our proposed Algorithm \ref{operator_T_algorithm1} and \ref{FPI_algorithm_modified_convex_case}. Next, we compare the behavior of the proposed Algorithm \ref{operator_T_algorithm1} with the method in \cite{davis2017three} \textcolor{black}{and \cite{ryu2020finding}}.  

\begin{example}\label{example_1}
Consider the following problem given by
\begin{eqnarray}\label{convex_problem1}
\underset{x\in\mathbb{R}^n}{\min}~f(x)+g(Lx)+h(x),    
\end{eqnarray}
where $f,~g,$ and $h$ are proper, closed, and convex functions, and $g$ is $(1/\beta)$-Lipschitz differentiable, and $L$ is linear mapping. The proposed Algorithm \ref{operator_T_algorithm1} applies here with the following monotone operators given by:
\[
\mathbb{A}=\partial f; ~\mathbb{B}=\nabla(g\circ L)= L^*\circ \nabla g\circ L; ~ \mathbb{C}=\partial h.
\]
If zer $(\partial f+ \nabla g+\partial h)\neq \emptyset$. Then, $x^k$ is a weakly minimal solution to \eqref{convex_problem1}. The modified form of Algorithm \ref{operator_T_algorithm1} for problem \eqref{convex_problem1} is discussed below:

\begin{algorithm}[!htb]
\caption{}\label{example_algo1}
Initialize an arbitrary $z^0\in\mathcal{X},$ stepsize $\gamma\in(0,2\beta/\lVert L\rVert^2)$ and $(\lambda_j)_{j\geq 0}\in (0, (4\beta-\gamma \lVert L\rVert^2)/2\beta)$. For $k=0,1,\ldots$
\begin{enumerate}
\item Compute $ x^{k+1}=\operatorname{prox}_{\gamma h}(z^k)$;
\item Compute  $y^{k+1}=L x^{k+1}$;
\item Compute $p^{k+1}= \operatorname{prox}_{\gamma f}  (2x^{k+1}-z^k-\gamma L^*\nabla g(y^{k+1}))$;
\item Compute $v^{k+1}=\operatorname{prox}_{\gamma g} (p^{k+1}+\gamma\nabla g(y^{k+1}))$;
\item Update $z^{k+1}=z^k+\lambda_k (v^{k+1}-w^{k+1}).$
\end{enumerate}
\end{algorithm}
\end{example}

\begin{example}\label{example1}
Consider the following problem given by
\begin{eqnarray}
\begin{rcases}
\underset{x\in \mathbb R^n}{\min}~&\tfrac{\alpha}{2}\lVert x- u\rVert^2_2+ {i}_{\Lambda_1}(x)+i_{\Lambda_2}(x),\hspace{0.5cm}\\
\text{subject to}& \Lambda_1=\{x:m\leq x_i\leq M,\quad \forall i\}\\
\text{and}&\Lambda_2 =\{x:Ax=b\},
\end{rcases}
\label{example-limiter}
\end{eqnarray}
where $A=[1,1,\ldots, 1]$, $b\in\mathbb{R}$, $u\in \mathbb R^n$ are given. Here $\alpha>0$ is a fixed constant. 
Such a simple problem can be used as a postprocessing step to enforce bounds for solving complicated PDEs \cite{liu2023simple, liu2024optimization, liu2024opt-2}. 
 
 Notice that this simple constrained minimization \eqref{example-limiter} can also be solved directly via the KKT system of the Lagrangian, which however might be less efficient than splitting methods for large problems, see a comparison of DRS with a direct solver of KKT system in \cite[Appendix]{liu2024optimization}. Moreover,  following the analysis in \cite{demanet2016eventual}, a sharp local linear convergence rate of Douglas-Rachford splitting for \eqref{example-limiter} can be derived, which can be further used to design optimal step size \cite{liu2023simple}.
Though \eqref{example-limiter}  can be solved by two-operator splitting, a more general version of \eqref{example-limiter} can no longer be easily solved by Douglas-Rachford splitting.
For example, for stabilizing numerical schemes solving gas dynamics equations \cite{zhang2012minimum,zhang2017positivity}, the bound-preserving constraint $\Lambda_1=\{x:m\leq x_i\leq M,\quad \forall i\}$ in \eqref{example-limiter} would be
replaced by the invariant domain preserving constraint
$\Lambda_1=\{x\in\mathbb R^{n\times 3}:  x_i\in G\subset \mathbb R^3 \quad \forall i\}$ for some convex invariant domain set $G$,  which a two-operator splitting cannot easily handle.  Instead, a three operator splitting like three-operator Davis-Yin splitting or the proposed splitting can be used.

\textcolor{black}{We conduct a comprehensive comparison of the major splitting schemes surveyed in \eqref{subsection_comparison} that are applicable to this problem. This includes Davis-Yin splitting, Forward-Douglas-Rachford-Forward (FDRF) \eqref{FDRF}, Forward-Reflected Douglas-Rachford (FRDR) \eqref{FRDR}, ADMM Dual Form \eqref{intro:admm}, and our proposed three-operator splitting \eqref{intro:news}.}
Let $f(x)=\tfrac{\alpha}{2}\lVert x- u\rVert^2_2, g(x)= {i}_{\Lambda_1}(x),$ and $h(x)=i_{\Lambda_2}(x)$.
The proposed scheme  \eqref{intro:news} applies here with the following operators given by:
\begin{eqnarray*}
    &&\mathbb{A}=\partial g=
    \begin{cases}
       [0,+\infty], & \text{if }x_i=M\\
       0, &\text{if } x_i\in (m,M) \\
       [-\infty,0], & \text{if }x_i=m,
    \end{cases},\\ &\text{and}&\mathbb{B}=\partial h =\mathcal{R}(A^\top), \\
    &\text{and}&\mathbb{C}=\partial f =\alpha(x-u).
\end{eqnarray*}
Equivalently, we can set $d_1=g, ~d_2=f,~d_3=h$ in Algorithm \ref{FPI_algorithm_modified_convex_case}  or  the scheme \eqref{intro:news} to obtain  
\begin{eqnarray*}
&&x^{k+\frac12}=A^+(b-Az^k)+z^k\\
&&p^{k+1}=\min(\max(2x^{k+\frac12}-z^k+\alpha\gamma(x^{k+\frac12}-u), m),M)\\
&&x^{k+1}=\tfrac{1}{\alpha\gamma+1}p^{k+1}+\tfrac{\alpha\gamma}{\alpha\gamma+1}u\\
&&z^{k+1}=z^k+x^{k+1}-x^{k+\frac12}.
\end{eqnarray*}
The Davis-Yin scheme \eqref{intro:dys} with $d_1=g, ~d_2=f,~d_3=h$ becomes
\begin{eqnarray*}
&&x^{k+\frac12}=A^+(b-Az^k)+z^k\\
&&x^{k+1}=\min(\max(2x^{k+\frac12}-z^k+\alpha\gamma(x^{k+\frac12}-u),m),M)\\
&&z^{k+1}=z^k+x^{k+1}-x^{k+\frac12}.
\end{eqnarray*}
\textcolor{black}{The FDRF scheme with $d_1=g, ~d_2=f,~d_3=h$ becomes
\begin{eqnarray*}
&&x^{k+\frac12}=A^+(b-Az^k)+z^k\\
&&x^{k+1}=\min(\max(2x^{k+\frac12}-z^k+\alpha\gamma(x^{k+\frac12}-u),m),M)\\
&&z^{k+1}=z^k+x^{k+1}-x^{k+\frac12}-\gamma(\alpha(x^{k+1}-u)-\alpha(x^{k+\frac12}-u)).
\end{eqnarray*}
}
\textcolor{black}{The FRDR scheme with $d_1=g, ~d_2=f,~d_3=h$ becomes
\begin{eqnarray*}
&&x^{k+1}=\min(\max(x^{k}-\gamma_{\text{FRDR}} z^k - \gamma_{\text{FRDR}}(2\alpha(x^{k}-u)-\alpha(x^{k-1}-u),m),M)\\
&&y^{k+1}=A^+(b-A(2x^{k+1}-x^k+\beta z^k)+2x^{k+1}-x^k+\beta z^k\\
&&z^{k+1}=z^k+\frac{1}{\beta}(2x^{k+1}-x^{k}-y^{k+1}).
\end{eqnarray*}
} 
We compare the proposed splitting scheme \eqref{intro:news}, Davis-Yin method \eqref{intro:dys} \textcolor{black}{FDRF method, and FRDR method, ADMM dual form} on the problem \eqref{example-limiter} with $\alpha=1, n=100, m=-1, M=1$, and $b=Au, u\in \mathbb R^n$ where $u$ is constructed by perturbing a sine profile by random noise: $$u_i=\sin(2\pi\frac{i}{n})+0.8*\mathcal N(0,1).$$ 
The number of entries in $u$ greater than $M=1$ is $17$ and 
the number of entries in $u$ less than $m=-1$ is $20$. The minimizer $x^*$ to \eqref{example-limiter} satisfies $x^*_i\in [m, M]$. The error measured by $\|x^{k+\frac12}-x^*\|$ is shown in Figure \ref{figure_1}. 
Let $L$ be the \textcolor{black}{true} Lipschitz constant of $\nabla f$. 
\textcolor{black}{Let $\mu$ be an   estimate of the Lipschitz constant $L$, and we use the following step sizes:
\begin{enumerate}
    \item Davis-Yin splitting, FDRF, ADMM dual form and the new splitting method all use the same step size $\gamma=\frac{1}{\mu}$.
  \item   The FRDR method has two parameters and it converges for any $\beta>0$ and $\gamma_{\text{FRDR}} <  \frac{\beta}{1 + 2L\beta}$.
  We use a small enough $\beta>0$ and $\gamma_{\text{FRDR}} =  \frac{\beta}{1 + 2\mu\beta}$.
\end{enumerate} 
For this simple test, we have  $L=\alpha$, but in practice   usually one has to estimate the Lipschitz constant. To this end, we show the performance using different estimates $\mu$
in Figure \ref{figure_1}, in which the optimal solution $x^*$ is generated by running Davis-Yin method enough number of iterations with $\mu = L$.
 The Davis-Yin method performs the best if using step size $\gamma=\frac{1}{L}$, i.e., when there is an accurate estimate of $L$, and the proposed method is not faster than the Davis-Yin method if $\gamma<\frac{2}{L}$. The Davis-Yin method,  FDRF and ADMM dual form will not converge using an underestimated Lipschitz constant like $\mu=\frac{L}{3}$ (or equivalently using a large step size like $\gamma=\frac{3}{L}$).
 Even with a significantly underestimated Lipschitz constant like $\mu=\frac{L}{40}$ (or equivalently using a large step size like $\gamma=\frac{40}{L}$ in the proposed new method),
the proposed new method  and FRDR converge. If using a significantly underestimated Lipschitz constant, FRDR is much faster with a properly tuned parameter $\beta$, but FRDR does not converge without tuning the extra parameter $\beta$.   
Though the proposed new method is slower than FRDR, it does not need an extra parameter thus it needs no further parameter tuning, which is its advantage. In other words, when there is no accurate estimate of the Lipschitz constant, this simple test suggests that the proposed new splitting is robust and easy to use. } 
\begin{figure}[htbp!]
    \centering
     \subfloat[\textcolor{black}{Lipschitz constant estimate $\mu=\frac{L}{0.3}$. The step size is $\gamma=\frac{0.3}{L}$.} ]{\includegraphics[width=0.48\textwidth]{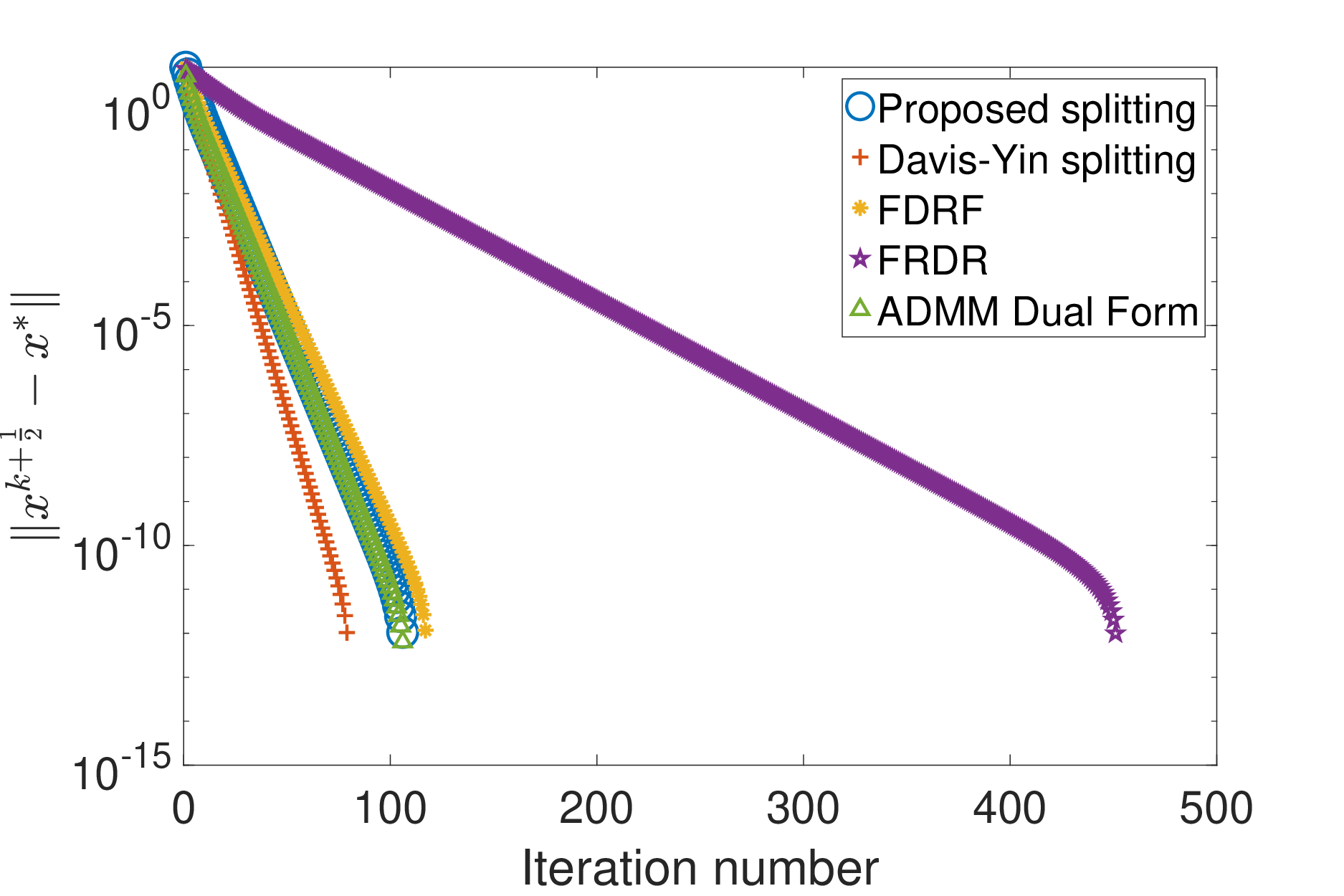}\label{fig1z}} \quad
    \subfloat[\textcolor{black}{Lipschitz constant estimate $\mu=0.99L$. The step size is $\gamma=\frac{0.99}{L}$.}]{\includegraphics[width=0.48\textwidth]{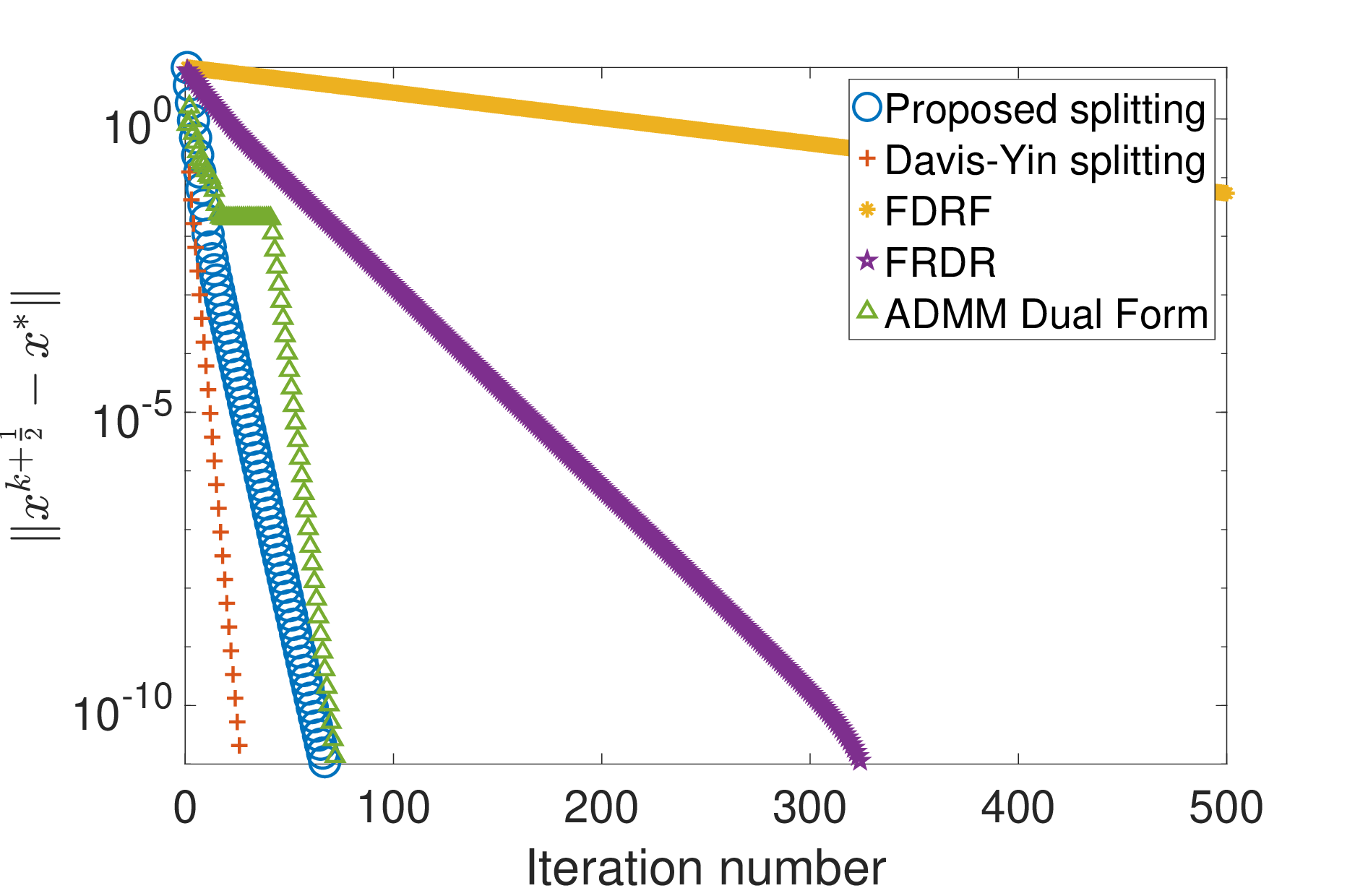}\label{fig1a}} \qquad
    \subfloat[\textcolor{black}{Lipschitz constant estimate $\mu=\frac{L}{1.8}$. The step size is $\gamma=\frac{1.8}{L}$.  ADMM Dual Form and FDRF fail to converge.}]{\includegraphics[width=0.48\textwidth]{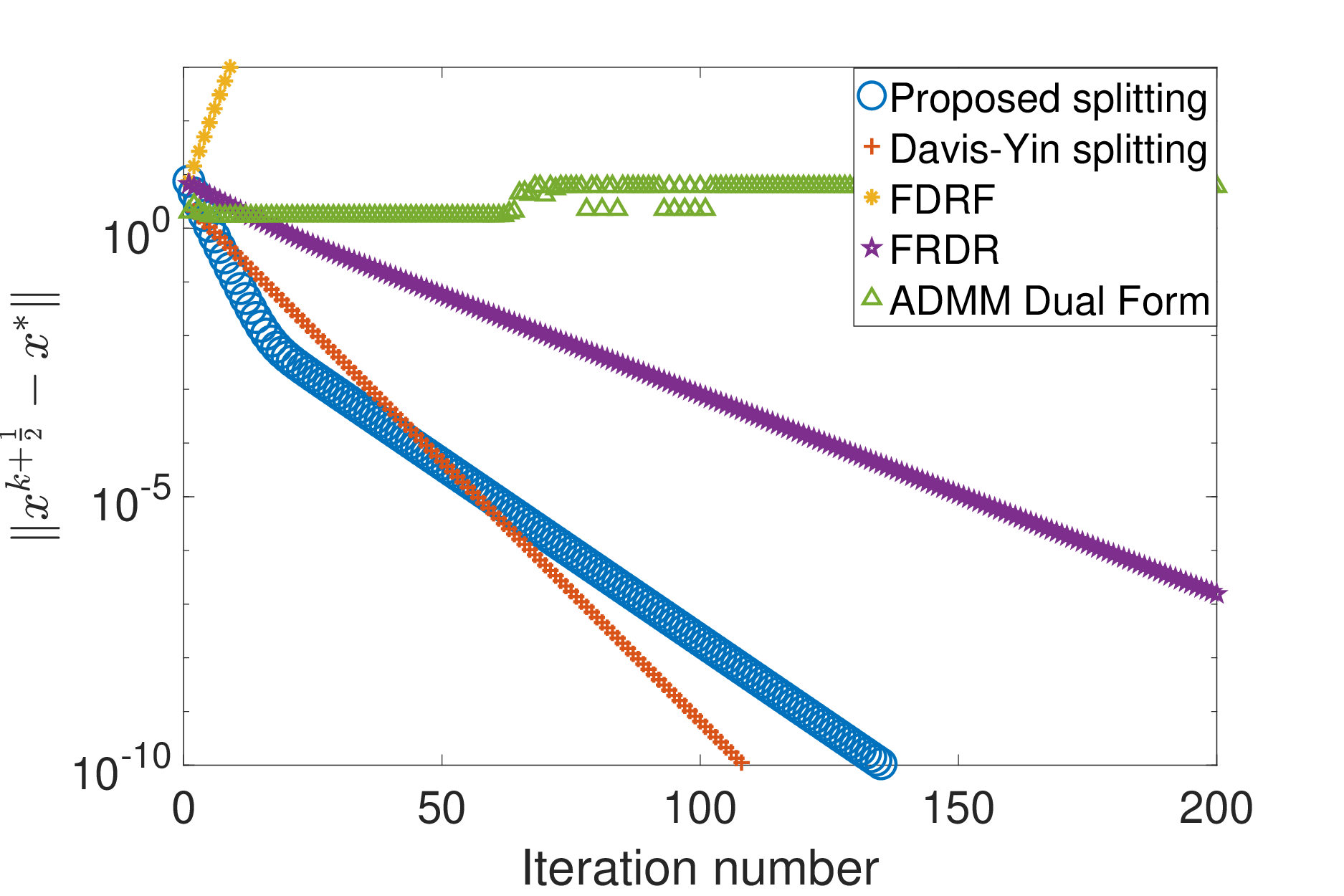} \label{fig1b}} \quad
    \subfloat[\textcolor{black}{Lipschitz constant estimate $\mu=\frac{L}{3}$. The step size is $\gamma=\frac{3}{L}$. ADMM Dual Form, Davis-Yin and FDRF fail to converge.}]{\includegraphics[width=0.48\textwidth]{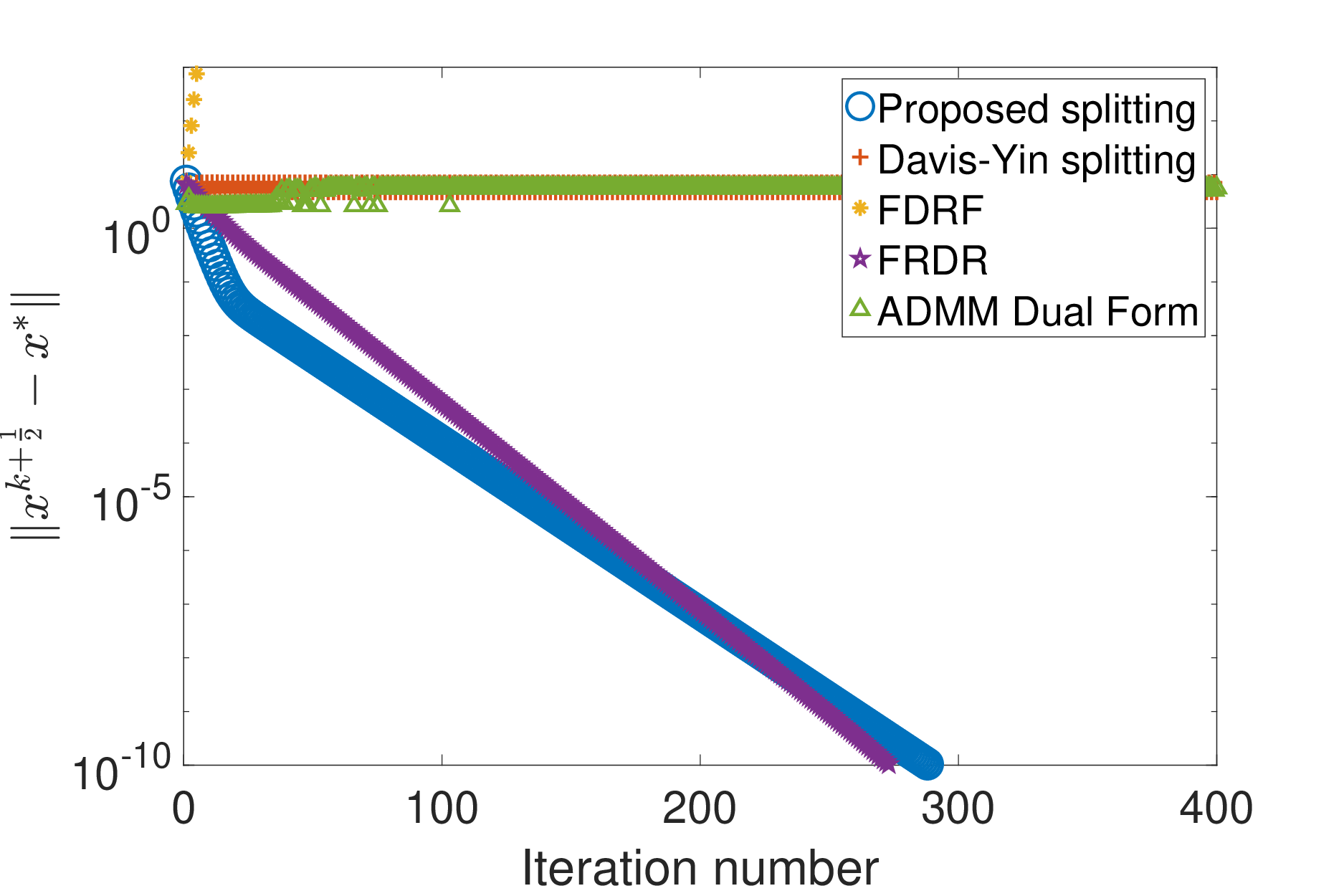}\label{fig1c}}\qquad
    \subfloat[\textcolor{black}{Lipschitz constant estimate $\mu=\frac{L}{20}$. The step size is $\gamma=\frac{20}{L}$. ADMM Dual Form, Davis-Yin and FDRF fail to converge. }]{\includegraphics[width=0.48\textwidth]{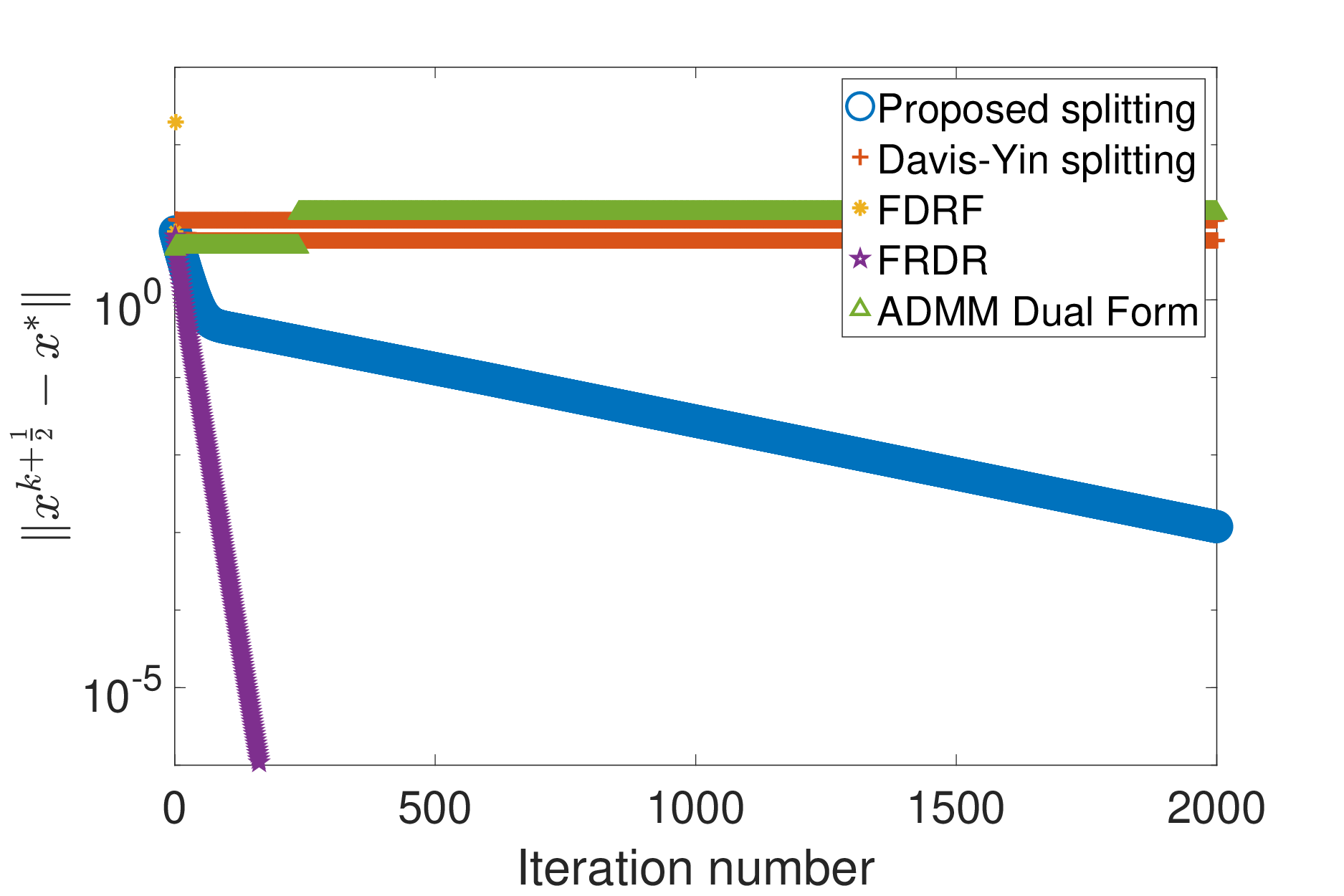}\label{fig1e}}\quad
    \subfloat[\textcolor{black}{Lipschitz constant estimate $\mu=\frac{L}{40}$. The step size is $\gamma=\frac{40}{L}$. ADMM Dual Form, Davis-Yin and FDRF fail to converge.}]{\includegraphics[width=0.48\textwidth]{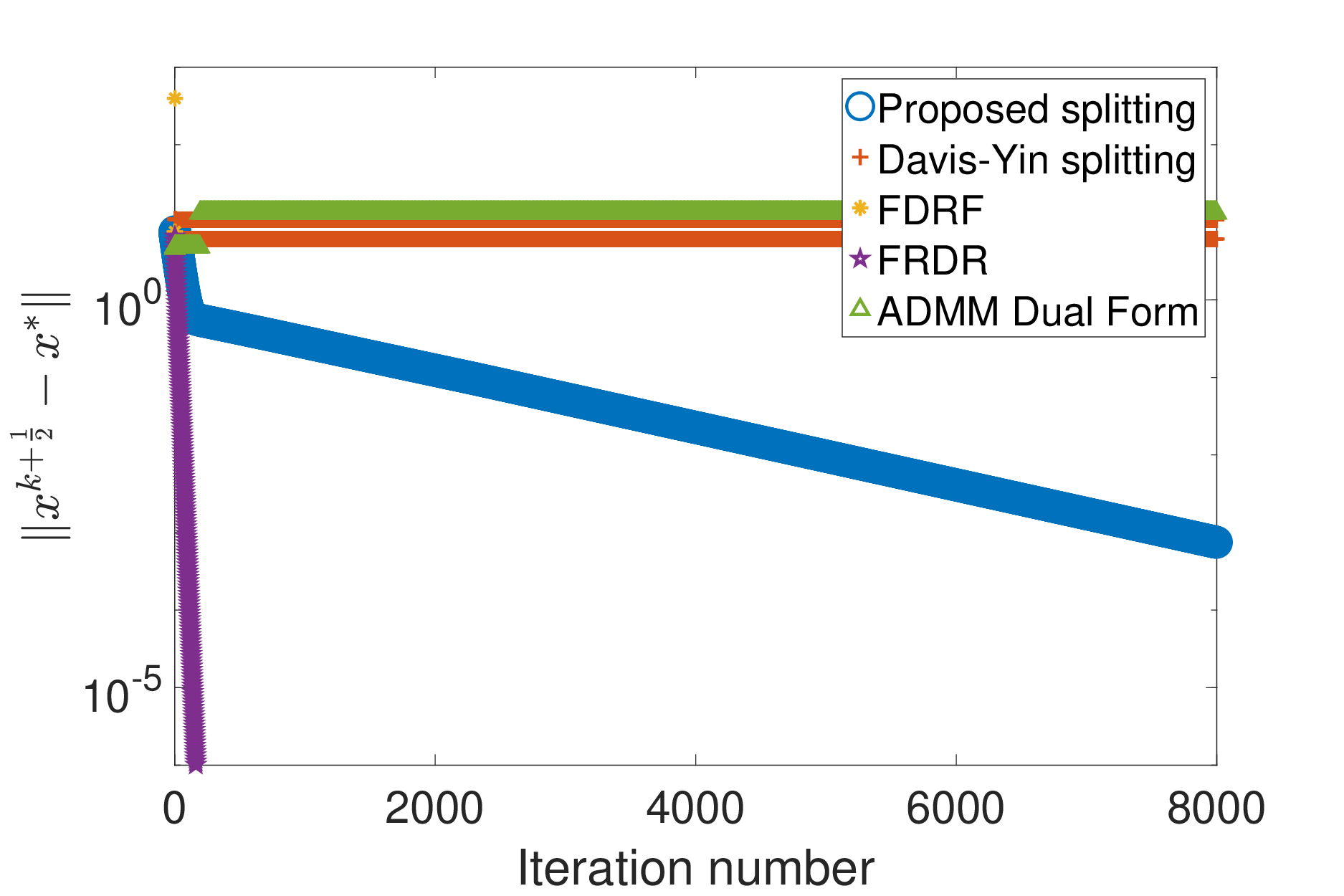}\label{fig1f}}
    \caption{ \textcolor{black}{ Results in Example \ref{example1}   for performance of using different values of $\mu$, which is an estimate of the Lipschitz constant $L=\alpha$.  For ADMM Dual Form, Davis-Yin splitting, FDRF and the proposed new splitting method, we use the same step size $\gamma=\frac{1}{\mu}$.
   For the FRDR method,
we use a small enough $\beta=0.1$ and $\gamma_{\text{FRDR}} =  \frac{\beta}{1 + 2\mu\beta}$. FRDR method will diverge without using small enough $\beta$ for $\mu\ll L$. }
    }
    \label{figure_1}
\end{figure}
\end{example}
\textcolor{black}{
\begin{example}\label{example2}
Next we test the performance of the proposed new splitting on a fused lasso in the following form
\begin{equation}
\underset{x\in\mathbb R^n}{\min}\frac{1}{2}\lVert Ax-b \rVert^2_2 + \mu_1\lVert x\rVert_1+\mu_2\lVert Bx\rVert_1\hspace{0.5cm}
\label{example-flasso}
\end{equation}
where $A\in\mathbb{R}^{r\times n},b\in\mathbb{R}^r$, and
$B = \begin{bmatrix}
-1 &  1 & & & \\
 & -1 &  1 \\
& & \dots & \dots \\
& & &-1 &  1
\end{bmatrix} \in\mathbb{R}^{(n-1)\times n}.$
The fused lasso has been used in signal processing, genomics, and machine learning to produce sparse and locally smoothness models.
When using a splitting method like Davis-Yin and the proposed new method in this paper for \eqref{example-flasso}, one needs the proximal operator for the TV-norm function $\|Bx\|_1$, which is expensive to approximate. The PD3O in \cite{yan2018new} avoids using the proximal operator of  $\|Bx\|_1$.
 For a very large scale fused lasso problem, e.g., TV-norm for 2D and 3D data, a method like PD3O should be used. {\bf We emphasize that we compare the proposed method with Davis-Yin splitting and PD3O for the fused lasso problem only for the purpose of validating the numerical performance of the new operator splitting.}  Let $f(x) = \frac{1}{2}\lVert Ax-b \rVert^2_2$, $g(x) = \mu_1\lVert x \rVert_1$, and $h(x) = \mu_2\lVert Bx \rVert_1$.
 It is well known that the proximal operator for $\|x\|_1$ with a step size $\gamma>0$ is the shrinkage operator $\operatorname{Shrinkage}_{\gamma}$. 
 Set $d_1=g, ~d_2=f,~d_3=h$ in Algorithm \ref{FPI_algorithm_modified_convex_case}  or  the scheme \eqref{intro:news} to obtain  
 \begin{subequations}
 \label{lasso-new}
     \begin{eqnarray}
&&x^{k+\frac12}=\operatorname{prox}_{\gamma \mu_2\lVert Bx\rVert_1}(z^k) = \underset{x}{\text{argmin}} (\mu_2 \lVert Bx\rVert_1 + \frac{1}{2\gamma}\lVert x - z^k\rVert^2)\label{lasso-new-TVprox}\\
&&p^{k+1}= \operatorname{Shrinkage}_{\gamma \mu_1}[2x^{k+\frac12}-z^k-\gamma(A^\top(Ax^{k+\frac12}-b))]\\
&&x^{k+1}=(A^\top A+\frac{1}{\gamma}I)^{-1}(A^\top b+\frac{1}{\gamma}(p^{k+1} + \gamma( A^\top( Ax^{k+\frac12})-b)) \label{lasso-new-extra}\\
&&z^{k+1}=z^k+x^{k+1}-x^{k+\frac12}.
\end{eqnarray} 
 \end{subequations}
 The Davis-Yin splitting with $d_1=g, ~d_2=f,~d_3=h$ becomes
\begin{subequations}
    \label{lasso-DY}
\begin{eqnarray}
&&x^{k+\frac12}=\operatorname{prox}_{\gamma \mu_2\lVert Bx\rVert_1}(z^k) = \underset{x\in\mathcal{X}}{\text{argmin}} (\mu_2 \lVert Bx\rVert_1 + \frac{1}{2\gamma}\lVert x - z^k\rVert^2)\\
&&x^{k+1}=\operatorname{Shrinkage}_{\gamma \mu_1}[2x^{k+\frac12}-z^k-\gamma(A^\top(Ax^{k+\frac12}-b))]\\
&&z^{k+1}=z^k+x^{k+1}-x^{k+\frac12}.
\end{eqnarray}
\end{subequations}
\newline \indent Notice that
 the subproblem \eqref{lasso-new-TVprox} is the classical TV norm minimization problem and it is has been known that such a problem can be efficiently approximated by the classical 2-block ADMM method. See \cite{torres2025asymptotic} for a recent result of local linear convergence of ADMM method solving TV norm minimization problems.
In each iteration of 2-block ADMM solving  \eqref{lasso-new-TVprox}, the main computation cost is the inversion $B^\top B+a I$ where $a$ is a positive constant. For 1D TV-norm, $B^\top B+a I$ is a tridiagonal matrix, and we use backslash in MATLAB to efficiently invert such a tridiagonal matrix
in our implementation. 
\newline \indent
For 
the PD3O scheme \eqref{scheme-PD3O}, we set $d_1(x)=\mu_1\|x\|_1, ~d_2=f,~d_3=\mu_2\|x\|_1$, then only matrix multiplication and shrinkage operator are needed thus it is much cheaper and easier to implement.  
\newline \indent We compare the the proposed splitting scheme \eqref{lasso-new} with Davis-Yin method \eqref{lasso-DY} and PD3O \eqref{scheme-PD3O} on the problem \eqref{example-flasso} with $A$ being a random matrix whose elements follow the standard Gaussian distribution, and $b$ is obtained by adding independent and identically distributed Gaussian noise with variance 0.01 onto $Ax$, $\mu_1 = 20$ and $\mu_2 = 200$.  
\newline \indent
We first consider a problem of size $A\in \mathbb R^{100\times 1000}$ shown in Figure \ref{figure_2}. In the 2-block ADMM iteration for the subproblem \eqref{lasso-new-TVprox},  we set
the step size as $1000$  with stopping thresholds as either the difference between two iterations being less than $1\times10^{-8}$ or the total number of iterations exceeding 10000.  
 As we can observe in Figure \ref{figure_2}, the more expensive methods \eqref{lasso-new} and \eqref{lasso-DY} converge faster than PD3O for certain precision regime, and the proposed new method \eqref{lasso-new} is more robust with larger step size $\gamma$ than Davis-Yin splitting \eqref{lasso-DY}. 
PD3O and  Davis-Yin, and ADMM Dual Form splitting fail to converge when the step size $\gamma$ is significantly larger than $\frac{10}{L}$.  
\newline \indent Next we consider a problem of larger size with $A\in \mathbb R^{400\times 20000}$. 
Figure \ref{figure_3} shows the performance. For the PD3O method, we use parameters $\delta=1, \gamma=\frac{1}{L}$, which may not be the optimal ones, but it is already faster in CPU time than the two expensive methods  \eqref{lasso-new} and \eqref{lasso-DY}. For  the subproblem \eqref{lasso-new-TVprox} in this test,
 we set
the step size as $1000$  with stopping thresholds as either the difference between two iterations being less than $1\times10^{-6}$ or the total number of iterations exceeding 500 for  the 2-block ADMM iteration.  
Though the iteration numbers of the proposed new method and the Davis-Yin splitting are almost the same as shown in Figure  \ref{figure_3}(a), the proposed new method is slower in CPU time due to the extra step \eqref{lasso-new-extra}. Nonetheless, the proposed scheme \eqref{lasso-new} is more robust for larger step size as shown in Figure \ref{figure_2}.
\begin{figure}[htbp!]
    \centering
    \subfloat[~$\gamma=\frac{1}{L}$]{\includegraphics[width=0.48\textwidth]{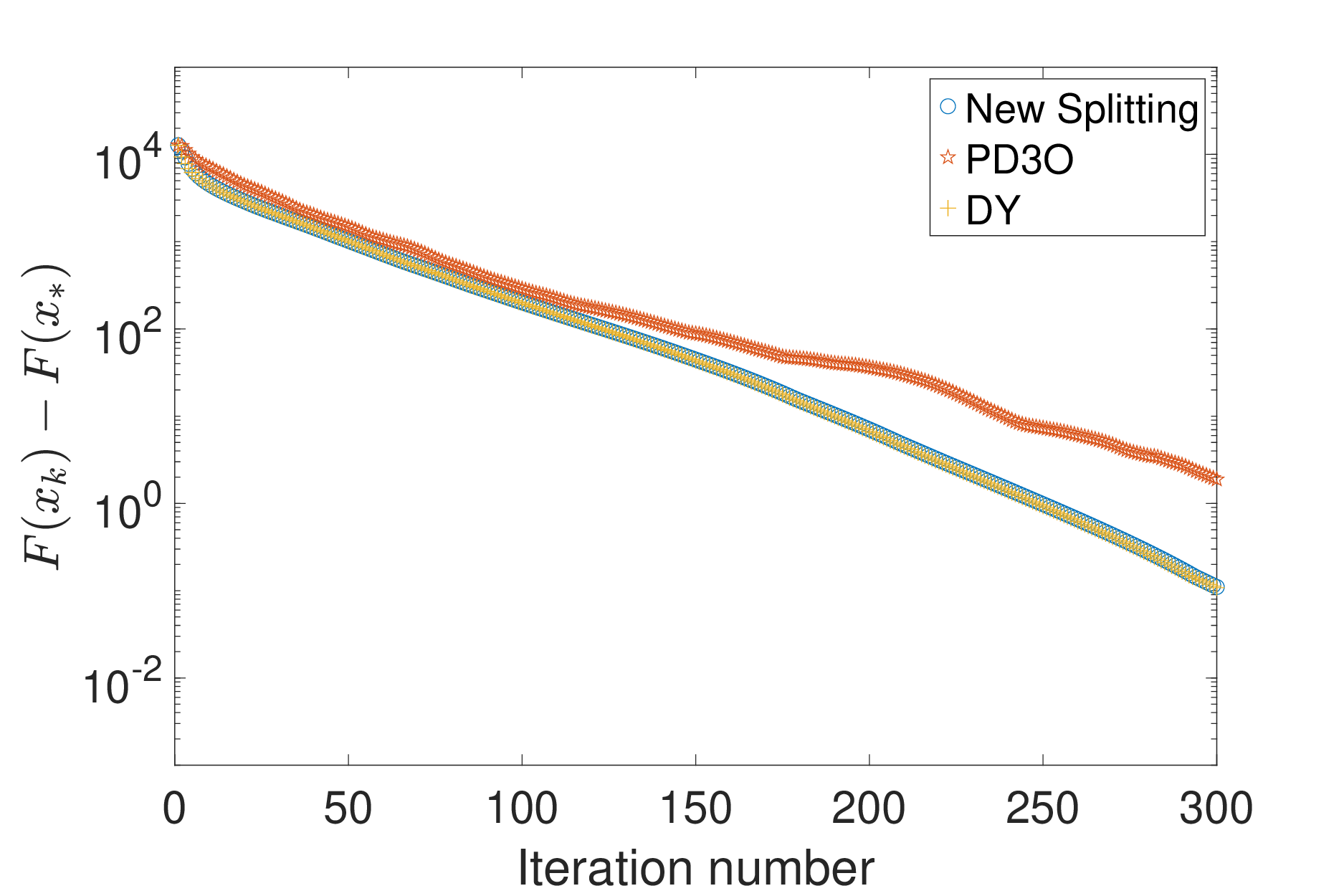}\label{fig2a}} \quad
    \subfloat[~$\gamma=\frac{10}{L}$, \textcolor{black}{Davis-Yin and PD3O fail to converge}]{\includegraphics[width=0.48\textwidth]{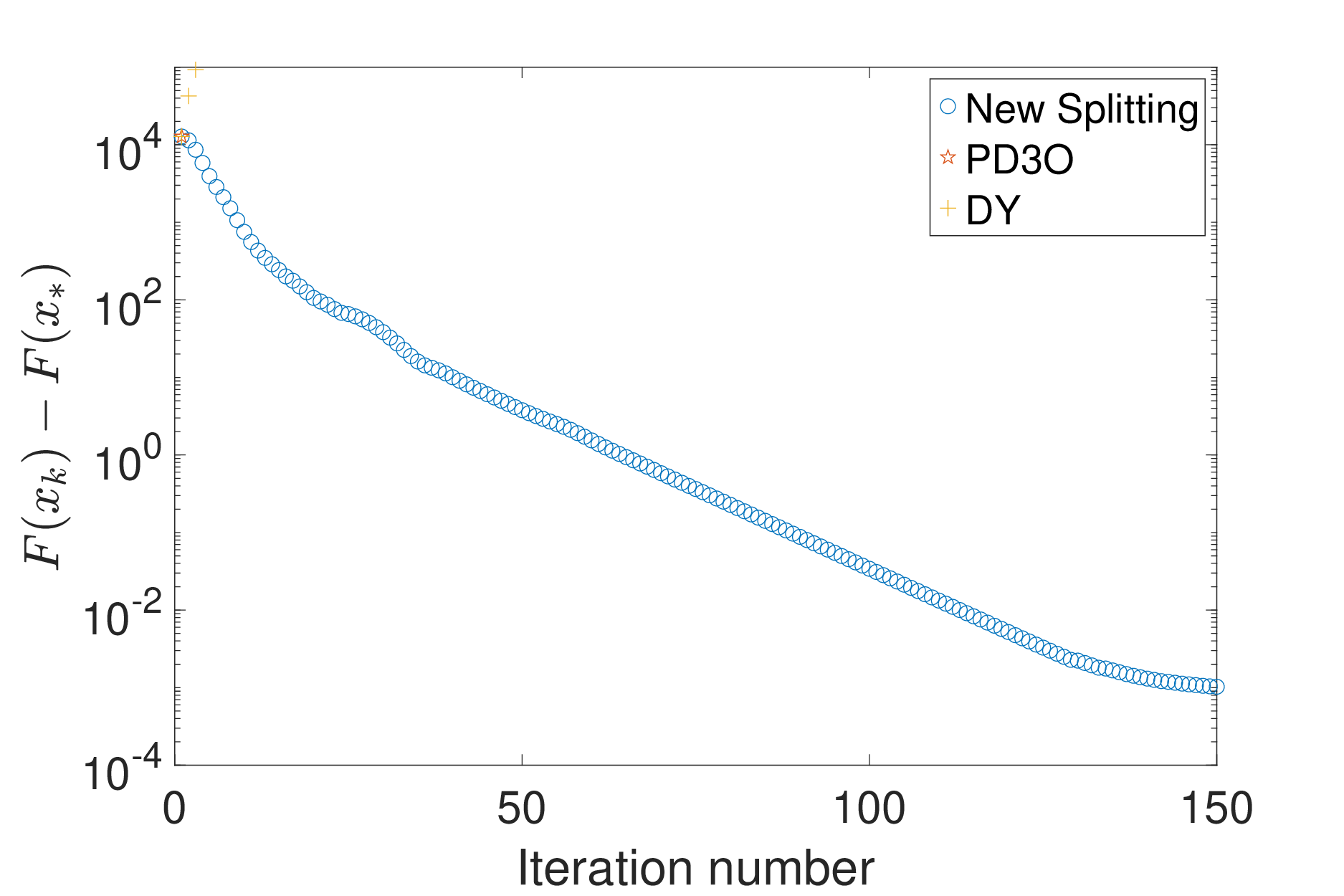}\label{fig2d}}\qquad
    \caption{ \textcolor{black}{Cost function value v.s. iteration numbers for a fused lasso problem with $r = 100$ $n = 1000$, $\mu_1=20$, $\mu=200$. Performance of the three methods using the same step size $\gamma$ and the extra parameter in PD3O is taken as $\delta=\frac{0.9}{\|BB^\top \|}\frac{1}{ \gamma}$. The reference solution $x_*$ is approximated by $40,000$ iterations of PD3O with $\gamma=\frac{1}{L}$.}}
    \label{figure_2}
\end{figure}
\begin{figure}[htbp!]
    \centering
    \subfloat[Cost function value v.s. iteration numbers. The difference between new method and DY is marginal. ]{\includegraphics[width=0.48\textwidth]{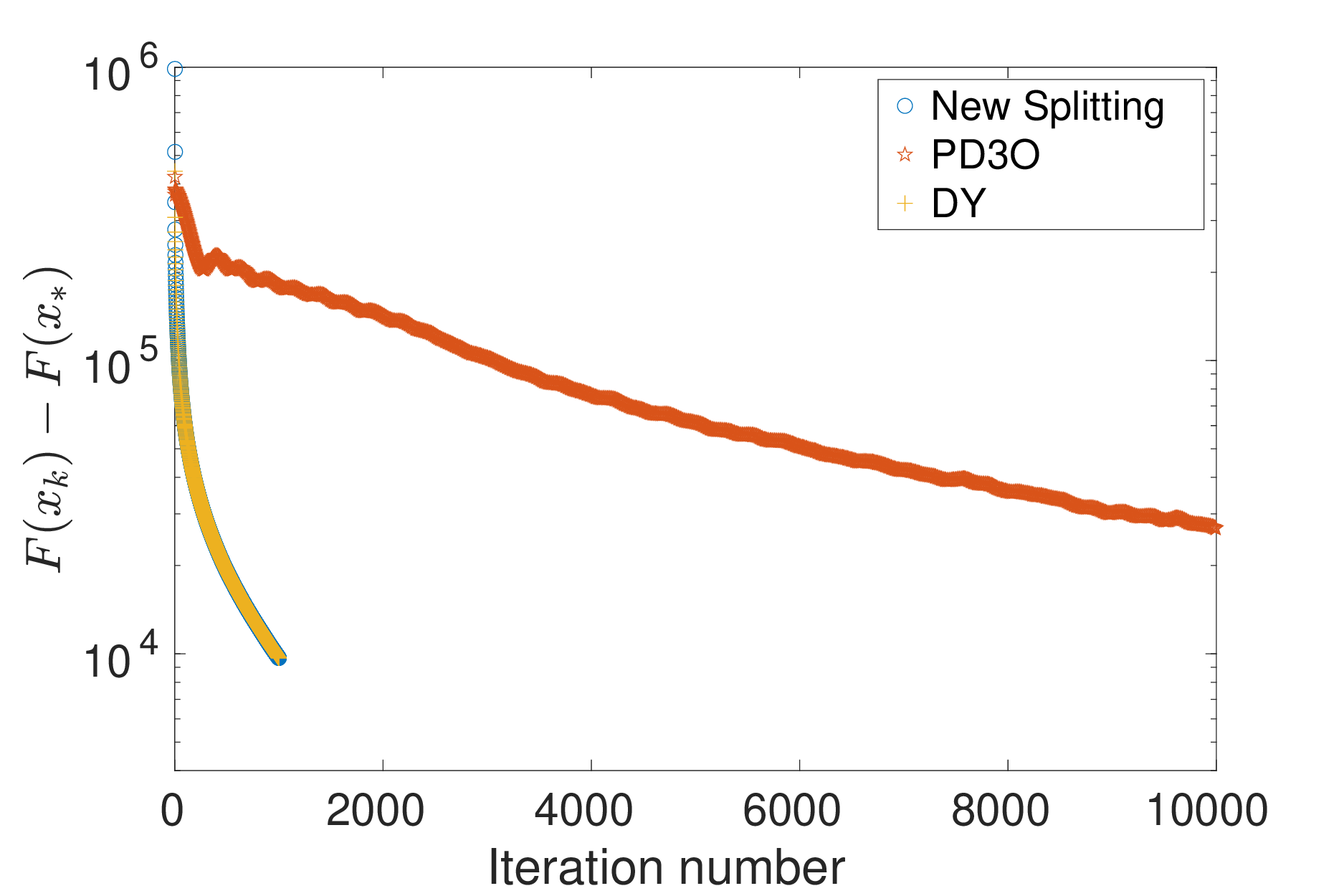}\label{figure5a}} \quad
    \subfloat[Cost function value v.s. CPU time. ]{\includegraphics[width=0.48\textwidth]{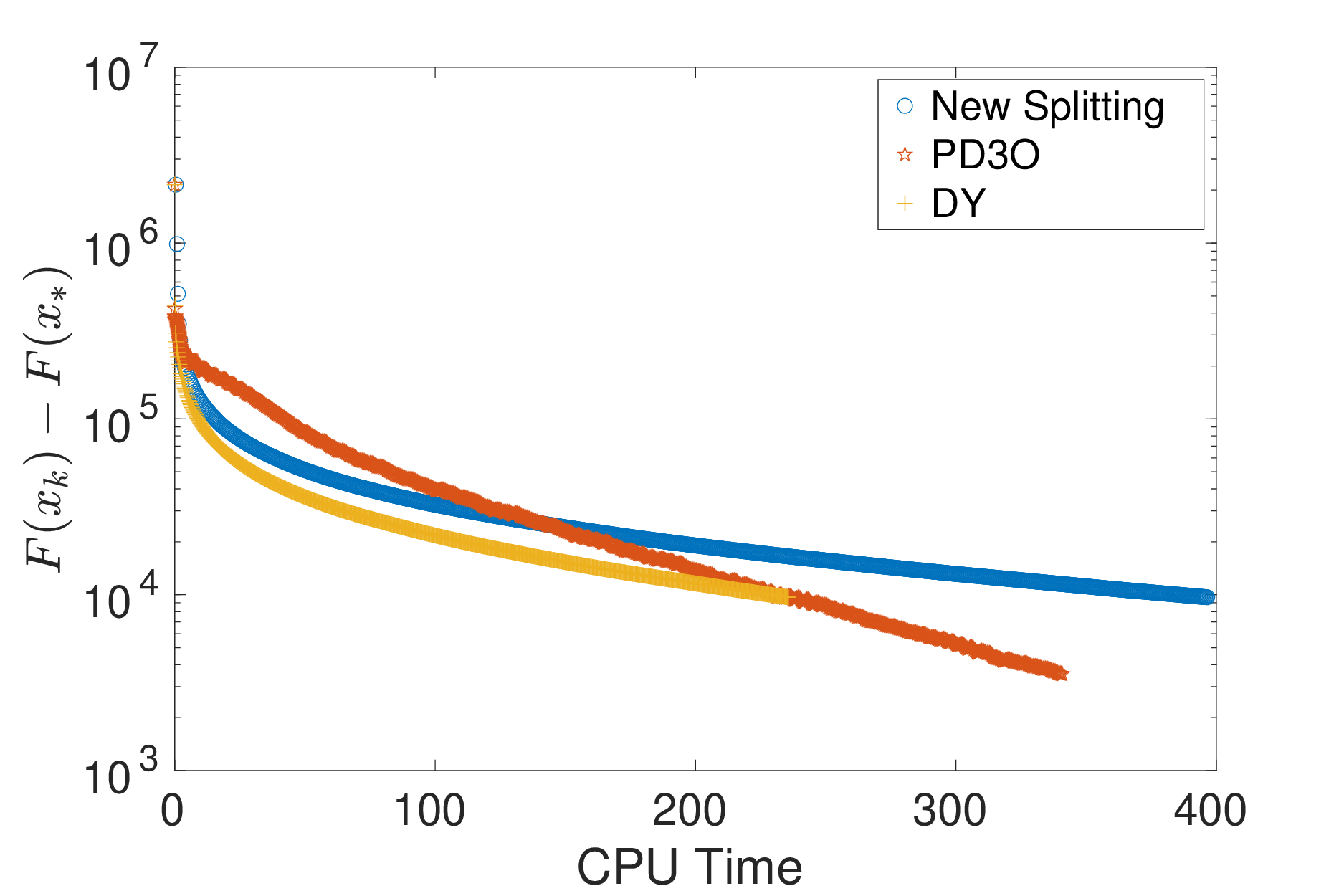}\label{figure5b}} \qquad
    \caption{ \textcolor{black}{A fused Lasso problem with $r = 400$, $n = 20000$, $\mu_1=20$, $\mu_2=200$. Step size is $\gamma=\frac{1}{L}$. The extra parameter in PD3O is taken as $\delta=1$.  The reference solution $x_*$ is approximated by $40,000$ iterations of PD3O.}}
    \label{figure_3}
\end{figure}
\end{example}
}
\section{Concluding Remarks}\label{section_conclusion}
In this paper, we have considered a three-operator splitting scheme for solving monotone inclusion problems, which is an extension of the Douglas-Rachford splitting.
\textcolor{black}{
It is similar to but different from Davis-Yin splitting, and the dual form of the standard three-block ADMM.
 In our numerical tests, it can allow a larger range of step size, compared to the  Davis-Yin Splitting and the dual form of the standard three-block ADMM.
  For solving a composite convex minimization problem, the splitting operator can be proven $1/2$-averaged if two functions have orthogonal domains. Future work includes exploration of nonsmooth problems as those recently studied in \cite{yin2024partial,bai2025inexact}}.\\

\noindent{\bf Acknowledgements} The authors are grateful to the Purdue-India partnership program to make the collaboration possible.   \\

\noindent{\bf Author Contributions} All authors have contributed to the conception and design of the study, analysis, and numerical tests. The manuscript was mainly written and edited by AA and XZ. 
The numerical tests were mainly performed by JL and XZ. All authors read and approved the final manuscript. \\

\noindent{\bf Funding} Anshika acknowledges a Science and Engineering Research Board-Overseas Visiting Doctoral Fellowship, Purdue-India partnership, with file number SB/S9/Z-03/2017-II (2022). X. Zhang is supported by NSF DMS-2208515. D. Ghosh acknowledges the financial support of the research grants MATRICS (MTR/2021/000696) and Core Research Grant (CRG/2022/001347) by the Science and Engineering Research Board, India. \\

\noindent{\bf  Data availability} Not applicable.

\section*{Declarations}

\noindent{\bf Conflict of interest} The authors declare no conflict of interest.\\

\noindent{\bf Ethics approval and consent to participate} Not applicable.\\

\noindent{\bf Consent for publication} Not applicable.\\

\begin{appendices}
\section*{Appendix: Derivation of a 3-block ADMM-like algorithm}
\label{appendix-primal form}
\textcolor{black}{
Since the proposed splitting scheme is different form the classical three-block ADMM, in this section we derive the ADMM-like form of the proposed scheme. For the proposed algorithm
\begin{eqnarray}
&& w^{k+1}=\operatorname{prox}_{\gamma d_3}(z^k)\label{relation_1}\\
&& p^{k+1}= \operatorname{prox}_{\gamma d_1} (2w^{k+1}-z^k-\gamma\partial d_2(v^{k}))\label{relation_2}\\
&& v^{k+1}=\operatorname{prox}_{\gamma d_2} (p^{k+1}+\gamma\partial d_2(v^{k})) \label{relation_3}\\
&& z^{k+1}=z^k+v^{k+1}-w^{k+1}, \label{relation_4}
\end{eqnarray}
we start with relation \eqref{relation_1}
\begin{eqnarray*}
  w^{k+1}=\operatorname{prox}_{\gamma d_3}(z^k)
&\iff&(I+\gamma \partial d_3)(w^{k+1})=z^k \nonumber \\
&\iff& w^{k+1}=z^k-\gamma \partial d_3(w^{k+1}) \nonumber\\
&\iff& w^{k+1}= z^k-\gamma (A_3\partial f_3^*(A_3^\top w^{k+1})-b) \nonumber\\
&\iff& A_3^Tw^{k+1}= A_3^T(z^k-\gamma (A_3\partial f_3^*(A_3^\top w^{k+1})-b)) \nonumber \\
&\iff& \partial f_3(\bar w_{k+1})= A_3^T(z^k-\gamma A_3(\bar w_{k+1})+\gamma b), \label{relation_9}
\end{eqnarray*}
 where   $\bar w^{k+1} = \partial f_3^*(A_3^Tw^{k+1})$, and $f^*$ denotes the convex conjugate of the proper, closed, and convex function $f$. By duality, they satisfy $A_3^Tw^{k+1} = \partial f_3(\bar w^{k+1})$. From relation \eqref{relation_2}, we have
\begin{eqnarray*}
&&p^{k+1}=\operatorname{prox}_{\gamma d_1}(2w^{k+1}-z^k-\gamma\partial d_2(w^{k+1})) \nonumber\\
&\iff& (I+\gamma\partial d_1)(p^{k+1})=2w^{k+1}-z^k-\gamma\partial d_2(w^{k+1}) \nonumber\\
&\iff& p^{k+1} = 2w^{k+1}-z^k-\gamma\partial d_2(w^{k+1}) -\gamma\partial d_1(p^{k+1})\\
&\iff& p^{k+1} = 2w^{k+1}-z^k-\gamma A_2\partial f_2^*(A_2^Tw^{k+1}) -\gamma A_1\partial f_1^*(A_1^Tp^{k+1})\\
&\iff& A_1^Tp^{k+1} = A_1^T(2w^{k+1}-z^k-\gamma A_2\partial f_2^*(A_2^Tw^{k+1}) -\gamma A_1\partial f_1^*(A_1^Tp^{k+1}))\\
&\iff& \partial f_1(\bar p_{k+1}) = A_1^T(2w^{k+1}-z^k-\gamma A_2\partial f_2^*(A_2^Tw^{k+1}) -\gamma A_1 \bar p^{k+1})\\
&\iff& \partial f_1(\bar p_{k+1}) = A_1^T(2(z^k-\gamma A_3 \bar w_{k+1} + \gamma b))-z^k-\gamma A_2\partial f_2^*(A_2^Tw^{k+1}) -\gamma A_1 \bar p^{k+1})\\
&\iff& \partial f_1(\bar p_{k+1}) = A_1^T(z^k-2\gamma A_3 \bar w_{k+1} -\gamma A_2 \tilde w^{k+1} -\gamma A_1 \bar p^{k+1}+2\gamma b),
\label{relation_11}
\end{eqnarray*}
where $\bar p^{k+1} = \partial f_1^*(A_1^Tp^{k+1})$ and $\tilde w^{k+1} = \partial f_2^*(A_2^Tw^{k+1})$. From relation \eqref{relation_3}, we have
\begin{eqnarray*}
&&v^{k+1}=\operatorname{prox}_{\gamma d_2}(p^{k+1}+\gamma\partial d_2(w^{k}))\nonumber \\
&\iff& (I+\gamma\partial d_2)(v^{k+1})=p^{k+1}+\gamma\partial d_2(w^{k}) \nonumber \\
&\iff& v^{k+1}=p^{k+1}+\gamma\partial d_2(w^{k}) - \gamma\partial d_2(v^{k+1}) \nonumber \\
&\iff& v^{k+1}=p^{k+1}+\gamma A_2\partial f_2^*(A_2^Tw^{k+1})-\gamma A_2\partial d_2(A_2^Tv^{k+1})\label{relation_12} \\
&\iff& A_2^Tv^{k+1}=A_2^T(p^{k+1}+\gamma A_2\partial f_2^*(A_2^Tw^{k+1})-\gamma A_2\partial d_2(A_2^Tv^{k+1}))\\
&\iff& \partial f_2(\bar v^{k+1})=A_2^T(p^{k+1}+\gamma A_2\partial f_2^*(A_2^Tw^{k+1})-\gamma A_2\bar v^{k+1})\\
&\iff& \partial f_2(\bar v^{k+1})=A_2^T(2w^{k+1}-z^k -\gamma A_1\bar p^{k+1}-\gamma A_2\bar v^{k+1})\\
&\iff& \partial f_2(\bar v^{k+1})=A_2^T(2(z^k-\gamma A_3(\bar w_{k+1})+\gamma b)-z^k -\gamma A_1\bar p^{k+1}-\gamma A_2\bar v^{k+1})\\
&\iff& \partial f_2(\bar v^{k+1})=A_2^T(z^k-\gamma A_1\bar p^{k+1}-\gamma A_2\bar v^{k+1}-2\gamma A_3(\bar w_{k+1})+2\gamma b),
\end{eqnarray*}
where  $\bar v^{k+1} = \partial f_2^*(A_2^Tv^{k+1})$. Finally, from \eqref{relation_4} we have
\begin{eqnarray*}
&& z^{k+1}=z^k+v^{k+1}-w^{k+1}\\
&&=z^k+p^{k+1}+\gamma A_2\tilde w^{k+1}-\gamma A_2\bar v^{k+1}-z^{k}+\gamma A_3\bar w^{k+1}-\gamma b\\
&&=z^k-\gamma(A_1\bar p^{k+1}+A_2\bar v^{k+1}+A_3\bar w^{k+1}-b)
\end{eqnarray*}
To summarize, after renaming the variables, we have 
\begin{eqnarray*}
&&\partial f_3(\bar w_{k+1})= A_3^T(z^k-\gamma A_3(\bar w_{k+1})+\gamma b)\\ \nonumber
&\iff& x_3^{k+1} = \bar w^{k+1} = \underset{x_3\in\mathcal{X}_3}{\operatorname{argmin}}\{f_3(x_3)+\tfrac{\gamma}{2}\lVert A_3x_3-b-\frac{z^k}{\gamma} \rVert^2\}\\
&&\partial f_1(\bar p_{k+1}) = A_1^T(z^k-2\gamma A_3 \bar w_{k+1} -\gamma A_2 \tilde w^{k+1} -\gamma A_1 \bar p^{k+1}+2\gamma b)\\ \nonumber
&\iff& x_1^{k+1} = \bar p^{k+1} = \underset{x_1\in\mathcal{X}_1}{\operatorname{argmin}}\{f_1(x_1)+\tfrac{\gamma}{2}\lVert A_1x_1 + A_2\tilde w^{k+1} -2A_3x_3^{k+1}-2b-\frac{z^k}{\gamma} \rVert^2\}\\
&&\partial f_2(\bar v^{k+1})=A_2^T(z^k-\gamma A_1\bar p^{k+1}-\gamma A_2\bar v^{k+1}-2\gamma A_3(\bar w_{k+1})+2\gamma b)\\ \nonumber
&\iff& x_2^{k+1} = \bar p^{k+1} = \underset{x_2\in\mathcal{X}_2}{\operatorname{argmin}}\{f_2(x_2)+\tfrac{\gamma}{2}\lVert A_1x_1^{k+1} + A_2x_2 -2A_3x_3^{k+1}-2b-\frac{z^k}{\gamma} \rVert^2\}\\
&&z^{k+1} = z^k-\gamma(A_1\bar p^{k+1}+A_2\bar v^{k+1}+A_3\bar w^{k+1}-b)\\ \nonumber
&\iff& z^{k+1} = z^k-\gamma(A_1x_1^{k+1}+A_2x_2^{k+1}+A_3x_3^{k+1}-b).
\end{eqnarray*} So the derived algorithm is
\begin{eqnarray}\label{primal form proposed splitting}
\begin{rcases}
& x_3^{k+1} = \underset{x_3\in\mathcal{X}_3}{\operatorname{argmin}}\{f_3(x_3)+\tfrac{\gamma}{2}\lVert A_3x_3-b-\frac{z^k}{\gamma} \rVert^2\}\\
& x_1^{k+1} = \underset{x_1\in\mathcal{X}_1}{\operatorname{argmin}}\{f_1(x_1)+\tfrac{\gamma}{2}\lVert A_1x_1 + A_2\tilde w^{k+1} -2A_3x_3^{k+1}-2b-\frac{z^k}{\gamma} \rVert^2\}\\
& x_2^{k+1} = \underset{x_2\in\mathcal{X}_2}{\operatorname{argmin}}\{f_2(x_2)+\tfrac{\gamma}{2}\lVert A_1x_1^{k+1} + A_2x_2 -2A_3x_3^{k+1}-2b-\frac{z^k}{\gamma} \rVert^2\}\\
&z^{k+1} = z^k-\gamma(A_1x_1^{k+1}+A_2x_2^{k+1}+A_3x_3^{k+1}-b)
\end{rcases},
\end{eqnarray}
where $\tilde w^{k+1} = \partial f_2^*(A_2^Tw^{k+1}) = \partial f_2^*(A_2^T(z^k-\gamma (A_3x_3^{k+1}-b))).$
}

\begin{figure}[htbp!]
    \centering
\includegraphics[width=0.48\textwidth]{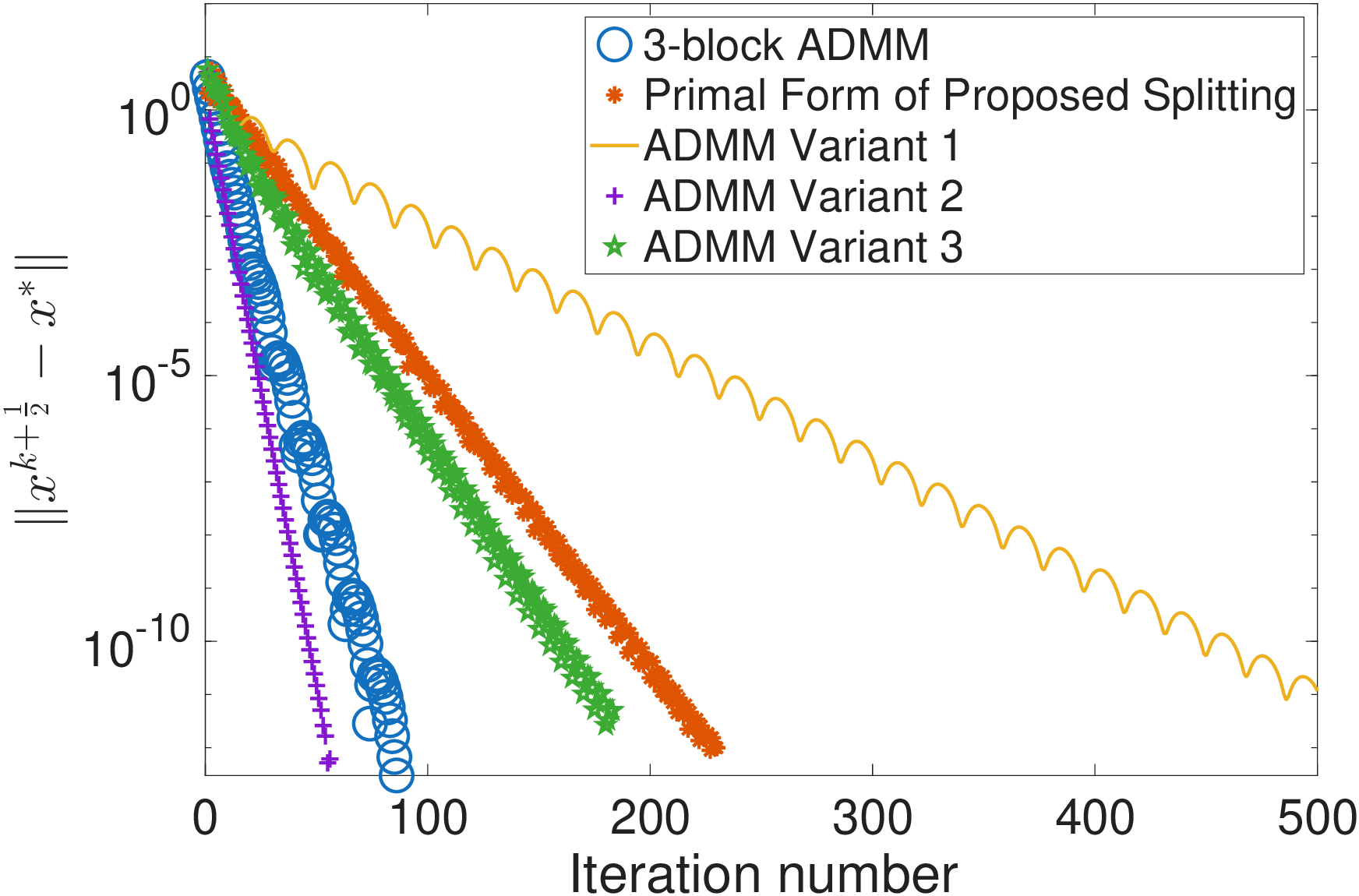} 
    \caption{\textcolor{black}{Comparison of the error $\|x^{k+\frac12}-x^*\|$ between the derived algorithm   \eqref{primal form proposed splitting} (the corresponding legend is {\it Primal Form of Proposed Splitting}) and other 3-block ADMM variants. The convergence rate of derived algorithm is between the ADMM variants.}}
    \label{figure_5}
\end{figure}

\textcolor{black}{
Since there are quite a few variants of 3-block ADMM methods, we consider a numerical comparison of some existing variants to the derived 3-block ADMM-like algorithm \eqref{primal form proposed splitting},  as a  demonstration of its practical efficiency.
We consider an example transformed from the problem \eqref{example-limiter}.
Let $f_2(x_2)=\tfrac{\alpha}{2}\lVert x_2- u\rVert^2_2, f_1(x_1)= {i}_{\Lambda_1}(x_1),$ and $f_3(x_3)=i_{\Lambda_2}(x_3)$, where the two indicator functions are for the same two sets as in  \eqref{example-limiter}. The derived scheme applies here with the following constraint:
\begin{eqnarray*}
A_1x_1 + A_2x_2 + A_3x_3 = d,
\end{eqnarray*} where $A_1 = \begin{bmatrix}I \\ 0 \end{bmatrix}\in\mathbb{R}^{2n\times n}$, $A_2 = \begin{bmatrix}-I \\ I \end{bmatrix}\in\mathbb{R}^{2n\times n}$, $A_3 = \begin{bmatrix}0 \\ -I \end{bmatrix}\in\mathbb{R}^{2n\times n}$, and $d =  \mathbf{0}_{2n}\in\mathbb{R}^{2n}$. For this problem, we obtain from \eqref{primal form proposed splitting}:
\begin{eqnarray*}
&&x_3^{k+1}= -z_2^{k}/\gamma + (\sum_{i=1}^n z_2^{k}(i)/\gamma + b)/n \cdot\mathbf{1}_n;\\
&&\tilde w^{k+1} = (-z_1^{k} + z_2^{k} + \gamma x_3^{k+1})/\alpha + u;\\
&&x_1^{k+1}=\min(\max(\tilde w^{k+1} + z_1/\gamma), m),M)\\
&&x_2^{k+1}=(\alpha u + \gamma x_1^{k+1} - 2\gamma x_3^{k+1} - z_1^{k} + z_2^{k})/(\alpha + 2\gamma)\\
&&z_1^{k+1}=z_1^k-\gamma(x_1^{k+1}-x_2^{k+1})\\
&&z_2^{k+1}=z_2^k-\gamma(x_2^{k+1}-x_3^{k+1}).\\
\end{eqnarray*}
The direct extension of 3-block ADMM \eqref{dual-ADMM} applied to this problem is:
\begin{eqnarray*}
&&x_1^{k+1}= \min(\max(\tilde x_2^{k+1} + w_1/\gamma), m),M);\\
&&x_2^{k+1}=(\alpha u + \gamma x_1^{k+1} + \gamma x_3^{k+1} - w_1^{k} + w_2^{k})/(\alpha + 2\gamma)\\
&&x_3^{k+1}=x_2 - w_2 / \gamma - (\sum_{i=1}^n (x_2(i) - w_2^{k}(i)/\gamma) + b)/n \cdot\mathbf{1}_n\\
&&w_1^{k+1}=w_1^k-\gamma(x_1^{k+1}-x_2^{k+1})\\
&&w_2^{k+1}=w_2^k-\gamma(x_2^{k+1}-x_3^{k+1}),\\
\end{eqnarray*}
where $\mathbf{1}_n$ is the n-dimensional vector with all components equal to $1$. We compare the derived algorithm with the direct extension of 3-block extension of ADMM with variant 1 \cite{he2018aclass} that uses prediction correction steps with parameter $\tau = \frac{1}{2}$ in the correction step, variant 2 in\cite{davis2017three} that only changes the $x_1^{k+1}$ update, and variant 3 in\cite{He2012GaussianBack} that uses prediction-correction steps with Gaussian back substitution with $\alpha = \frac{1}{2}$ in the correction step.
For simplicity we use step size $\gamma = 1$ for algorithms.
All methods are tested on the problem \eqref{example-limiter} with the same parameters $\alpha=1, n=100, m=-1, M=1$, and $b=Au$, where $u$ is constructed by perturbing a sine profile by random noise: $$u_i=\sin(2\pi\frac{i}{n})+0.8*\mathcal N(0,1).$$ The comparison is shown in Figure \ref{figure_5}.
}
\end{appendices}

\bibliographystyle{spmpsci}
\bibliography{mybib}

\end{document}